\documentclass[letter,preprint,11pt]{elsarticle}
\biboptions{number}

\usepackage[dvips]{color}
\usepackage{amsmath,amssymb,amsthm,latexsym,mathrsfs,euscript,verbatim}
\usepackage{graphicx,xtab}
\usepackage{hyperref}
\usepackage{pstricks,amscd}

\newtheoremstyle{thm}
{10pt}% this gives how much space is left above theorem
{6pt}% space below theorem
{\itshape}% font used (italics for theorems)
{}% Indent amount (empty = no indent, \parindent = paragraph indent)
{\bf}% Tital font
{.}% Punctuation after thm head
{.5em}% Space after thm head
{} % This one should just be left empty

\theoremstyle{thm}

\newtheorem{thm}{Theorem}[section]
\newtheorem{prop}[thm]{Proposition}
\newtheorem{lemma}[thm]{Lemma}
\newtheorem{cor}[thm]{Corollary}

\theoremstyle{definition}
\newtheorem{defn}[thm]{Definition}

\newcommand{\nc}{\newcommand}
\newcommand{\rnc}{\renewcommand}
\renewcommand{\l}{\lambda}
\renewcommand{\d}{\delta}
\rnc{\k}{\kappa}

\newcommand{\e}{\epsilon}
\newcommand{\noi}{\noindent}
\nc{\bmw}{\mathscr{C}}
\nc{\ckt}[1]{\mathbb{KT}_{#1}^k}
\nc{\akt}[1]{\widehat{\mathbb{KT}}_{#1}}
\nc{\cl}[1]{\mathrm{cl}_{#1}}
\newcommand{\bnk}{\mathscr{B}_n^k}
\newcommand{\bmk}{\mathscr{B}_{n-1}^k}
\newcommand{\bt}{\mathscr{B}_2^k}
\rnc{\b}[1]{\,\overline{\!#1}}
\nc{\blk}{\mathscr{B}_l^k}
\nc{\bl}[1]{\widetilde{\mathscr{B}}_{#1}^k}
\nc{\ak}[1]{\mathfrak{h}_{#1,k}}
\newcommand{\At}[1]{\widetilde{\mathfrak{h}}_{#1,k}}
\newcommand{\ix}{X_i^{-1}}
\newcommand{\iy}{Y^{-1}}
\nc{\kum}[1]{\sum_{#1=0}^{k-1}}
\nc{\kkum}[1]{\sum_{#1=0}^{k}}
\nc{\pum}[1]{\sum_{#1=1}^p}
\nc{\ppum}[1]{\sum_{#1=0}^{p-1}}
\nc{\tw}[1]{0 \leq #1 \leq k-1}
\nc{\tf}{\text{,\quad for all }}
\nc{\iif}{\text{,\quad if }}
\nc{\la}[1]{\left\langle #1 \right\rangle}
\nc{\rel}[1]{\stackrel{\tiny{(\ref{#1})}}{=}}
\newcommand{\stack}[1]{\stackrel{\tiny{#1}}{=}}
\nc{\sstack}[1]{\stackrel{\tiny{#1}}{\subseteq}}
\nc{\f}{&\hphantom{=}\,\,}
\nc{\p}{\stackrel{\hphantom{(1)}}{=}}
\nc{\h}{\stackrel{\hphantom{(19)}}{=}}
\newcommand{\y}[2]{Y^{\prime \,#2}_{#1}}

\rnc{\a}[2]{\alpha_{#1}^{#2}}
\rnc{\o}{\ldots}
\nc{\floor}[1]{\left\lfloor #1 \right\rfloor}
\nc{\ceil}[1]{\left\lceil \frac{#1}{2} \right\rceil}
\rnc{\ll}[1]{\left( #1 \right)}
\nc{\s}{\subseteq}
\nc{\Z}{\mathbb{Z}}
\nc{\dg}{\psi} %map from \bnk to CKT
\nc{\wh}[1]{\widehat{#1}}

\nc{\M}{\mathbb{M}}
\nc{\X}{\mathfrak{W}}

\def\be{\begin{align*}}
\def\ee{\end{align*}}
\rnc{\=}{&=&}
\rnc{\-}{&=}
\newcommand{\beq}{\begin{align}}
\newcommand{\eeq}{\end{align}}
\rnc{\thefootnote}{\fnsymbol{footnote}}

\journal{\,}

\begin{document}
\begin{frontmatter}

\title{On the freeness of the cyclotomic BMW algebras: \\ admissibility and an isomorphism with the cyclotomic Kauffman tangle algebras}
\makeatletter
\def\authorbugfix[#1]#2{\g@addto@macro\elsauthors{%
   \def\baselinestretch{1}%
   \authorsep#2\unskip\textsuperscript{%#1%
     \@for\@@affmark:=#1\do{%
      \edef\affnum{\@ifundefined{X@\@@affmark}{1}{\Ref{\@@affmark}}}%
                  \unskip\comma\affnum\let\comma,}%
     %\ifx\@cormark\@empty\else\unskip,\@cormark\fi
     \ifx\@fnmark\@empty\else\unskip\@fnmark\fi}%
   \def\authorsep{\unskip,\space}%
   \global\let\@cormark\@empty
   \global\let\@fnmark\@empty}%
   \@eadauthor={#2}
}
\makeatother

\authorbugfix[S1]{Stewart Wilcox}
%\author[S1]{Stewart Wilcox}
\ead{swilcox@fas.harvard.edu}

\author[S2]{Shona Yu\corref{fn1}}
\ead{s.h.yu@tue.nl}
\cortext[fn1]{Corresponding author}

\address[S1]{Stewart Wilcox, Department of Mathematics, Harvard University, 1 Oxford Street, \\ Cambridge MA 02318, USA.}
\address[S2]{Shona Yu, Den Dolech 2, Faculteit Wiskunde en Informatica, Technische Universiteit Eindhoven, \\ 5600 MB Eindhoven, The Netherlands.}

\begin{abstract}
The cyclotomic Birman-Murakami-Wenzl (or BMW) algebras $\bnk$, introduced by R.~H\"aring-Oldenburg,
are a generalisation of the BMW algebras associated with the cyclotomic Hecke algebras of type $G(k,1,n)$ (also known as Ariki-Koike algebras) and type $B$ knot theory.
In this paper, we prove the algebra is free and of rank $k^n (2n-1)!!$ over ground rings with parameters satisfying so-called ``admissibility conditions". These conditions are necessary in order for these results to hold and arise from the representation theory of $\bt$, which is analysed by the authors in a previous paper. Furthermore, we obtain a geometric realisation of $\bnk$ as a cyclotomic version of the Kauffman tangle algebra, in terms of affine $n$-tangles in the solid torus, and produce explicit bases that may be described both algebraically and diagrammatically.
\end{abstract}

\begin{keyword}
cyclotomic BMW algebras; cyclotomic Hecke algebras; Ariki-Koike algebras; Brauer algebras; affine tangles; Kauffman link invariant; admissibility.

\MSC 16G99, 20F36, 81R05, 57M25
\end{keyword}

\end{frontmatter}

%% \linenumbers

%% main text
\section{Introduction}

The cyclotomic BMW algebras were introduced in \cite{HO01} by H\"aring-Oldenburg  as a generalisation of the BMW algebras associated with the cyclotomic Hecke algebras of type $G(k,1,n)$ (also known as Ariki-Koike algebras) and type $B$ knot theory involving affine tangles.

The motivation behind the definition of the BMW algebras may be traced back to an important problem in knot theory; namely, that of classifying knots (and links) up to isotopy.
The BMW algebras (conceived independently by Murakami \cite{M87} and Birman and Wenzl \cite{BW89}) are algebraically defined by generators and relations modelled on certain tangle diagrams appearing in the skein relation for the Kauffman link invariant of \cite{K90}, an invariant of regular isotopy for links in $S^3$: 
\vspace{0.2cm}
  \begin{center}
    \includegraphics{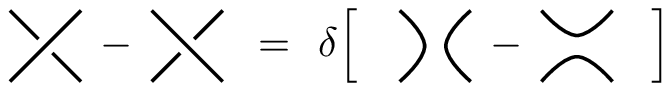}
  \end{center} 
%\vspace{0.2cm}
  
\begin{defn} \label{defn:bmw}
Fix a natural number $n$. Let $R$ be a unital commutative ring 
containing an element $A_0$ and units $q$ and $\l$ such that $\l - \l^{-1} = \d(1-A_0)$ holds, where $\d := q - q^{-1}$. 
The \emph{\textbf{BMW algebra}} 
$\bmw_n := \bmw_n(q,\l,A_0)$ is defined to be the unital associative $R$-algebra 
generated by $X_1^{\pm 1}, \ldots, X_{n-1}^{\pm 1}$ and 
$e_1, \ldots, e_{n-1}$ subject to the following relations, which hold for all possible values of $i$ unless otherwise stated.
\[
\begin{array}{rcll}
\hspace{2cm} X_i - X_i^{-1} &=& \d(1-e_i)& \\
X_iX_j &=& X_jX_i &\qquad \text{for } |i-j| \geq 2 \\
X_iX_{i+1}X_i &=& X_{i+1}X_iX_{i+1}& \\
X_ie_j &=& e_jX_i &\qquad \text{for } |i-j| \geq 2 \\
e_ie_j &=& e_je_i &\qquad \text{for } |i-j| \geq 2 \\
X_ie_i &=& e_iX_i \,\,\,=\,\,\, \l e_i& \\
X_i X_j e_i &=& e_j e_i \,\,\,=\,\,\, e_j X_i X_j &\qquad \text{for } |i-j| = 1 \\
e_i e_{i\pm 1} e_i &=& e_i& \\
e_i^2 &=& A_0 e_i.& 
\end{array}
\]
\end{defn}

In particular, the defining relations were originally inspired by the diagrammatic relations satisfied by these tangle diagrams and the relations $X_i - X_i^{-1} = \d(1-e_i)$ seen in Definition~\ref{defn:bmw} above reflects the Kauffman skein relation. Furthermore, the Kauffman link polynomial may be recovered from a nondegenerate Markov trace function on the BMW algebras, in a way analogous to the relationship between the Jones polynomial and the Temperley-Lieb algebras.

Naturally, one would expect the BMW algebras to have a geometric realisation in terms of tangles. Indeed, under the maps illustrated below, Morton and Wasserman \cite{MW89} proved the BMW algebra $\bmw_n$ is isomorphic to the Kauffman tangle algebra $\mathbb{KT}_{n}$, an algebra of (regular isotopy equivalence classes of) tangles on $n$ strands in the disc cross the interval (that is, a solid cylinder) modulo the Kauffman skein relation (see Kauffman \cite{K90} and Morton and Traczyk \cite{MT90}). As a result, they also show the algebra $\bmw_n$ is free of rank $(2n-1)!! = (2n-1) \cdot (2n-3) \cdots 3 \cdot 1$.
\vspace{0.2cm}
  \begin{center}  
 \includegraphics{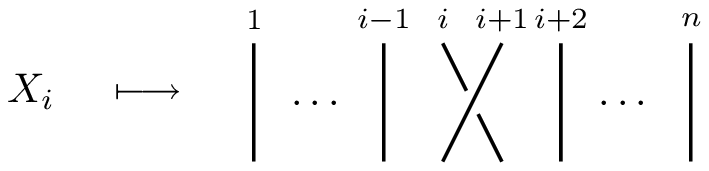} 
    \vskip 0.7cm
    \hspace*{-0.07 cm} 
    \includegraphics{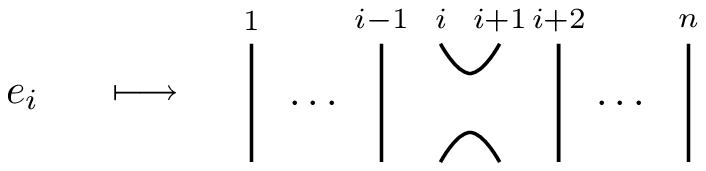}
\end{center}
\vspace{0.2cm} 

The BMW algebras are closely connected with the Artin braid groups of type~$A$, Iwahori-Hecke algebras of type $A$, and with many diagram algebras, such as the Brauer and Temperley-Lieb algebras. In fact, they may be construed as deformations of the Brauer algebras obtained by replacing the symmetric group algebras with the corresponding Iwahori-Hecke algebras. These various algebras also feature prominently in the theory of quantum groups, subfactors, statistical mechanics, and topological quantum field theory.

In view of these relationships between the BMW algebras and several objects of ``type $A$'', several authors have since naturally generalised the BMW algebras for other types of Artin groups.
Motivated by knot theory associated with the Artin braid group of type $B$, H\"aring-Oldenburg  introduced the ``cyclotomic BMW algebras'' in \cite{HO01}. They are so named because the cyclotomic Hecke algebras of type $G(k,1,n)$ from \cite{AK94,BM93}, which are also known as Ariki-Koike algebras, arise as quotients of the cyclotomic BMW algebras~$\bnk$ in the same way the Iwahori-Hecke algebras arise as quotients of the BMW algebras. They are obtained from the original BMW algebras by adding an extra generator $Y$ satisfying a polynomial relation of finite order $k$ and imposing several further relations modelled on type $B$ knot theory. For example, $Y$ satisfies the Artin braid relations of type $B$ with the generators $X_1$, \o, $X_{n-1}$ of the ordinary BMW algebra. 
When this $k^\mathrm{th}$ order relation on the generator $Y$ is omitted, one obtains the infinite dimensional affine BMW algebras, studied by Goodman and Hauschild in \cite{GH06}. This extra affine generator $Y$ may be visualised as the affine braid of type $B$ illustrated below.

\begin{center}
  \includegraphics{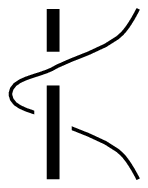}
\end{center}

Given what has already been established for the BMW algebras, it is then conceivable that the cyclotomic and affine BMW algebras be isomorphic to appropriate analogues of the
Kauffman tangle algebras. Indeed, by utilising the results and techniques of Morton and Wasserman~\cite{MW89} for the ordinary BMW algebras, this was shown to be the case for the affine version, by Goodman and Hauschild in \cite{GH06}. The topological realisation of the affine BMW algebra is as an algebra of (regular isotopy equivalence classes of) affine tangles on $n$ strands in the annulus cross the interval (that is, the solid torus) modulo Kauffman skein relations. 

In this paper, we prove the cyclotomic BMW algebras $\bnk(R)$ are $R$-free of rank $k^n (2n-1)!!$ and show they have a topological realisation as a certain cyclotomic analogue of the Kauffman tangle algebra in terms of affine $n$-tangles (see Definition \ref{defn:cycloKT}). Furthermore, we obtain bases that may be explicitly described both algebraically and diagrammatically in terms of affine tangles. One may visualise the basis given in Theorem \ref{thm:span} and Corollary \ref{cor:basis} as a kind of ``inflation" of bases of smaller Ariki-Koike algebras by `dangles', as seen in Xi \cite{X00}, with powers of Jucy-Murphy type elements $Y_i'$ attached. This is illustrated in Figure \ref{fig:inflatedbasis}.

Unlike the BMW and Ariki-Koike algebras, one needs to impose extra so-called ``admissibility conditions'' (see Definition \ref{defn:adm}) on the parameters of the ground ring in order for these results to hold. This is due to potential torsion on elements associated with certain tangles on two strands, caused by the polynomial relation of order $k$ imposed on $Y$. It turns out that the representation theory of~$\mathscr{B}_2^k$, analysed in detail by the authors in \cite{WY06}, is crucial in determining these conditions precisely. A particular result in \cite{WY06} shows that admissibility ensures freeness of the algebra $\bt(R)$ over~$R$. These results are stated but incompletely proved in H\"aring-Oldenburg \cite{HO01}. Moreover, it turns out that admissibility as defined in this paper (not \cite{WY06}) is \emph{necessary and sufficient} for freeness results for \emph{general}~$n$.

The results presented here are proved in the Ph.D. thesis \cite{Y07}, completed at the University of Sydney in 2007 by the second author, in which these bases are shown to lead to a cellular basis, in the sense of Graham and Lehrer \cite{GL96} (see also \cite{WY09b}). When $k = 1$, all results specialise to those previously established for the BMW algebras by Morton and Wasserman \cite{MW89}, Enyang~\cite{E04} and Xi \cite{X00}.

Since the submission of this thesis, new preprints were released in which Goodman and Hauschild Mosley \cite{GH108,GH208} use alternative topological and Jones basic construction theory type arguments to establish freeness of $\mathscr{B}_n^k$ and an isomorphism with the cyclotomic Kauffman tangle algebra. However, they require their ground rings to be integral domains with parameters satisfying the stronger conditions introduced by the authors in \cite{WY06}. In \cite{G08}, Goodman has also obtained cellularity results.

Rui and Xu \cite{RX08} have also proved freeness and cellularity of $\bnk$ when $k$ is odd, and later Rui and Si \cite{RS08} for general $k$, under the extra assumption that $\delta$ is invertible and using another condition called ``$\mathbf{u}$-admissibility''. The methods and arguments employed are strongly influenced by those used by Ariki, Mathas and Rui \cite{AMR06} for the cyclotomic Nazarov-Wenzl algebras, which are degenerate versions of the cyclotomic BMW algebras, and involve the construction of seminormal representations. 

It is not straightforward to compare the various notions of admissibility due to the few small but important differences in the assumptions on certain parameters of the ground ring (see Remarks below Definitions \ref{defn:bnk}, \ref{defn:adm} and \ref{defn:GHadm}). The classes of ground rings covered in the freeness and cellularity results of \cite{GH108,GH208, G08, RX08,RS08} are \emph{subsets} of the set of rings with admissible parameters as defined in this paper and \cite{Y07}. Moreover, Goodman and Rui-Si-Xu use a weaker definition of cellularity, to bypass a problem discovered in their original proofs relating to the anti-involution axiom of the original Graham-Lehrer definition.

The structure of the paper is as follows. In Section \ref{sec:prelims}, we introduce the cyclotomic BMW algebras and derive some straightforward identities and formulas pertinent to the next section.
%A natural problem to address first is whether these algebras are always free. We expect that the cyclotomic BMW algebras $\bnk$ should be free of rank $k^n(2n-1)!!$.
Section \ref{sect:spanning} is concerned with obtaining a spanning set of $\bnk$ of cardinality $k^n(2n-1)!!$. These two sections omit certain straightforward but tedious calculations, which can be found in \cite{Y07}. In Section \ref{sect:admissible}, we give the needed admissibility conditions explicitly (see Definition \ref{defn:adm}) and construct a ``generic'' ground ring $R_0$, in the sense that for any ring $R$ with admissible parameters there is a unique map $R_0 \rightarrow R$ which respects the parameters. We also shed some light on the relationships between the various admissibility conditions appearing in the literature at the end of Section~\ref{sect:admissible}. In Section \ref{sect:ckt}, we introduce the cyclotomic Kauffman tangle algebras.
The admissibility conditions are closely related to the existence of a nondegenerate (unnormalised) Markov trace function of $\mathscr{B}_n^k$, constructed in Section \ref{sect:traceconstruction}, which is then used together with the cyclotomic Brauer algebras in the linear independency arguments contained in Section \ref{sect:basis}. These nondegenerate Markov trace functions on $\bnk$ yields a family of Kauffman-type invariants of links in the solid torus; cf. Turaev \cite{T88}, tom Dieck \cite{T93}, Lambropoulou \cite{L99}.

\section{The cyclotomic BMW algebras $\bnk$} \label{sec:prelims}

In this section, we define the cyclotomic BMW algebras $\bnk$ and, through straightforward calculations and induction arguments, we establish several useful formulas and identities between special elements of the algebra. As seen in Definition \ref{defn:bnk} below, the defining relations of the algebra $\bnk$ consist of those for the BMW algebra $\bmw_n$, from Definition \ref{defn:bmw}, and further relations involving an extra generator $Y$ which satisfies a polynomial relation of order $k$. Throughout let us fix natural numbers $n$ and $k$. 

\begin{defn} \label{defn:bnk}
Let $R$ be a unital commutative ring 
containing units $q_0, q,\l$ and further 
elements $q_1, \ldots, q_{k-1}$ and $A_0, A_1, \ldots, A_{k-1}$ such that $\l - \l^{-1} = \d(1-A_0)$ holds, where $\d := q - q^{-1}$. \\ The \emph{\textbf{cyclotomic BMW algebra}} 
$\bnk := \bnk(q,\l,A_i,q_i)$ is the unital associative $R$-algebra 
generated by $Y^{\pm 1}, X_1^{\pm 1}, \ldots, X_{n-1}^{\pm 1}$ and 
$e_1, \ldots, e_{n-1}$ subject to the following relations, which hold for all possible values of $i$ unless otherwise stated.
\begin{eqnarray}
X_i - X_i^{-1} &=& \d(1-e_i) \label{eqn:ix} \\
X_iX_j &=& X_jX_i \qquad\qquad\qquad\quad \mspace{-3mu} \text{for } |i-j| \geq 2 \label{eqn:braid2} \\
X_iX_{i+1}X_i &=& X_{i+1}X_iX_{i+1}  \label{eqn:braid1} \\
X_ie_j &=& e_jX_i \ \qquad\qquad\qquad\quad \text{for } |i-j| \geq 2 \label{eqn:xecomm} \\
e_ie_j &=& e_je_i \qquad\qquad\qquad\quad\,\,\,\mspace{2mu} \text{for } |i-j| \geq 2 \label{eqn:eecomm}
\end{eqnarray}

\begin{eqnarray}
X_ie_i &=& e_iX_i \,\,\,=\,\,\, \l e_i \label{eqn:untwist} \\
X_i X_j e_i &=& e_j e_i \,\,\,=\,\,\, e_j X_i X_j \mspace{1mu}\,\qquad \text{for } |i-j| = 1 \label{eqn:xxe} \\
e_i e_{i\pm 1} e_i &=& e_i \label{eqn:eee} \\
e_i^2 &=& A_0e_i \\
Y^k \= \sum_{i=0}^{k-1} q_i Y^i \label{eqn:ycyclo} \\
X_1YX_1Y &=& YX_1YX_1  \label{eqn:braidb}\\
YX_i &=& X_iY \qquad\qquad\qquad\quad\mspace{5mu} \text{for } i > 1 \label{eqn:ycommx}\\
Ye_i &=& e_iY \qquad\qquad\qquad\quad\mspace{12mu} \text{for } i > 1 \label{eqn:ycomme}\\
YX_1Ye_1 &=& \lambda^{-1}e_1 \,\,\,=\,\,\, e_1YX_1Y \label{eqn:eyxy} \\
e_1Y^me_1 &=& A_m e_1 \qquad\qquad\qquad\quad\mspace{-1mu} \text{for } 0 \leq m \leq k-1. \label{eqn:eye}
\end{eqnarray}
If $R$ is a ring as in Definition \ref{defn:bnk}, we may use $\bnk(R)$ or $\bnk$ for short to denote the algebra $\bnk(q,\l,A_i,q_i)$.
\end{defn}

\textbf{Remarks:}\quad (1)\,\, There are slight but important differences in the parameters and the assumptions imposed on them in the literature.
The original definition of $\bnk$ given in H\"aring-Oldenburg~\cite{HO01} supposes that the $k^\mathrm{th}$ order polynomial relation (\ref{eqn:ycyclo}) splits over $R$; that is, $\prod_{i=0}^{k-1}(Y-p_i) = 0$, where the $p_i$ are units in the ground ring $R$. Under this relation, the $q_i$ in relation (\ref{eqn:ycyclo}) become the signed elementary symmetric polynomials in the $p_i$, where $q_0 = (-1)^{k-1} \prod_i p_i$ is invertible. However, we need not impose this stronger polynomial relation on $Y$ in the present work.
In addition to the splitting assumption, Goodman and Hauschild Mosley \cite{GH108,GH208} also assumes the invertibility of $A_0$ and that $\delta $ is not a zero divisor, and Rui-Si-Xu \cite{RX08,RS08} assume the invertibility of $\delta$ in $R$. Also, the assumption in \cite{Y07} that $A_0$ is invertible has been removed in this paper.

(2)\,\, When the relation (\ref{eqn:ycyclo}) is omitted, one obtains the affine BMW algebras, as studied by Goodman and Hauschild in \cite{GH06} algebras. In fact, in \cite{GH108,GH208,RX08,RS08}, the affine BMW algebra is initially considered with infinite parameters $\{A_j \mid j \geq 0\}$ instead and $e_1 Y^{j} e_1 = A_{j} e_1$, for all $j \geq 0$. The cyclotomic BMW algebra is then defined to be the quotient of this algebra by the ideal generated by the $k^\mathrm{th}$ order relation $\prod_{i=0}^{k-1}(Y-p_i) = 0$.

(3)\,\, Observe that, by relations (\ref{eqn:ix}) and (\ref{eqn:ycyclo}), it is unnecessary to include the inverses of $Y$ and $X$ as generators of $\bnk$ in Definition \ref{defn:bnk}. 

(4)\,\, Define $q_k := -1$. Then $\sum_{j=0}^k q_j Y^j = 0$ and the inverse of $Y$ may then be expressed as $Y^{-1} = - q_0^{-1} \sum_{i=0}^{k-1} q_{i+1}  Y^i$.

Using the defining $k^\mathrm{th}$ order relation on $Y$ and (\ref{eqn:eye}), there exists elements $A_{m}$ of $R$, \emph{for all} $m \in \mathbb{Z}$, such that 
\begin{equation} \label{eqn:negeye}
e_1 Y^{m} e_1 = A_{m} e_1.
\end{equation}
We will see later that, in order for our algebras to be ``well-behaved'', the $A_m$ cannot be chosen independently of the other parameters of the algebra.

(5)\,\, Observe that there is an unique anti-involution $^*: \bnk \rightarrow \bnk$ such that
\renewcommand{\theequation}{$\ast$}
\begin{equation} \label{eqn:star}
  Y^* = Y, \quad X_i^* = X_i \quad \text{and} \quad e_i^* = e_i,
\end{equation}
for every $i = 1,\ldots,n-1$. Here an anti-involution always means an involutary $R$-algebra anti-automorphism. \\
%Also, the following defines an $R$-algebra isomorphism \[\bnk(q,\l,A_i,q_i) \rightarrow %\bnk(q^{-1},\l^{-1},A_{-i},-q_{k-i}q_0^{-1}):\]
%\[  
%Y \mapsto \iy, \quad X_i \mapsto \ix, \quad e_i \mapsto e_i.
%\]
\vskip 0.2cm
\noi For all $i = 1,\ldots, n$, define the following elements of $\bnk$:
\[ Y_i':= X_{i-1} \ldots X_2 X_1 Y X_1 X_2 \ldots X_{i-1}.\]
Observe that these elements are fixed under the (\ref{eqn:star}) anti-involution. We now establish several identities in the algebra which will be used frequently in future proofs, including the pairwise commutativity of the $Y_i'$, which is their most important and useful property. Let us fix $n$ and~$k$. The following calculations are valid over a general ring $R$ with any choice of the above parameters $A_0, \o, A_{k-1}, q_0, \o, q_{k-1}, q, \l$.
\renewcommand{\theequation}{\arabic{equation}}
\addtocounter{equation}{-1}
\begin{prop} \label{prop:conseq}
The following relations hold in $\bnk$, for all $i$, $j$, and $p$ unless otherwise stated.
\begin{equation} \label{eqn:xi2}
  X_i^2 = 1+ \d X_i - \d \l e_i
\end{equation}
\begin{equation} \label{eqn:exe} 
e_iX_{i\pm 1}e_i = \l^{-1} e_i 
\end{equation}
\begin{equation} \label{eqn:prop1b}
  X_i Y_j' = Y_j' X_i \quad\text{and}\quad e_i Y_j' = Y_j' e_i \text{,\quad when $i \neq j$ or $j-1$,}
\end{equation}
\begin{equation} \label{eqn:prop1c}
Y_i' Y_j' = Y_j' Y_i' \text{,\quad for all } i,j,
\end{equation} 
\begin{equation} \label{eqn:prop1d}
Y_i' X_i Y_i' e_i = \l^{-1} e_i = e_i Y_i' X_i Y_i'
\end{equation}
\begin{equation} \label{eqn:m9}
    e_i \y{i+1}{p} = e_i \y{i}{-p} \quad \text{ and } \quad \y{i+1}{p} e_i = \y{i}{-p} e_i. 
\end{equation}
\end{prop}

\begin{proof}
The quadratic relation (\ref{eqn:xi2}) follows by multiplying relation (\ref{eqn:ix}) by $X_i$ and applying relation (\ref{eqn:untwist}) to simplify. Equation (\ref{eqn:exe}) is proved below.
\[
  e_i X_{i\pm 1} e_i \rel{eqn:untwist} \l^{-1} e_i X_{i+1} X_i e_i \rel{eqn:xxe} \l^{-1} e_i e_{i+1} e_i \rel{eqn:eee} \l^{-1} e_i.
\]
\noi The first equation in (\ref{eqn:prop1b}) follows from the braid relations (\ref{eqn:braid2}), (\ref{eqn:braid1}) and (\ref{eqn:ycommx}) and the second follows from relations (\ref{eqn:xecomm}), (\ref{eqn:xxe}) and (\ref{eqn:ycomme}).

\noi Equation (\ref{eqn:prop1c}) follows from (\ref{eqn:prop1b}) and the braid relation (\ref{eqn:braidb}). 

\noi We prove (\ref{eqn:prop1d}) by induction on $i \geq 1$. The case where $i=1$ is simply relation (\ref{eqn:eyxy}). Now assume (\ref{eqn:prop1d}) holds for a fixed $i$. Thus, remembering that $Y_{i+1}' = X_i Y_i' X_i$ and applying equations (\ref{eqn:braid1}), (\ref{eqn:prop1b}) then (\ref{eqn:xxe}) gives
\[
Y_{i+1}' X_{i+1} Y_{i+1}' e_{i+1} 
= X_i X_{i+1} Y_i' X_i Y_i' e_i e_{i+1}
\stackrel{\tiny{\text{ind. hypo.}}}{=} \l^{-1} X_i X_{i+1} e_i e_{i+1}
\stackrel{\tiny{(\ref{eqn:xxe}),(\ref{eqn:eee})}}{=} \l^{-1} e_{i+1}. 
\]
The second equality of (\ref{eqn:prop1d}) now follows immediately by applying the anti-involution (\ref{eqn:star}) to the first.
Moreover, (\ref{eqn:m9}) follows from parts (\ref{eqn:prop1c}) and (\ref{eqn:prop1d}), remembering that $Y_{j+1}' = X_j Y_j' X_j$.
\end{proof}
 
In the remainder of this section, we present some useful identities involving the $Y_i'$, $X_i$ and $e_i$ which shall be used extensively throughout later proofs. The proof of the following Proposition involves straightforward application of the relations in Definition \ref{defn:bnk} and shall be left as an exercise to the reader; full details can be found in Proposition 1.3 of \cite{Y07}.
\begin{prop} \label{prop:XXp1}
The following equations hold for all $i$:
\begin{eqnarray}
e_i e_{i+1} e_{i+2} \gamma_i \= \gamma_{i+2} e_i e_{i+1} e_{i+2}, \quad \text{where $\gamma_i = X_i,e_i$ or $Y_i'$}; \label{eqn:ep2} \\
X_i X_{i+1} \gamma_i  \= \gamma_{i+1} X_i X_{i+1}, \quad \text{where $\gamma_i = X_i$ or $e_i$}. \label{eqn:XXp1}
\end{eqnarray} 
\end{prop}
\begin{comment}
\begin{proof}
We prove (\ref{eqn:ep2}) in the three stated cases.
\[
e_i e_{i+1} e_{i+2} X_i \rel{eqn:xecomm} e_i e_{i+1} X_i e_{i+2} 
\stack{(\ref{eqn:xxe}),(\ref{eqn:eee})} e_i X_{i+1}^{-1} e_{i+2} 
\stack{(\ref{eqn:xxe}),(\ref{eqn:eee})} e_i X_{i+2} e_{i+1} e_{i+2} 
\rel{eqn:xecomm} X_{i+2} e_i e_{i+1} e_{i+2}. 
\]
\[
e_i e_{i+1} e_{i+2} e_i \rel{eqn:eecomm} e_i e_{i+1} e_i e_{i+2} 
\rel{eqn:eee} e_i e_{i+2} 
\rel{eqn:eee} e_i e_{i+2} e_{i+1} e_{i+2} 
\rel{eqn:eecomm} e_{i+2} e_i e_{i+1} e_{i+2}. 
\]
\[
e_i e_{i+1} e_{i+2} Y_i' \rel{eqn:prop1b} e_i Y_i' e_{i+1} e_{i+2} 
\rel{eqn:m9} e_i \y{i+1}{-1} e_{i+1}e_{i+2} 
\rel{eqn:m9} e_i Y_{i+2}' e_{i+1}e_{i+2}
\rel{eqn:prop1b} Y_{i+2}' e_i e_{i+1} e_{i+2}. 
\]
Equation (\ref{eqn:XXp1}) follows clearly from relations (\ref{eqn:braid1}) and (\ref{eqn:xxe}). 
\end{proof}
\end{comment}

\begin{lemma}\label{lemma:mult}
The following hold for any $i$ and non-negative integer $p$: 
\begin{align}
X_i \y{i}{p} &= \y{i+1}{p} X_i - \d \pum{s} \y{i+1}{s} \y{i}{p-s} + \d \pum{s} \y{i+1}{s} e_i \y{i}{p-s} \label{eqn:magic} \\
X_i  \y{i}{-p} &= \y{i+1}{-p} X_i + \d \pum{s} \y{i+1}{s-p} \y{i}{-s} - \d \pum{s} \y{i+1}{s-p} e_i \y{i}{-s} \label{eqn:m2} \\
%\ix \y{i}{p} &= \y{i+1}{p} \ix - \d \pum{s} \y{i+1}{p-s} \y{i}{s} + \d \pum{s} \y{i+1}{p-s} e_i \y{i}{s} \label{eqn:m3} \\
%\ix \y{i}{-p} &= \y{i+1}{-p} \ix + \d \pum{s} \y{i+1}{-s} \y{i}{-(p-s)} - \d \pum{s} \y{i+1}{-s} e_i \y{i}{-(p-s)} \label{eqn:m4} \\
X_i \y{i+1}{p} &= \y{i}{p} X_i + \d \pum{s} \y{i}{p-s} \y{i+1}{s} - \d \pum{s} \y{i}{p-s} e_i \y{i+1}{s} \label{eqn:m5} \\
X_i \y{i+1}{-p} &= \y{i}{-p} X_i - \d \pum{s} \y{i}{-s} \y{i+1}{s-p} + \d \pum{s} \y{i}{-s} e_i \y{i+1}{s-p} \label{eqn:m6} \\
%\ix \y{i+1}{p} &= \y{i}{p} \ix + \d \pum{s} \y{i}{s} \y{i+1}{p-s} - \d \pum{s} \y{i}{s} e_i \y{i+1}{p-s} \label{eqn:m7} \\
%\ix \y{i+1}{-p} &= \y{i}{-p} \ix - \d \pum{s} \y{i}{-(p-s)} \y{i+1}{-s} + \d \pum{s} \y{i}{-(p-s)} e_i \y{i+1}{-s} \label{eqn:m8} \\
X_i \y{i}{p} X_i &= \y{i+1}{p} - \d \sum_{s=1}^{p-1} \y{i+1}{s} \y{i}{p-s} X_i + \d \sum_{s=1}^{p-1} \y{i+1}{s} e_i \y{i}{p-s} X_i \label{eqn:m10} \\
X_i \y{i}{p} X_i &= \y{i+1}{p} - \d \sum_{s=1}^{p-1} X_i \y{i}{s} \y{i+1}{p-s} + \d \sum_{s=1}^{p-1} X_i \y{i}{s} e_i \y{i+1}{p-s} \label{eqn:m11} \\
X_i \y{i}{-p} X_i &= \y{i+1}{-p} + \d \sum_{s=0}^{p} \y{i+1}{s-p} \y{i}{-s} X_i - \d \sum_{s=0}^{p} \y{i+1}{s-p} e_i \y{i}{-s} X_i \label{eqn:m12} \\
X_i \y{i}{-p} X_i &= \y{i+1}{-p} + \d \sum_{s=0}^{p} X_i \y{i}{-s} \y{i+1}{s-p} - \d \sum_{s=0}^{p} X_i \y{i}{-s} e_i \y{i+1}{s-p}. \label{eqn:m13} 
\end{align}
\end{lemma}
\begin{proof}
We obtain the first equation through the following straightforward calculation. For all $p \geq 0$,
\be %\begin{align*}
  X_i \y{i}{p} &\p Y'_{i+1} \ix \y{i}{p-1} \\ 
  &\rel{eqn:ix} Y'_{i+1} X_i \y{i}{p-1} - \d Y'_{i+1} \y{i}{p-1} + \d Y'_{i+1} e_i \y{i}{p-1} \\ 
&\p \y{i+1}{2}X_i \y{i}{p-2} - \d \y{i+1}{2} \y{i}{p-2} + \d \y{i+1}{2} e_i \y{i}{p-2} - \d Y'_{i+1} \y{i}{p-1} + \d Y'_{i+1} e_i \y{i}{p-1} \\
&\p \ldots \, = \\
&\p \y{i+1}{p-1} X_i Y_i' - \d \sum_{s=1}^{p-1} \y{i+1}{s} \y{i}{p-s} + \d \sum_{s=1}^{p-1} \y{i+1}{s} e_i \y{i}{p-s} \\
&\p \y{i+1}{p} X_i - \d \sum_{s=1}^{p} \y{i+1}{s} \y{i}{p-s} + \d \sum_{s=1}^{p} \y{i+1}{s} e_i \y{i}{p-s}\text{,\qquad proving }(\ref{eqn:magic}).
\end{align*}

Multiplying equation (\ref{eqn:magic}) on the left by $\y{i+1}{-p}$ and the right by $\y{i}{-p}$ and rearranging gives equation (\ref{eqn:m2}). 
Applying (\ref{eqn:star}) to equations (\ref{eqn:magic}) and (\ref{eqn:m2}) and rearranging then produces equations (\ref{eqn:m5}) and (\ref{eqn:m6}), respectively. Equation (\ref{eqn:m10}) follows as an easy consequence of equations (\ref{eqn:magic}) and (\ref{eqn:xi2}). Furthermore, applying (\ref{eqn:star}) to (\ref{eqn:m10}) and a straightforward change of summation now gives equation (\ref{eqn:m11}). Similarly, using equation (\ref{eqn:m2}) and (\ref{eqn:star}), one obtains equations (\ref{eqn:m12}) and (\ref{eqn:m13}).
\end{proof}

\textbf{Notation.}\quad
In this paper, we shall adopt the following notation conventions.
If $J$ is a subset of an $R$-module, $\la{J}$ is used to denote the $R$-span of $J$.
Finally, for a subset $S \subseteq R$, we denote $\la{S}_R$ to be the \emph{ideal} generated by $S$ in $R$ and only omit the subscript $R$ if it does not create any ambiguity in the current context. 

\begin{lemma} \label{lemma:eype}
For all integers $p$, the following hold:
\begin{enumerate}
  \item[(I)] $e_i \y{i}{p} e_i \in \la{Y^{s_1} \y{2}{s_2} \ldots \y{i-1}{s_{i-1}} e_i}$;
  \item[(II)] $X_i \y{i}{p} e_i \in \la{Y^{s_1} \y{2}{s_2} \ldots \y{i-1}{s_{i-1}} \y{i}{s_i} e_i \,\big|\,  |s_i| \leq |p|\,}$;
  \item[(III)] $e_i \y{i}{p} X_i \in \la{e_i Y^{s_1} \y{2}{s_2} \ldots \y{i-1}{s_{i-1}} \y{i}{s_i} \,\big|\,  |s_i| \leq |p|\,}$.
\end{enumerate}
\end{lemma}

\begin{proof}
The $k^\mathrm{th}$ order relation on $Y$ and relation (\ref{eqn:eye}) tells us that $e_1 Y^{p} e_1$ is always a scalar multiple of $e_1$, for any integer $p$, hence showing part (I) of the lemma for the case $i=1$.

\noi Now, for all $p \geq 0$, equation (\ref{eqn:magic}) and (\ref{eqn:untwist}) implies that
\begin{eqnarray*}
X_1 Y^p e_1 &=& \l \y{2}{p} e_1 - \d \pum{s} \y{2}{s}Y^{p-s} e_1 + \d \pum{s} \y{2}{s} e_1 Y^{p-s} e_1 \\
&\stack{(\ref{eqn:m9}),(\ref{eqn:eye})}& \l Y^{-p} e_1 - \d \pum{s} Y^{p-2s} e_1 + \d \pum{s} A_{p-s} Y^{-s} e_1.
\end{eqnarray*}

\noi Similarly, by equations (\ref{eqn:m2}), (\ref{eqn:untwist}) and (\ref{eqn:negeye}),
\[
X_1 Y^{-p} e_1 = \l Y^p e_1 + \d \pum{s} Y^{p-2s} e_1 - \d \pum{s} A_{-s} Y^{p-s} e_1 \text{, for all }p \geq 0.
\]

\noi Observe that $|\!-\! p\kern1.1pt| = |p|$ and when $1 \leq s \leq p$, we have $|s|$,$|p-s|$, $|p-2s| \leq |p|$. Hence $X_1 Y^p e_1 \in \la{Y^m e_1 \,\big|\, |m| \leq |p|\,}$, for all $p \in \mathbb{Z}$, proving part (II) of the lemma for the case $i=1$.

We are now able to prove (I) and (II), for all integers $p \geq 0$, together by induction on~$i$, which in turn involves inducting on $p \geq 0$. By relations (\ref{eqn:eye}) and (\ref{eqn:untwist}), both hold clearly for~$p = 0$. \\ Now let us assume that: \\ 
$X_i \y{i}{r} e_i \in \la{Y^{s_1} \y{2}{s_2} \ldots \y{i-1}{s_{i-1}} \y{i}{s_i} e_i \,\big|\,  |s_i| \leq |r|\,}$ and $e_{i} \y{i}{r} e_{i} \in \la{Y^{s_1} \y{2}{s_2} \ldots \y{i-1}{s_{i-1}} e_{i}}$, for all $r < p$, and $X_{i-1} \y{i-1}{r} e_{i-1} \in \la{Y^{s_1} \y{2}{s_2} \ldots \y{i-2}{s_{i-2}} \y{i-1}{s_{i-1}} e_{i-1} \,\big|\,  |s_{i-1}| \leq |r|\,}$ and $e_{i-1} \y{i-1}{r} e_{i-1} \in \la{Y^{s_1} \y{2}{s_2} \ldots \y{i-2}{s_{i-2}} e_{i-1}}$, for all $r \geq 0$.%since $e_i^2 = A_0 e_i$ and $X_i e_i = \l e_i$, for all $i$. Now let us assume that: $X_i \y{i}{r} e_i \in \la{Y^{s_1} \y{2}{s_2} \ldots \y{i-1}{s_{i-1}} \y{i}{s_i} e_i \,\big|\,  |s_i| \leq |r|\,}$ and $e_{i} \y{i}{r} e_{i} \in \la{Y^{s_1} \y{2}{s_2} \ldots \y{i-1}{s_{i-1}} e_{i}}$, for all $r < p$, and $X_{i-1} \y{i-1}{r} e_{i-1} \in \la{Y^{s_1} \y{2}{s_2} \ldots \y{i-2}{s_{i-2}} \y{i-1}{s_{i-1}} e_{i-1} \,\big|\,  |s_{i-1}| \leq |r|\,}$ and $e_{i-1} \y{i-1}{r} e_{i-1} \in \la{Y^{s_1} \y{2}{s_2} \ldots \y{i-2}{s_{i-2}} e_{i-1}}$, for all $r \geq 0$.

\noi Recall that $Y_i' = X_{i-1} Y_{i-1}' X_{i-1}$. Using this and equation (\ref{eqn:m5}), followed by relations (\ref{eqn:ix}), (\ref{eqn:xxe}) then (\ref{eqn:prop1b}) we see that, for all $p > 0$,
\begin{align}
e_i \y{i}{p} e_i &\h e_i e_{i-1} \y{i-1}{p} \ix X_{i-1}^{-1} e_i + \d \ppum{s} e_i X_{i-1} \y{i-1}{p-s} \y{i}{s} e_i - \d \ppum{s} e_i X_{i-1} \y{i-1}{p-s} e_{i-1} \y{i}{s} e_i \nonumber \\
&\h e_i e_{i-1} \y{i-1}{p} e_{i-1} e_i + \d \ppum{s} e_i X_{i-1} \y{i}{s} e_i \y{i-1}{p-s}  - \d \ppum{s} e_i X_{i-1} \y{i-1}{p-s} e_{i-1} \y{i}{s} e_i, \label{eqn:eype}
\end{align}
by relation (\ref{eqn:xxe}) and Proposition \ref{prop:conseq}.

\noi Let us consider the first term in the latter equation above. By induction on $i$, 
\[
e_{i-1} \y{i-1}{p} e_{i-1} \in \la{Y^{s_1} \y{2}{s_2} \ldots \y{i-2}{s_{i-2}} e_{i-1}}.
\]
Therefore, by relation (\ref{eqn:eee}),
$e_i e_{i-1} \y{i-1}{p} e_{i-1} e_i \in \la{Y^{s_1} \y{2}{s_2} \ldots \y{i-2}{s_{i-2}} e_i}$.
Now let us consider the second term in the RHS of (\ref{eqn:eype}). Fix $0 \leq s \leq p-1$.
\[
e_i X_{i-1} \y{i}{s} e_i \y{i-1}{p-s}
\stack{(\ref{eqn:xxe}), (\ref{eqn:ix})} e_i e_{i-1} X_i \y{i}{s} e_i \y{i-1}{p-s} - \d e_i e_{i-1} \y{i}{s} e_i \y{i-1}{p-s} + \d e_i e_{i-1} e_i \y{i}{s} e_i \y{i-1}{p-s}. 
\]

\noi By induction on $p$ and equation (\ref{eqn:prop1b}),\vspace{-1mm}
\begin{align*}
e_i e_{i-1} X_i \y{i}{s} e_i \y{i-1}{p-s} 
&\rel{eqn:m9} \la{e_i Y^{m_1} \y{2}{m_2} \ldots \y{i-2}{m_{i-2}} e_{i-1} \y{i-1}{m_{i-1}-m_i} e_i \y{i-1}{p-s} \,\big|\, |m_i| \leq |s|\,} \\
&\h \la{Y^{m_1} \y{2}{m_2} \ldots \y{i-2}{m_{i-2}} e_i e_{i-1} e_i \y{i-1}{p-s+m_{i-1}-m_i} \,\big|\, |m_i| \leq |s|}.
\end{align*}

\noi Therefore $e_i e_{i-1} X_i \y{i}{s} e_i \y{i-1}{p-s} \in \la{Y^{m_1} \y{2}{m_2} \ldots \y{i-2}{m_{i-2}} \y{i-1}{p-s+m_{i-1}-m_i} e_i}$.

\noi Also, by (\ref{eqn:m9}), (\ref{eqn:prop1b}) and (\ref{eqn:eee}),
$e_i e_{i-1} \y{i}{s} e_i \y{i-1}{p-s} = e_i \y{i-1}{p-2s}$ and, by (\ref{eqn:eee}) and induction on $p$, we have $e_i e_{i-1} e_i \y{i}{s} e_i \y{i-1}{p-s} \stackrel{\phantom{(1)}}{\in} \la{Y^{m_1} \y{2}{m_2} \ldots \y{i-1}{m_{i-1}+p-s} e_i}$.

\noi Thus, for all $0 \leq s \leq p-1$,
$
e_i X_{i-1} \y{i}{s} e_i \y{i-1}{p-s} \in \la{Y^{s_1} \y{2}{s_2} \ldots \y{i-1}{s_{i-1}} e_i}.
$
Hence the second term in the RHS of equation (\ref{eqn:eype}) is in $\la{Y^{s_1} \y{2}{s_2} \ldots \y{i-1}{s_{i-1}} e_i}$.

\noi Finally, by induction on $i$ and using (\ref{eqn:m9}), (\ref{eqn:prop1b}) and (\ref{eqn:eee}),
\be
e_i X_{i-1} \y{i-1}{p-s} e_{i-1} \y{i}{s} e_i
&\in \la{e_i Y^{m_1} \y{2}{m_2} \ldots \y{i-1}{m_{i-1}} e_{i-1} \y{i}{s} e_i \,\big|\, |m_{i-1}| \leq |p-s|\,} \\
&\in \la{Y^{m_1} \y{2}{m_2} \ldots \y{i-1}{m_{i-1}-s} e_i \,\big|\, |m_{i-1}| \leq |p-s|\,}.
\end{align*}
Thus the third term in the RHS of equation (\ref{eqn:eype}) is in $\la{Y^{s_1} \y{2}{s_2} \ldots \y{i-1}{s_{i-1}} e_i}$.
%Hence, for all $p \geq 0$, $e_i \y{i}{p} e_i \in \la{Y^{s_1} \y{2}{s_2} \ldots \y{i-1}{s_{i-1}} e_i \,\big|\, s_j \in \mathbb{Z}}$.

\noi Also, for all $p \geq 0$, equation (\ref{eqn:magic}) implies that \vspace{-2mm}
\begin{eqnarray*}
  X_i \y{i}{p} e_i \= \y{i+1}{p} X_i e_i - \d \pum{s} \y{i+1}{s} \y{i}{p-s} e_i + \d \pum{s} \y{i+1}{s} e_i \y{i}{p-s} e_i \\ 
&\stackrel{\tiny{(\ref{eqn:untwist}),(\ref{eqn:m9})}}{=}& \l \y{i}{-p} e_i - \d \pum{s} \y{i}{p-2s} e_i + \d \pum{s} \y{i}{-s} e_i \y{i}{p-s} e_i.
\end{eqnarray*}
The first term above is clearly in $\la{Y^{s_1} \y{2}{s_2} \ldots \y{i-1}{s_{i-1}} \y{i}{s_i} e_i \,\big|\, |s_i| \leq |p|\,}$, since $|\!-\! p\kern1pt| = |p\,|$. \\
Regarding the second term above, since $1 \leq s \leq p$, $|p-2s| \leq |p|$, so it is also an element of $\la{Y^{s_1} \y{2}{s_2} \ldots \y{i-1}{s_{i-1}} \y{i}{s_i} e_i \,\big|\, |s_i| \leq |p|\,}$.
Moreover, $0 \leq p-s \leq p-1 < p$, so by induction on $p$,\vspace{-3mm}
\[
  \y{i}{-s} e_i \y{i}{p-s} e_i 
\in \la{\y{i}{-s} Y^{m_1} \y{2}{m_2} \ldots \y{i-1}{m_{i-1}} e_i} \subseteq \la{Y^{s_1} \y{2}{s_2} \ldots \y{i-1}{s_{i-1}} \y{i}{s_i} e_i \,\big|\, |s_i| \leq |p|\,}.
\]
Therefore, for all $p \geq 0$, \vspace{-1mm}
\[
X_i \y{i}{p} e_i \in \la{Y^{s_1} \y{2}{s_2} \ldots \y{i-1}{s_{i-1}} \y{i}{s_i} e_i \,\big|\, |s_i| \leq |p|\,}
\text{\quad and \quad}
e_i \y{i}{p} e_i \in \la{Y^{s_1} \y{2}{s_2} \ldots \y{i-1}{s_{i-1}} e_i}.
\]

Let us denote $\dagger: \bnk(q^{-1},\l^{-1},A_{-i},-q_{k-i}q_0^{-1}) \rightarrow \bnk(q,\l,A_i,q_i)$ to be the isomorphism of $R$-algebras defined by 
\[
  Y \mapsto \iy, \quad X_i \mapsto \ix, \quad e_i \mapsto e_i.
%  q &\mapsto& q^{-1}, \quad \d \mapsto -\d, \quad \l \mapsto \l^{-1}, \quad A_i \mapsto A_{-i} \quad \text{and} \quad p_i \mapsto p_i^{-1}.
\]
Note that $\dagger$ maps $Y_i'$ to its inverse. \\
For all $p \geq 0$, we have shown above that $e_i \y{i}{p} e_i \in \la{Y^{s_1} \y{2}{s_2} \ldots \y{i-1}{s_{i-1}} e_i}$, as an element of $\bnk(q^{-1},\l^{-1},A_{-i},-q_{k-i}q_0^{-1})$. Therefore, using $\dagger$,
\begin{equation} \label{eqn:eypeneg} \vspace{-3mm}
e_i \y{i}{-p} e_i \in \la{Y^{s_1} \y{2}{
s_2} \ldots \y{i-1}{s_{i-1}} e_i},
\end{equation} 
as an element of $\bnk(q,\l,A_i,q_i)$, for all $p \geq 0$. 

\noi Furthermore, our previous work also shows that,
\[
X_i \y{i}{p} e_i \in \la{Y^{s_1} \y{2}{s_2} \ldots \y{i-1}{s_{i-1}} \y{i}{s_i} e_i \,\big|\, |s_i| \leq |p|\,},
\]
as an element of $\bnk(q^{-1},\l^{-1},A_{-i},-q_{k-i}q_0^{-1})$. Thus, applying $\dagger$ implies
\[
\ix \y{i}{-p} e_i \in \la{Y^{s_1} \y{2}{s_2} \ldots \y{i-1}{s_{i-1}} \y{i}{s_i} e_i \,\big|\, |s_i| \leq |p|\,},
\]
as an element of $\bnk(q,\l,A_i,q_i)$.

\noi However, by relation (\ref{eqn:ix}), $\ix \y{i}{-p} e_i = X_i \y{i}{-p} e_i - \d \y{i}{-p} e_i + \d e_i \y{i}{-p} e_i$. By (\ref{eqn:eypeneg}), the last two terms are clearly in $\la{Y^{s_1} \y{2}{s_2} \ldots \y{i-1}{s_{i-1}} \y{i}{s_i} e_i \,\big|\, |s_i| \leq |p|\,}$.
So, as an element of $\bnk(q,\l,A_i,q_i)$,
\[
X_i \y{i}{-p} e_i \in \la{Y^{s_1} \y{2}{s_2} \ldots \y{i-1}{s_{i-1}} \y{i}{s_i} e_i \,\big|\, |s_i| \leq |p|\,}.
\]
Hence this concludes the proof of (I) and (II) for all integers $p$. Applying (\ref{eqn:star}) to part (II) immediately gives part (III) of the Lemma.
\end{proof}

\section{Spanning sets of {$\mathscr{B}_{n}^{k}$}} \label{sect:spanning}

In this section, we produce a spanning set of $\bnk(R)$ for any ring $R$, as in Definition \ref{defn:bnk}, of cardinality $k^n(2n-1)!! = k^n (2n-1) \cdot (2n-3) \cdots 3 \cdot 1$. Hence this shows the rank of $\bnk$ is at most $k^n(2n-1)!!$. The spanning set we obtain involves picking any basis of the Ariki-Koike algebras, which we define below. In Section \ref{sect:basis}, we will see that these spanning sets are linearly independent provided we impose ``admissibility conditions" on the parameters of $R$, which shall be analysed in the next section. This section contains many straightforward but lengthy calculations; for full complete details, we refer the reader to \cite{Y07}. Finally, we note here that our spanning sets differs from that obtained by Goodman and Hauschild Mosley in~\cite{GH108}. 

\begin{defn}
For any unital commutative ring $R$ and $q', q_0, \o, q_{k-1} \in R$. The \textbf{Ariki-Koike algebra} $\mathfrak{h}_{n,k}(R)$ denote the unital associative $R$-algebra with generators $T_0^{\pm 1}$, $T_1^{\pm 1}$, \o, $T_{n-1}^{\pm 1}$ and relations
\[
\begin{array}{rcll}
  T_0 T_1 T_0 T_1 &=& T_1 T_0 T_1 T_0& \\
  T_i T_{i \pm 1} T_i \= T_{i \pm 1} T_i T_{i \pm 1}&\quad \text{for } i = 1,\o,n-2\\
  T_i T_j \= T_j T_i&\quad \text{for } |i-j|\geq 2 \\
%  \prod_{i=0}^{k-1}(T_0-p_i) &=& 0 &\\
T_0^k \= {\displaystyle \sum_{i=0}^{k-1} q_i T_0^i} \\
   T_i^2 \= (q'-1)T_i + q'&\quad \text{for } i = 1,\o,n-2.
%  T_i^2 &=& \d T_i + 1. 
\end{array}
\]
\end{defn}

The algebras $\mathfrak{h}_{n,k}$ are also referred to as the `cyclotomic Hecke algebras of type $G(k,1,n)$' and were introduced independently by Ariki and Koike \cite{AK94} and Brou\'e and Malle \cite{BM93}. They may be thought of as the Iwahori-Hecke algebras corresponding to the complex reflection group $(\mathbb{Z}/k\mathbb{Z}) \wr \mathfrak{S}_n$, the wreath product of the cyclic group $\mathbb{Z}/k\mathbb{Z}$ of order $k$ with the symmetric group $\mathfrak{S}_n$ of degree $n$. Indeed, by considering the case when $q' = 1$, $q_0 = 1$ and $q_i = 0$, one recovers the group algebra of $(\mathbb{Z}/k\mathbb{Z}) \wr \mathfrak{S}_n$. Also, it is isomorphic to the Iwahori-Hecke algebra of type $A_{n-1}$ or $B_n$, when $k = 1$ or $2$, respectively. Ariki and Koike \cite{AK94} prove that it is $R$-free of rank $k^n n!$, the cardinality of $(\mathbb{Z}/k\mathbb{Z}) \wr \mathfrak{S}_n$. In addition, they classify its irreducible representations and give explicit matrix representations in the generic semisimple setting. Also, Graham and Lehrer \cite{GL96} and Dipper, James and Mathas \cite{DJM98} prove that the algebra is cellular. The modular representation theory of these algebras have also been studied extensively in the literature.
 
Now suppose $R$ is a ring as in the definition of $\bnk$ and let $q' := q^2$. Then, from the given presentations of the algebras, it is straightforward to show that $\ak{n}(R)$ is a quotient of $\bnk(R)$ under the following projection
\begin{eqnarray*}
  \hspace{2cm}\pi_n: \bnk &\rightarrow& \ak{n} \\
  Y &\mapsto& T_0, \\
  X_i &\mapsto& q^{-1} T_i, \text{ \quad for }1 \leq i \leq n-1  \\
  e_i &\mapsto& 0.
\end{eqnarray*}
Indeed, $\bnk/I \cong \mathfrak{h}_{n,k}$ as $R$-algebras, where $I$ is the two-sided ideal generated by the $e_i$'s in $\bnk(R)$. (Remark: due to relation (\ref{eqn:eee}), it is clear that $I$ is actually equal to the two-sided ideal generated by just a single fixed $e_j$).

Our aim in this section is to obtain a spanning set of $\bnk$, for any choice of basis $\X_{m,k}$ for any $\mathfrak{h}_{m,k}$. 
For any basis $\X_{n,k}$ of $\mathfrak{h}_{n,k}$, let $\overline{\X}_{n,k}$ be any subset of $\bnk$ mapping onto $\X_{n,k}$ of the same cardinality. 
Also, for any $l \leq n$, there is a natural map $\blk \rightarrow \bnk$. Let $\bl{l}$ denote the image of $\blk$ under this map; that is, $\bl{l}$ is the subalgebra of $\bnk$ generated by $Y, X_1, \o, X_{l-1}, e_1, \o, e_{l-1}$. Note that \emph{a priori} it is not clear that this map is injective; i.e., that $\bl{l}$ is isomorphic to $\blk$. In fact, over a specific class of ground rings, this will follow as a consequence of freeness of $\bnk$, which is established in Section~\ref{sect:basis}.

Finally, define $\widetilde{\X}_{l,k}$ to be the image of $\overline{\X}_{l,k}$ in $\bnk$.
The goal of this section is to prove the following theorem.

\begin{thm} \label{thm:span}
The set of elements of the following form spans $\bnk$.
\vspace{-0.1cm}
\begin{eqnarray*}
\y{i_1}{s_1} \y{i_2}{s_2} \, \ldots \, \y{i_m}{s_m} \!\!\!\!\!&&\!\!(X_{i_1} \ldots X_{j_{\scriptscriptstyle 1}-1} e_{j_1} \ldots e_{n-2} e_{n-1}) \\
% \!\!\!\!\!&&\!\! (X_{i_2} \ldots X_{j_{\scriptscriptstyle 2}-1} e_{j_2} \ldots e_{n-4} e_{n-3}) \ldots \\
 \ldots \!\!\!\!\!&&\!\! (X_{i_m} \ldots X_{j_{\scriptscriptstyle m}-1} e_{j_m} \ldots e_{n-2m} e_{n-2m+1})\,\, \chi^{(n-2m)} \\
\!\!\!\!\!&&\!\! (e_{n-2m} \ldots e_{h_{\scriptscriptstyle m}} X_{h_{\scriptscriptstyle m} - 1} \ldots X_{g_m}) \ldots \\
%\ldots \!\!\!\!\!&&\!\! (e_{n-3} e_{n-4} \ldots e_{h_{\scriptscriptstyle 2}} \ldots X_{g_2}) \\
\!\!\!\!\!&&\!\! (e_{n-2} \ldots e_{h_{\scriptscriptstyle 1}} X_{h_{\scriptscriptstyle 1} - 1} \ldots X_{g_1}) \,\,\, \y{g_m}{t_m} \y{g_2}{t_2} \, \ldots \, \y{g_1}{t_1},
\end{eqnarray*}
where $m = 1,2,\ldots, \lfloor \frac{n}{2} \rfloor$, $i_1 > i_2 > \ldots > i_m$, $g_1 > g_2 > \ldots > g_m$ and, for each $f = 1,2,\ldots m$, we require $1 \leq i_f\leq j_f \leq n-2f+1$, $1 \leq g_f\leq h_f \leq n-2f+1$, $s_f, t_f \in \left\{ \floor{\frac{k}{2}} - (k-1), \o, \floor{\frac{k}{2}} \right\}$ and $\chi^{(n-2m)}$ is an element of $\widetilde{\X}_{n-2m,\,k}$. 
\end{thm}

To make the spanning set above more palatable for now, we introduce the following notation and relate parts of the expression diagrammatically where possible.
Suppose $l \geq 1$. Let $i$ and $j$ be such that $i \leq j \leq l+1$ and $p$ be any integer. Define
\[
\a{ijl}{p} := \y{i}{p}  X_i \o X_{j-1} e_j \o e_l.
\]
Then Theorem \ref{thm:span} states that the algebra $\bnk$ is spanned by the set of elements
\[
\alpha_{i_1,j_1,n-1}^{s_1}\ldots\alpha_{i_m,j_m,n-2m+1}^{s_m} \chi^{(n-2m)}(\alpha_{g_m,h_m,n-2m}^{t_m})^*\ldots (\alpha_{g_1,h_1,n-2}^{t_1})^*,
\]
with conditions specified as above. From this point, we will always assume $j\leq l$ in the expression $\alpha_{ijl}^p$ (that is, there should be at least one $e$ in the product) unless specified otherwise. Diagrammatically, in the Kauffman tangle algebra on $n$ strands, the product $\a{ijl}{0}$ may be visualised as a `tangle diagram' with $n$ points on the top and bottom row such that the $i^\mathrm{th}$ and $(j+1)^\mathrm{th}$ are joined by a horizontal strand in the top row. The rest of the diagram consists of vertical strands, which cross over this horizontal strand but not each other, and a horizontal strand joining the $l^\mathrm{th}$ and $(l+1)^\mathrm{th}$ points in the bottom row. We illustrate this roughly in Figure~\ref{fig:onealpha} below.
\vspace{-2mm}
 \begin{figure}[h!]
 \begin{center}
   \includegraphics[scale=0.9]{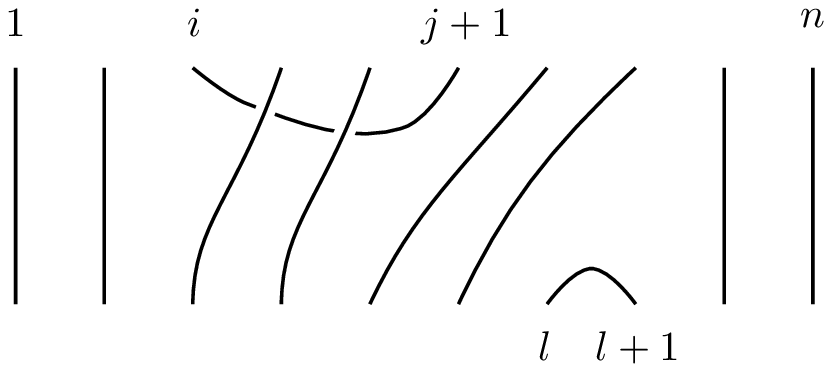}
   \caption{A diagrammatic interpretation of $\a{ijl}{0} = X_i \o X_{j-1} e_j \o e_l$ as a tangle on $n$ strands.} \label{fig:onealpha}
\end{center}
\end{figure}
\\ \vskip -1cm Thus one should think of the set in Theorem \ref{thm:span} as an ``inflation'' of a basis of $\mathfrak{h}_{n-2f,k}$, for each $f = 1,2,\ldots, \lfloor \frac{n}{2} \rfloor$, by ``dangles" with $f$ horizontal arcs (as seen in Xi \cite{X00}), one for each $\alpha$ chain occurring, with powers of $Y_i'$ elements attached. This will be further illustrated by Figure~\ref{fig:inflatedbasis} of Section \ref{sect:basis}. Using this pictorial visualisation, one may then use a straightforward calculation to show that the spanning set of Theorem \ref{thm:span} has cardinality $k^n (2n-1)!!$. Our eventual goal is to prove that this spanning set is in fact a basis of $\bnk$.

The following lemma essentially states that left multiplication of an $\a{mjl}{p}$ chain by a generator of $\bl{l+1}$ yields another $\alpha$ chain multiplied by `residue' terms in the smaller $\bl{l-1}$ subalgebra. Specifically, it helps us to prove that the $R$-submodule spanned by $\{\a{ijl}{p} \bl{l-1}\}$ is a left ideal of $\bl{l+1}$, in particular when $p$ is restricted to be within a certain range of $k$ consecutive integers.
\begin{lemma} \label{lemma:leftideal}
  For $\gamma \in \{X,e\}$, $m \leq l$ and $p \in \mathbb{Z}$,
  \[
  \gamma_m \alpha_{ijl}^p \in \la{\a{i'j'l}{p'} \bl{l-1} \,\Big|\, i' \geq \mathrm{min}(i,m), |p'| \leq |p| \text{ and } p'p \geq 0 \text{ unless } i'=m}.
  \]
  In fact, the only case in which $p'p < 0$ occurs is the case $X_i \alpha_{ijl}^p$, where $i \leq j \leq l$ and $p \in \mathbb{Z}$.
\end{lemma}

\begin{proof}
  Let $p$ be any integer and fix $m$, $i$, $j$ and $l$. \\
Henceforth, let $T := \la{\a{i'j'l}{p'} \bl{l-1} \,\big|\, i' \geq \mathrm{min}(i,m), |p'| \leq |p| \text{ and } p'p \geq 0 \text{ unless } i'=m}$.
In all of the following calculations, it is straightforward in each case to check that the resulting elements satisfy the minimality condition required to be a member of $T$.

\noi  \textbf{The action of $e_m$}. \\
  The action of $e_m$ on $\a{ijl}{p}$ falls into the following four cases:
\begin{enumerate}
  \item[(1)] $e_m \cdot \y{m}{p} X_m \ldots e_j \ldots e_l$, where $m < j \leq l$,\vspace{2.7mm}
  \item[(2)] $e_m \cdot \y{m}{p} e_m \ldots e_l$, where $m = j \leq l$,\vspace{2.7mm}
  \item[(3)] $e_m \cdot \a{m+1,j,l}{p}$, where $m+1 \leq j \leq l$, \vspace{2.7mm}
  \item[(4a)] $e_m \cdot \a{ijl}{p}$, where $m < i-1$ and $i \leq j \leq l$, \vspace{2.7mm}
  \item[(4b)] $e_m \cdot \a{ijl}{p}$, where $m > i$ and $i \leq j \leq l$.
\end{enumerate}
\textbf{(1).} By Lemma \ref{lemma:eype} (III),
\begin{align*}
  e_m \y{m}{p} X_m \ldots e_j \ldots e_l &\stackrel{\hphantom{(1)}}{\in} \la{e_m Y^{s_1} \ldots \y{m-1}{s_{m-1}}\y{m}{s_m} X_{m+1} \ldots e_j \ldots e_l \,\big|\, |s_m| \leq |p|\,} \\
&\stackrel{\hphantom{(1)}}{\in} \la{e_m  X_{m+1} \ldots e_j \ldots e_l \left( Y^{s_1} \ldots \y{m-1}{s_{m-1}}\y{m}{s_m} \right) \,\big|\, |s_m| \leq |p|\,} \\
&\stackrel{(\ref{eqn:xxe})}{\in} \la{e_m \ldots e_l \left(X_{j-2}^{-1} \ldots X_{m+1}^{-1} X_m^{-1} \y{m}{s_m} \y{m-1}{s_{m-1}} \ldots Y^{s_1}\right) \,\big|\, |s_m| \leq |p|\,}\!.
\end{align*}
Here $m\leq j-1\leq l-1$, so the term $X_{j-2}^{-1} \ldots X_{m+1}^{-1} X_m^{-1} \y{m}{s_m} \y{m-1}{s_{m-1}} \ldots Y^{s_1}$ in the brackets above is in $\bl{l-1}$. 
Hence, $e_m \y{m}{p} X_m \ldots e_j \ldots e_l \in \la{\a{mml}{0}\bl{l-1}} \subseteq T$.

\noi \textbf{(2).} By Lemma \ref{lemma:eype} (I),
\[
  e_m \y{m}{p} e_m \ldots e_j \ldots e_l \in \la{Y^{s_1} \y{2}{s_2} \ldots \y{m-1}{s_{m-1}} e_m \o e_l}
\subseteq \la{\a{mml}{0} \bl{l-1}} \subseteq T,
\]
since $Y^{s_1} \y{2}{s_2} \ldots \y{m-1}{s_{m-1}} \in \bl{l-1}$, as $m \leq l$ in this case. 

\noi \textbf{(3).} By equations (\ref{eqn:m9}), (\ref{eqn:prop1b}) and (\ref{eqn:xxe}),
\[
  e_m \a{m+1,j,l}{p} = e_m X_{m+1} \o e_j \o e_l \y{m}{-p} e_m \o e_l \left(X_{j-2}^{-1} \ldots X_{m+1}^{-1} X_m^{-1} \y{m}{-p} \right) \in T.
\]

\noi \textbf{(4).} We want to prove that $e_m \a{ijl}{p} \in T$, when $i \neq m$, $m+1$ and $i \leq j \leq l$.
This is separated into the following two cases. 

(a) If $m \leq i-2 \leq l-2$, then $e_m \in \bl{l-1}$ commutes past $\a{ijl}{p}$ so $e_m \a{ijl}{p} = \a{ijl}{p} e_m \in T$.

(b) On the other hand, if $m \geq i+1$ then: 

When $m < j$, using (\ref{eqn:xxe}) and the commuting relations gives the following
\begin{align*}
  e_m \a{ijl}{p} &\rel{eqn:xxe} e_m \o e_l \left(X_{j-2}^{-1} \ldots X_{m+1}^{-1} X_m^{-1}\right) \left(\y{i}{p} X_i \o X_{m-2} e_{m-1} \right),
\end{align*}
which is an element of $T$, since $m \leq l-1$ and $j-2 \leq l-2$.

When $m=j$,
\[
e_m \a{ijl}{p} \rel{eqn:exe} \l^{-1} \y{i}{p} X_i \o X_{m-2} e_m \o e_l \p e_m \o e_l (\l^{-1} \y{i}{p} X_i \o X_{m-2}),
\]
which is an element of $T$, since $m-2 = j-2 \leq l-2$.

When $m>j$,
\begin{align*}
  e_m \a{ijl}{p} &\p \y{i}{p} X_i \o \gamma_{m-2} e_m e_{m-1} e_m e_{m+1} \o e_l, \quad \text{where $\gamma$ could be $X$ or $e$,} \\
&\rel{eqn:eee} \y{i}{p} X_i \o \gamma_{m-2} e_m e_{m+1} \o e_l \,=\, e_m \o e_l \left( \y{i}{p} X_i \o \gamma_{m-2} \right) \in T.
\end{align*}

\noi We have now proved that $e_m \a{ijl}{p} \in T$, for all $m \leq l$, $i \leq j \leq l$, and $p \in \mathbb{Z}$.
\\
\\
\textbf{The action of $X_m$.}
\\
The action of $X_m$ on $\a{ijl}{p}$ falls into the following four cases:
\begin{enumerate}
\item[(A)] $X_m \cdot \a{m+1,j,l}{p}$, where $m+1 \leq j \leq l$, \vspace{2.9mm}
\item[(B)] $X_m \cdot \y{m}{p}e_m \o e_l$, where $m = j \leq l$, \vspace{2.9mm}
\item[(C)] $X_m \cdot \y{m}{p}X_m \o e_j \o e_l$, where $m < j \leq l-1$,  \vspace{2.9mm}
\item[(D1)] $X_m \cdot \a{ijl}{p}$, where $m < i-1$ and $i \leq j \leq l$, \vspace{2.9mm}
\item[(D2)] $X_m \cdot \a{ijl}{p}$, where $m > i$ and $i \leq j \leq l$.
\end{enumerate}
\vskip .3cm
\noi \textbf{(A).} When $p$ is a non-negative integer, using equations (\ref{eqn:m5}), (\ref{eqn:m9}), (\ref{eqn:prop1b}) and (\ref{eqn:prop1c}) gives 
\begin{eqnarray*}
X_m \a{m+1,j,l}{p} 
%\- \y{m}{p} X_m \gamma_{m+1} \o e_j \o e_l + \d \pum{s} \y{m}{p-s} \y{m+1}{s} \gamma_{m+1} \o e_j \o e_l \\
%&\phantom{= \y{m}{p} X_m \gamma_{m+1} \o e_j \o e_l}\,\, - \d \pum{s} \y{m}{p-s} e_m \y{m+1}{s}  \gamma_{m+1} \o e_j \o e_l \\
\= \y{m}{p} X_m \o e_j \o e_l + \d \pum{s} \y{m+1}{s} \gamma_{m+1} \o e_j \o e_l \left( \y{m}{p-s} \right) \\
&&\phantom{\y{m}{p} X_m \o e_j \o e_l} - \d \pum{s} \y{m}{p-s} e_m \o e_l \left(X_{j-2}^{-1} \ldots X_{m+1}^{-1} X_m^{-1} \y{m}{-s}\right),
\end{eqnarray*}
where $\gamma_{m+1}$ could be either $X_{m+1}$ or $e_{m+1}$.
Observe that if $1 \leq s \leq p$ then $0 \leq p-s \leq p-1$, hence $|s|, |p-s| \leq |p|$. Also, in this case, $m \leq l-1$ and $j-2 \leq l-2$ so the expressions in the brackets above are indeed elements of $\bl{l-1}$. Hence, for all $m+1 \leq j \leq l$, we have $X_m \cdot \a{m+1,j,l}{p} \in T$.

\noi Also, by equations (\ref{eqn:m6}), (\ref{eqn:m9}) and (\ref{eqn:prop1b}), 
\begin{eqnarray*}
X_m \a{m+1,j,l}{-p} 
\= \y{m}{-p} X_m \o e_j \o e_l - \d \pum{s} \y{m+1}{s-p} \gamma_{m+1} \o e_j \o e_l \ll{\y{m}{-s}} \\
&&\phantom{\y{m}{-p} X_m \o e_j \o e_l} + \d \pum{s} \y{m}{-s} e_m \o e_l \left(X_{j-2}^{-1} \ldots X_{m+1}^{-1} X_m^{-1} \y{m}{p-s}\right).
\end{eqnarray*}
Observe that when $1 \leq s \leq p$, we have that $1-p \leq s-p \leq 0$ so certainly $|\!-\! s\kern1.1pt| = |s|$ and $|s-p| \leq |p|$. Again, since $m \leq l-1$ and $j-2 \leq l-2$, the expressions in the brackets above are indeed elements of $\bl{l-1}$. Hence, for all $m+1 \leq j \leq l$, $X_m \cdot \a{m+1,j,l}{-p} \in T$.

\noi \textbf{(B).} By Lemma \ref{lemma:eype}(II) and equation (\ref{eqn:prop1b}),
\begin{align*}
  X_m \y{m}{p}e_m \o e_l &\in \la{\a{mml}{s_m} \left(Y^{s_1} \o \y{m-1}{s_{m-1}} \right) \,\big|\, |s_m| \leq |p|\,} \subseteq T, \quad \text{since $m \leq l$}.
\end{align*}

\noi \textbf{(C).} 
When $p$ is a non-negative integer, using equations (\ref{eqn:m11}), (\ref{eqn:m9}), (\ref{eqn:prop1b}) and (\ref{eqn:prop1c}),
\begin{align*}
X_m & \y{m}{p}X_m \o e_j \o e_l \\
&\rel{eqn:xxe} \y{m+1}{p} X_{m+1} \ldots e_j \o e_l - \d \sum_{s=1}^{p-1} X_m \y{m+1}{p-s} X_{m+1} \ldots e_j \o e_l \left( \y{m}{s} \right) \\
&\phantom{\p \y{m+1}{p} X_{m+1} \ldots e_j \o e_l}\,\, + \d \sum_{s=1}^{p-1} X_m \y{m}{s} e_m \o e_l \left( X_{j-2}^{-1} \o X_{m+1}^{-1} X_m^{-1} \y{m}{s-p} \right).
\end{align*}

\noi The first term in the above equation is $\a{m+1,j,l}{p} \in T$. In the second summation term, we have elements of the form $X_m \y{m+1}{u} X_{m+1} \ldots e_j \o e_l \bl{l-1}$, where $1 \leq u \leq p-1$. By case (A) above and since $|u| \leq |p|$, we know therefore the second term is in $T$. Moreover, by case (B), $X_m \y{m}{s} e_m \o e_l \in T$ for all $1\leq s \leq p-1$, so the third term is also in $T$. \\
Hence, for all $m \leq l-1$, we have $X_m \y{m}{p}X_m \o e_j \o e_l \in T$.

\noi Similarly, using equations (\ref{eqn:m13}), (\ref{eqn:m9}) and (\ref{eqn:prop1b}),
\begin{align*}
X_m & \y{m}{-p}X_m \o e_j \o e_l \\
\- \y{m+1}{-p} X_{m+1} \ldots e_j \o e_l + \d \sum_{s=0}^{p} X_m \y{m}{s-p} X_{m+1} \ldots e_j \o e_l \left( \y{m}{-s} \right) \\
&\phantom{= \y{m+1}{-p} X_{m+1} \ldots e_j \o e_l}\,\, + \d \sum_{s=0}^{p} X_m \y{m}{-s} e_m \o e_l \left( X_{j-2}^{-1} \o X_{m+1}^{-1} X_m^{-1} \y{m}{p-s} \right).
\end{align*}

\noi The first term in the above equation is $\a{m+1,j,l}{-p} \in T$. The second summation term involves elements of the form $X_m \y{m+1}{u} X_{m+1} \ldots e_j \o e_l \bl{l-1}$, where $-p \leq u \leq 0$. Thus the 2nd term is in $T$, by case (A) above and since $|u| \leq |p|$. Moreover, by case (B), $X_m \y{m}{-s} e_m \o e_l \in T$ for all $0 \leq s \leq p$, so the third term is also in $T$. \\
Hence, for all $m \leq l-1$, $X_m \y{m}{-p}X_m \o e_j \o e_l \in T$. We have now proved that, whether
$p$ is positive or negative, $X_m \a{mjl}{p} \in T$.

\noi \textbf{(D).} We want to prove that $X_m \a{ijl}{p} \in T$, when $i \neq m$, $m+1$, $i \leq j \leq l$ and $p$ is any integer. This is separated into the following two cases.

(D1). If $m \leq i-2 \leq l-2$, then $X_m \in \bl{l-1}$ commutes past $\a{ijl}{p}$. Hence $X_m \a{ijl}{p} = \a{ijl}{p} X_m \in T$.

(D2). On the other hand, if $m \geq i+1$ then again we have the following three cases to consider:

When $m <j \leq l$,
\[
  X_m \a{ijl}{p} \stack{(\ref{eqn:braid1}),(\ref{eqn:braid2})} \y{i}{p} X_i \o X_{m-2} X_{m-1} X_m \o e_j \o e_l (X_{m-1}).
\]
This is an element of $T$ as $X_{m-1} \in \bl{l-1}$, since $m < j \leq l$ in this case.

When $m = j \leq l$,
we have $X_m \a{ijl}{p} \rel{eqn:xxe} \y{i}{p} X_i \o X_{m-2} e_{m-1} e_m \o e_l = \a{i,j-1,l}{p} \in T$.

When $m>j$,
\begin{eqnarray*}
  X_m \a{ijl}{p} \= \y{i}{p} X_i \o \gamma_{m-2} X_m e_{m-1} e_m e_{m+1} \o e_l, \quad \text{where $\gamma$ could be $X$ or $e$,} \\
&\stack{(\ref{eqn:xxe}),(\ref{eqn:ix})}& \y{i}{p} X_i \o \gamma_{m-2} X_{m-1} e_m e_{m+1} \o e_l - \d \y{i}{p} X_i \o \gamma_{m-2} e_m e_{m+1} \o e_l \\
&&\phantom{\y{i}{p} X_i \o \gamma_{m-2} X_{m-1} e_m e_{m+1} \o e_l} + \d \y{i}{p} X_i \o \gamma_{m-2} e_{m-1} e_m e_{m+1} \o e_l \\
&\rel{eqn:braid2}& \y{i}{p} X_i \o \gamma_{m-2} X_{m-1} e_m e_{m+1} \o e_l - \d e_m e_{m+1} \o e_l \left( \y{i}{p} X_i \o \gamma_{m-2} \right) \\
&&\phantom{\y{i}{p} X_i \o \gamma_{m-2} X_{m-1} e_m e_{m+1} \o e_l} + \d \y{i}{p} X_i \o \gamma_{m-2} e_{m-1} e_m e_{m+1} \o e_l.
\end{eqnarray*}
Observe that $\y{i}{p} X_i \o \gamma_{m-2} \in \bl{l-1}$, as $m\leq l$, in this case. 

\noi Furthermore, if $m-2 \geq j$, then $\gamma_{m-2} = e_{m-2}$ and 
\begin{align*}
X_m \a{ijl}{p} \- \y{i}{p} X_i \o e_{m-2} e_{m-1} e_m \o e_l \left( X_{m-2}^{-1} \right) - \d e_m e_{m+1} \o e_l \left( \y{i}{p} X_i \o e_{m-2} \right)\\ 
&\phantom{= \y{i}{p} X_i \o e_{m-2} e_{m-1} e_m \o e_l \left( X_{m-2}^{-1} \right)}\,\, + \d \y{i}{p} X_i \o e_{m-2} e_{m-1} e_m \o e_l.
\end{align*}
Otherwise, if $m-1 = j$, then $\gamma_{m-2} = X_{m-2}$ and
\begin{align*}
X_m \a{ijl}{p} \- \y{i}{p} X_i \o X_{m-2} e_{m-1} e_m \o e_l - \d e_m e_{m+1} \o e_l \left( \y{i}{p} X_i \o X_{m-2} \right) \\ 
&\phantom{= \y{i}{p} X_i \o X_{m-2} e_{m-1} e_m \o e_l}\,\, + \d \y{i}{p} X_i \o X_{m-2} e_{m-1} e_m \o e_l.
\end{align*}
We have now proved that for all $m \leq l$ and $i \leq j \leq l$ and $p \in \mathbb{Z}$, $X_m \a{ijl}{p} \in T$.
\end{proof}

The following lemma says, for a fixed $l$, the $R$-span of all $\a{ijl}{p} \bl{l-1}$ is a \emph{left ideal} of $\bl{l+1}$, when $p$ lies in a range of $k$ consecutive integers.

\begin{lemma} \label{lemma:angels}
  Fix some $l$. Suppose $0 \leq K < k$ and let \[P = \{-K, -K+1, \o, k-K-1\}.\] The $R$-submodule 
\[
L := \la{\a{ijl}{p} \bl{l-1} \,\big| \, p \in P}
\]
is a left ideal of $\bl{l+1}$.
\end{lemma}

\begin{proof}
We want to prove that $L$ is invariant under left multiplication by the generators of $\bl{l+1}$, namely $Y, X_1, \o, X_l, e_1, \o, e_l$.

If $i > 1$, $Y$ commutes with $\a{ijl}{p}$. Otherwise, when $i = 1$, $Y \a{ijl}{p} = \a{ijl}{p+1}$ so clearly $L$ is invariant under left multiplication by $Y^{\pm 1}$, due to the $k^{\text{th}}$ order relation on $Y$. We will show by induction on $m \leq l$ that $L$ is invariant under $X_m$ and $e_m$.

Suppose $L$ is invariant under $X_{m'}$ and $e_{m'}$ for $m'<m$. Note that when $m=1$, this assumption is vacuous. 
Then in particular, $L$ is invariant under $X_{m'}^{-1} = X_{m'} - \d + \d e_{m'}$  for all $m'<m$. Moreover, this implies $L$ is invariant under $\y{m}{\pm 1}$. Thus, for all $p' \in \mathbb{Z}$,
\begin{equation} \label{eqn:angel}
\a{m,j',l}{p'} = \y{m}{p'} \a{m,j',l}{0} \in L.
\end{equation}
For $\gamma_m \in \{X_m,e_m\}$ and $p \in P$, Lemma \ref{lemma:leftideal} implies that
\begin{align*}
  \gamma_m \a{ijl}{p} &\in  \la{\a{i'j'l}{p'} \bl{l-1} \,\big|\, i' \geq \mathrm{min}(m,i), |p'| \leq |p| \text{ and } p'p \geq 0 \text{ unless } i'=m}\\
&\subseteq \la{\a{i'j'l}{p'} \bl{l-1} \,\big|\, |p'| \leq |p| \text{ and } p'p \geq 0}+ \la{\a{mj'l}{p} \bl{l-1}}.
\end{align*}
The first set lies in $L$, as if $|p'| \leq |p|$ and $p'p \geq 0$, then $p \in P$ implies $p' \in P$. By (\ref{eqn:angel}) above, $\la{\a{mj'l}{p} \bl{l-1}} \subseteq L$. Thus $\gamma_m \a{ijl}{p} \in L$, whence $\gamma_m L \subseteq L$ and $L$ is a left ideal of $\bl{l+1}$.
\end{proof}

Now we fix $K := \left\lfloor \frac{k-1}{2} \right\rfloor$. %When $k$ is odd, $K = \frac{k-1}{2}$ and when $k$ is even, $K = \frac{k-2}{2}$.
The range $P$ in Lemma \ref{lemma:angels} becomes
\[P = \left\{ -\left\lfloor \frac{k-1}{2} \right\rfloor, \o, k - \left\lfloor \frac{k-1}{2} \right\rfloor - 1\right\} = \left\{ \floor{\frac{k}{2}} - (k-1), \o, \floor{\frac{k}{2}} \right\}.\]
For $k$ odd, $P = \{ -K, \o, K \}$ and for $k$ even, $P = \{ -K, \o, K+1 \}$.

We are now almost ready to prove Theorem \ref{thm:span}. A standard way to show that a set which contains the identity element spans the entire algebra is to show it spans a left ideal of the algebra or, equivalently, show that its span is invariant under left multiplication by the generators of the algebra. We demonstrate this in stages, almost as if by pushing through one $\alpha$ chain at a time. With the previous lemma in mind, we observe that `pushing' a generator through each $\alpha$ chain may distort the `ordering' of the $\alpha$ chains (the $i_1 > i_2 > \ldots > i_m$ requirement in the statement of Theorem \ref{thm:span}). Motivated by this, we first prove the following Lemma. 

\begin{lemma} \label{lemma:beer}
  If $i \leq g$ and $p,r \in P$,
\[
\a{i,j,l}{p} \a{g,h,l-2}{r} \in \la{\a{i',j',l}{p'} \a{g',h',l-2}{r'} \bl{l-3} \,|\, i' > i \text{ and } p',r' \in P}.
\]
\end{lemma}
\begin{proof}
  Observe that, by Lemma \ref{lemma:angels}, $L = \la{\a{g'h'l-2}{r'} \bl{l-3} | r' \in P}$ is a left ideal of $\bl{l-1}$ therefore it suffices to prove that, for all $i \leq g$ and $p,r \in P$,
\begin{align*}
\a{ijl}{p} \a{g,h,l-2}{r} 
&\in \la{\a{i',j',l}{p'} \bl{l-1} \a{g',h',l-2}{r'} \bl{l-3} \,|\, i' > i \text{ and } p',r' \in P}.
\end{align*}

\noi Let us denote $\la{\a{i',j',l}{p'} \bl{l-1} \a{g',h',l-2}{r'} \bl{l-3} \,\big|\, i' > i \text{ and } p',r' \in P}$ by $S$.

\noi If $g \geq j$, then 
$ e_j \o e_l \ll{\y{g}{r} X_g \o e_h \o e_{l-2}} \stack{(\ref{eqn:ep2}),(\ref{eqn:eee})} \y{g+2}{r} X_{g+2} \o e_{h+2} \o e_l e_j \o e_{l-2}$.
Thus, using the commuting relations (\ref{eqn:braid2}), (\ref{eqn:xecomm}) and equation (\ref{eqn:prop1b}),
\[
  \a{ijl}{p} \a{g,h,l-2}{r} =  \y{g+2}{r} X_{g+2} \o e_{h+2} \o e_l \y{i}{p} X_i \o X_{j-1} e_j \o e_{l-2}  = \a{g+2,h+2,l}{r} \a{i,j,l-2}{p}.
\]
Note that $g+2 > j \geq i$ in this case.
Hence, when $g \geq j$, we have $\a{ijl}{p} \a{g,h,l-2}{r} \in  S$. 

\noi Now suppose on the contrary $g < j$. 
When $r$ is non-negative, we have the following:

\begin{eqnarray*}
  \a{ijl}{p} \a{g,h,l-2}{r}
&\rel{eqn:prop1b}& \y{i}{p} X_i \o X_g \y{g}{r} \a{g+1,j,l}{0} \a{g,h,l-2}{0} \\
&\rel{eqn:magic}& \y{i}{p} X_i \o X_{g-1} \ll{\y{g+1}{r} X_g} \a{g+1,j,l}{0} \a{g,h,l-2}{0} \\
&& - \,\, \d \sum_{s=1}^{r} \y{i}{p} X_i \o X_{g-1} \ll{\y{g+1}{s} \y{g}{r-s}} \a{g+1,j,l}{0} \a{g,h,l-2}{0} \\
&& + \,\, \d \sum_{s=1}^{r} \y{i}{p} X_i \o X_{g-1} \ll{\y{g+1}{s} e_g \y{g}{r-s}} \a{g+1,j,l}{0} \a{g,h,l-2}{0} \\
&\rel{eqn:prop1b}& \y{g+1}{r} \a{ijl}{p} \a{g,h,l-2}{0} \\
&& - \,\, \d \sum_{s=1}^{r} \a{g+1,j,l}{s} \y{i}{p} X_i \o X_{g-1} \a{g,h,l-2}{r-s} \\
&& + \,\, \d \sum_{s=1}^{r} \y{i}{p} X_i \o X_{g-1} \y{g+1}{s} e_g \a{g+1,j,l}{0} \a{g,h,l-2}{r-s}.
\end{eqnarray*}

\noi Observe that if $r \in P$ is non-negative, then because $1 \leq s \leq r$, it is clear that $s \in P$ and $r-s \leq r-1 \leq K$, hence $r-s \in P$ and $|r-s| \leq K$.

\noi On the other hand,
\begin{eqnarray*}
  \a{ijl}{p} \a{g,h,l-2}{-r} &\rel{eqn:m2}& \y{i}{p} X_i \o X_{g-1} \ll{\y{g+1}{-r} X_g} \a{g+1,j,l}{0} \a{g,h,l-2}{0} \\
&& + \,\, \d \sum_{s=1}^{r} \y{i}{p} X_i \o X_{g-1} \ll{\y{g+1}{s-r} \y{g}{-s}} \a{g+1,j,l}{0} \a{g,h,l-2}{0} \\
&& - \,\,  \d \sum_{s=1}^{r} \y{i}{p} X_i \o X_{g-1} \ll{\y{g+1}{s-r} e_g \y{g}{-s}} \a{g+1,j,l}{0} \a{g,h,l-2}{0} \\
&\rel{eqn:prop1b}& \y{g+1}{-r} \a{ijl}{p} \a{g,h,l-2}{0} \\
&& + \,\, \d \sum_{s=1}^{r} \a{g+1,j,l}{s-r} \y{i}{p} X_i \o X_{g-1} \a{g,h,l-2}{-s} \\
&& - \,\, \d \sum_{s=1}^{r} \y{i}{p} X_i \o X_{g-1} \y{g+1}{s-r} e_g \a{g+1,j,l}{0} \a{g,h,l-2}{-s}.
\end{eqnarray*}

If $-r \in P$, this means $-r \in \{-K, \o, -1\}$. So $1 \leq s \leq r$ implies $-s \in P$. Moreover, $|\!-\! s\kern1.1pt| \leq K$ and $s-r \in P$. To summarise, whether $r$ is positive \emph{or} negative,
\begin{align}
&  \a{ijl}{p} \a{g,h,l-2}{r} \in \y{g+1}{r} \a{ijl}{p} \a{g,h,l-2}{0} + \la{\a{g+1,j,l}{s} \y{i}{p} X_i \o X_{g-1} \a{g,h,l-2}{r-s} \,\big|\, s \in P, |r-s| \leq K} \nonumber \\
&\qquad \qquad \qquad \qquad \quad + \la{\y{i}{p} X_i \o X_{g-1} \y{g+1}{s} e_g \a{g+1,j,l}{0}\a{g,h,l-2}{r-s} \,\big|\, s \in P, |r-s| \leq K}. \label{eqn:twochains}
\end{align}
We now deal with each term of (\ref{eqn:twochains}) separately. \\
The first term is $\y{g+1}{r} \a{ijl}{p} \a{g,h,l-2}{0} = \y{g+1}{r} \y{i}{p} \a{ijl}{0} \a{g,h,l-2}{0}$.

If $h \leq j-2$ (so $i \leq g\leq h \leq j-2$), then 
\begin{align*}
  \a{ijl}{0} \a{g,h,l-2}{0} &\rel{eqn:XXp1} \ll{X_{g+1} \o X_h e_{h+1} \o e_{j-2}} \ll{X_i \o X_{j-2} X_{j-1} e_j \o e_l} \ll{e_{j-2} \o e_{l-2}} \\
&\kern.18em\rel{eqn:xxe} \ll{X_{g+1} \o X_h e_{h+1} \o e_{j-2} e_{j-1} \o e_l} \ll{X_i \o X_{j-3} e_{j-2} \o e_l} \\
&\h \a{g+1,h+1,l}{0} \a{i,j-2,l-2}{0}.
\end{align*}
As $i < g+1$, $\y{i}{p}$ commutes with $\a{g+1,h+1,l}{0}$, by equation (\ref{eqn:prop1b}).
Thus we have shown that $\y{g+1}{r} \a{ijl}{p} \a{g,h,l-2}{0} =  \a{g+1,h+1,l}{r} \a{i,j-2,l-2}{p}$ is an element of $S$, when $h \leq j-2$.

Now suppose $h \geq j-1$. Then, using equation (\ref{eqn:XXp1}) for $h \geq j$ and (\ref{eqn:ix}), we have the following:
\begin{eqnarray*}
  \a{ijl}{0} \a{g,h,l-2}{0} 
  &=& \ll{X_{g+1} \o X_{j-1}} \ll{X_i \o X_g \o X_{j-1}^{-1} e_j \o e_h \o e_l} \ll{X_{j-1} \o e_h \o e_{l-2}} \\
&& + \, \d \ll{X_{g+1} \o X_{j-1}} \ll{X_i \o X_g \o X_{j-2} e_j \o e_h \o e_l} \ll{X_{j-1} \o e_h \o e_{l-2}} \\
&& - \, \d \ll{X_{g+1} \o X_{j-1}} \ll{X_i \o X_g \o X_{j-2} e_{j-1} e_j \o e_h \o e_l} \ll{X_{j-1} \o e_h \o e_{l-2}} \\
&\stack{(\ref{eqn:xxe}),(\ref{eqn:ep2})}& \ll{X_{g+1} \o X_j} \ll{X_i \o X_g \o X_{j-2} e_{j-1} e_j \o e_h \o e_l} \ll{X_{j-1} \o e_h \o e_{l-2}} \\
&& + \, \d \ll{X_{g+1} \o X_{j-1} e_j \o e_l} \ll{X_i \o X_{j-2} X_{j-1} \o e_h \o e_{l-2}} \\
&& - \, \d \ll{X_{g+1} \o X_{j-1}} \ll{X_{j+1} \o e_{h+2} \o e_l} \ll{X_i \o X_{j-2} e_{j-1} \o e_h \o e_l} \\
&\stack{(\ref{eqn:ep2}),(\ref{eqn:eee})}& \ll{X_{g+1} \o X_j X_{j+1} \o e_l} \ll{X_i \o X_{j-2} e_{j-1} \o e_h \o e_l} \\
&& + \, \d \ll{X_{g+1} \o X_{j-1} e_j \o e_l} \ll{X_i \o X_{j-2} X_{j-1} \o e_h \o e_{l-2}} \\
&& - \, \d {X_{j+1} \o e_{h+2} \o e_l} \ll{X_{g+1} \o X_{j-1}} \ll{X_i \o X_{j-2} e_{j-1} \o e_h \o e_{l-2}} \\
&\rel{eqn:xxe}& \a{g+1,h+2,l}{0} \a{i,j-1,l-2}{0} + \d \a{g+1,j,l}{0} \a{i,h,l-2}{0} - \d \a{j+1,h+2,l}{0} \a{i,g,l-2}{0}.
\end{eqnarray*}
Therefore, when $h \geq j-1$,
\[
\y{g+1}{r} \a{ijl}{p} \a{g,h,l-2}{0} = \a{g+1,h+2,l}{r} \a{i,j-1,l-2}{p} + \d \, \a{g+1,j,l}{r} \a{i,h,l-2}{p} - \d \,\y{g+1}{r} \a{j+1,h+2,l}{0} \a{i,g,l-2}{p}.
\]
However, since $g < j$,  $\y{g+1}{r}$ commutes with $\a{j+1,h+2,l}{0}$. Moreover, as $g \leq l-2$, $\y{g+1}{r} \in \bl{l-1}$. Thus $\y{g+1}{r} \a{ijl}{p} \a{g,h,l-2}{0} \in S$.
So far we have proved that the first term of (\ref{eqn:twochains}) is a member of $S$, for all possibilities $i,j,g,h$ where $i \leq g < j$.

We now need to show $\la{\a{g+1,j,l}{s} \y{i}{p} X_i \o X_{g-1} \a{g,h,l-2}{r-s} \,\big|\, s \in P, |r-s| \leq K} \subseteq S$. But this follows immediately from the definition of $S$, as $r-s \in P$ and $\y{i}{p} X_i \o X_{g-1} \in \bl{l-1}$, since $i \leq g-1 \leq l-3$. \\
Finally, we now prove $\la{\y{i}{p} X_i \o X_{g-1} \y{g+1}{s} e_g \a{g+1,j,l}{0}\a{g,h,l-2}{r-s} \,\big|\, s \in P, |r-s| \leq K} \subseteq S$.  

\noi Let $\sigma := \y{i}{p} X_i \o X_{g-1} \y{g+1}{s} e_g \a{g+1,j,l}{0}\a{g,h,l-2}{r-s}$, where $s \in P$ and $|r-s| \leq K$.

\noi By equation (\ref{eqn:xxe}), $e_g \a{g+1,j,l}{0} = e_g \o e_l \ll{X_{j-2}^{-1} \o X_{g+1}^{-1} X_g^{-1}}$. Then
\[
\sigma = \y{i}{p} X_i X_{i+1} \o X_{g-1} \y{g+1}{s} \a{ggl}{0} \ll{X_{j-2}^{-1} \o X_{g+1}^{-1} X_g^{-1}} \a{g,h,l-2}{r-s}.
\] 
Lemma \ref{lemma:leftideal} implies that 
\[
X_g^{-1} \a{g,h,l-2}{r'} \in \la{\a{g'h'l-2}{p'} \bl{l-3} \,|\, g' \geq g, |p'| \leq |r'| \text{ and } p'r' \geq 0 \text{ unless } g' = g},
\]
where $r' = r-s$.
Observe that $p'$ and $r'$ could be of different signs, but as $|r-s| \leq K$, we know that $|p'| \leq K$.
%hence $p' \in P$ as well.
This allows us to apply Lemma \ref{lemma:leftideal} repeatedly to get
\[
X_{j-2}^{-1} \o X_{g+1}^{-1} X_g^{-1} \a{g,h,l-2}{r-s} \in \la{\a{g'h'l-2}{p'} \bl{l-3} \,|\, g' \geq g, |p'| \leq K}.
\]
Therefore
\begin{eqnarray*}
  \sigma &\in& \la{\y{i}{p} X_i X_{i+1} \o X_{g-1} \y{g+1}{s} \a{ggl}{0} \a{g'h'l-2}{p'} \bl{l-3} \,|\, g' \geq g, |p'| \leq K} \\
&\stack{(\ref{eqn:ep2}),(\ref{eqn:eee})}& \la{\y{i}{p} X_i X_{i+1} \o X_{g-1} \y{g+1}{s} \a{g'+2,h'+2,l}{p'} \a{g,g,l-2}{0} \bl{l-3} \,|\, g' \geq g, |p'| \leq K} \\
&\in& \la{\a{g'+2,h'+2,l}{p'} \y{i}{p} X_i X_{i+1} \o X_{g-1} \y{g+1}{s} \a{g,g,l-2}{0} \bl{l-3} \,|\, g' \geq g, |p'| \leq K}.
%&\stackrel{\phantom{(19)}}{\subseteq} \la{\a{g'+2,h'+2,l}{p'} \y{i}{p} X_i X_{i+1} \o X_{g-1} \y{g+1}{s} L \,|\, g' \geq g, |p'| \leq K}.
%&\stackrel{\phantom{(19)}}{\subseteq} \la{\a{g'+2,h'+2,l}{p'} L \,|\, g' \geq g, |p'| \leq K}.
\end{eqnarray*}
Now $i \leq g \leq l-2$, so $\y{i}{p} X_i X_{i+1} \o X_{g-1} \y{g+1}{s} \in \bl{l-1}$. 
Note that $g'+2 \geq g+2 > i$ and $p' \in P$. Thus $\sigma \in S$.

Therefore we have now proven that each term arising in (\ref{eqn:twochains}) is in S, for all $i,j,g,h$ where $i \leq g < j$.
This concludes the proof of Lemma \ref{lemma:beer}.
\end{proof}

Henceforth, we implicitly require $p \in P$ and $j \leq l$ whenever $\alpha_{ijl}^p$ is written, unless stated otherwise. For all $m \geq 0$ and $l \geq 2m$, let us define the following subsets of $\bnk$. Note that these are \emph{not} $R$-submodules.
\begin{eqnarray*}
  V_{l,m}^g &:=& \left\{ \a{i_1,j_1,l-1}{s_1} \a{i_2,j_2,l-3}{s_2} \o \a{i_m,j_m,l-2m+1}{s_m} \,\big|\,
		i_s\leq j_s\leq l-2s+1 \right\} \quad \text{and} \\
  V_{l,m} &:=& \left\{ \a{i_1,j_1,l-1}{s_1} \a{i_2,j_2,l-3}{s_2} \o \a{i_m,j_m,l-2m+1}{s_m} \,\big|\, i_1 > i_2 > \ldots > i_m \text{ and } i_s\leq j_s\leq l-2s+1\right\}.
\end{eqnarray*}
If $m > l/2$, we let $V_{l,m}^g = V_{l,m}$ be the empty set.

\noi 
If $U_1$ and $U_2$ are subsets of $\bnk$, let $U_1 U_2 := \la{ u_1u_2 | u_1 \in U_1, u_2 \in U_2}$; in other words, the $R$-span of the set $\{u_1u_2 | u_1 \in U_1, u_2 \in U_2\}$.

\begin{lemma} \label{lemma:Vgleftideal}
For all $l$ and $m$, $V_{l,m}^g\bl{l-2m}$ is a left ideal of $\bl{l}$.
\end{lemma}

\begin{proof}
We prove the statement by induction on $m$. When $m=0$, we have $V_{l,m}^g = \{1\}$ so $V_{l,m}^g\bl{l-2m}=\bl{l}$ and the statement then follows trivially.

Suppose that $m\geq1$ and assume the statement is true for $m-1$. Note that $l\geq2m\geq2$. By the definition of $V_{l,m}^g$, we have
\begin{align*}
    \bl{l}V_{l,m}^g\bl{l-2m}\-\langle\bl{l}\alpha_{ij,l-1}^pV_{l-2,m-1}^g\bl{l-2m}\rangle\\
        &\subseteq \langle\alpha_{ij,l-1}^p\bl{l-2}V_{l-2,m-1}^g\bl{l-2m}\rangle\text{,\quad by Lemma \ref{lemma:angels},}\\
        &\subseteq \langle\alpha_{ij,l-1}^pV_{l-2,m-1}^g\bl{l-2m}\rangle\text{,\quad by induction,}\\
        \-V_{l,m}^g\bl{l-2m},
\end{align*}
as required.
\end{proof}

\begin{lemma}\label{lemma:VgV}
For all $l$ and $m$, we have $V_{l,m}^g\bl{l-2m} = V_{l,m}\bl{l-2m}$.
\end{lemma}

\begin{proof}
By definition, $V_{l,m} \subseteq V_{l,m}^g$ hence $V_{l,m}^g\bl{l-2m}\supseteq V_{l,m}\bl{l-2m}$. It now remains to prove the reverse inclusion. We again proceed by induction on $m$. In the case $m=0$, the statement merely says $\bl{l} = \bl{l}$. Furthermore, the statement is clearly satisfied when $m=1$, as $V_{l,1}^g = V_{l,1}$.

Suppose $m\geq2$ and the statement is true for $m'<m$. Then, by the definition~of~$V_{l,m}^g$,
\[
V_{l,m}^g\bl{l-2m}=\langle\alpha_{ij,l-1}^pV_{l-2,m-1}^g\bl{l-2m}\rangle.
\]
It therefore suffices to show
\begin{equation}\label{eqn:order}
    \alpha_{ij,l-1}^pV_{l-2,m-1}^g\bl{l-2m}\subseteq V_{l,m}\bl{l-2m},
\end{equation}
for $1\leq i<l$. We will prove (\ref{eqn:order}) by descending induction on $i$. Suppose (\ref{eqn:order}) holds for all $i'$ such that $i<i'<l$. Observe that when $i=l-1$, the inductive hypothesis is vacuous. By induction on $m$, $V_{l-2,m-1}^g\bl{l-2m}=V_{l-2,m-1}\bl{l-2m}$, thus the LHS of (\ref{eqn:order}) is spanned by the set of elements of the form
\[
\alpha_{ij,l-1}^p\alpha_{i_2j_2l-3}^{s_2}\ldots\alpha_{i_mj_ml-2m+1}^{s_m}\bl{l-2m},
\]
where $i_2>i_3>\ldots>i_m$. If $i>i_2$, then we already have $i>i_2>\ldots>i_m$, so this is a subset of $V_{l,m}\bl{l-2m}$ by definition. On the 
other hand, if $i\leq i_2$ then
\begin{align*}
    \alpha_{ij,l-1}^p\alpha_{i_2j_2l-3}^{s_2}\ldots\alpha_{i_mj_ml-2m+1}^{s_m} \bl{l-2m} \hspace{-30mm}\\
        &\subseteq \alpha_{ij,l-1}^p\alpha_{i_2j_2l-3}^{s_2}V_{l-4,m-2}^g\bl{l-2m}\\
        &\subseteq \langle\alpha_{i'j'l-1}^{p'}\alpha_{g'h'l-3}^{r'}\bl{l-4}V_{l-4,m-2}^g\bl{l-2m}\mid i'>i\rangle\text{, \quad by Lemma \ref{lemma:beer}},\\
        &\subseteq \langle\alpha_{i'j'l-1}^{p'}\alpha_{g'h'l-3}^{r'}V_{l-4,m-2}^g\bl{l-2m}\mid i'>i\rangle\text{, \quad by Lemma \ref{lemma:Vgleftideal}}, \\
        &\subseteq \langle\alpha_{i'j'l-1}^{p'}V_{l-2,m-1}^g\bl{l-2m}\mid i'>i\rangle\\
        &\subseteq V_{l,m}\bl{l-2m}\text{, \quad by induction on }i.
\end{align*}
Thus (1) holds. Hence $V_{l,m}^g\bl{l-2m} = V_{l,m}\bl{l-2m}$.
\end{proof}

\noi Recall
\[
\blk/\blk e_{l-1}\blk\cong \ak{l}
\]
and $\pi_l:\blk \twoheadrightarrow \ak{l}$ is the corresponding projection. Recall $\overline{\X}_{l,k}$ was an arbitrary subset of $\blk$ mapping onto a basis $\X_{l,k}$ of $\ak{l}$ and $|\overline{\X}_{l,k}| = |\X_{l,k}|$. We can define an $R$-module homomorphism $\phi_l:\ak{l} \rightarrow \blk$ by sending each element of $\X_{l,k}$ to the corresponding element of $\overline{\X}_{l,k}$.
Thus $\pi_l\phi_l={\rm id}_{\ak{l}}$.
%Such maps exist because $\ak{l}$ is a free $R$-module.
Note that when $l=0$ or $1$, we have an isomorphism $\pi_l:\blk\rightarrow \ak{l}$, with inverse $\phi_l$. And, for $l \geq 2$,
\[
    \blk=\phi_l(\ak{l})+\blk e_{l-1}\blk.
\]
Thus
\begin{equation}\label{Bsplit}
    \bl{l}=\At{l}+\bl{l}e_{l-1}\bl{l}, \quad\text{ for }l\geq2,
\end{equation}
where $\At{l}$ is the image of $\phi_l(\ak{l})$ in $\bnk$.\label{defn:phiak} 
Also, \begin{equation}\label{Bsplit01}
    \At{0}=\bl{0}\qquad \text{and}\qquad \At{1}=\bl{1}.
\end{equation}
%We define
%\begin{equation}\label{Bsplit01}
%    \At{0}=\bl{0}\qquad \text{and}\qquad %\At{1}=\bl{1}.
%\end{equation}
%For any $l\geq0$, let $\{\chi^{(l)}_i\}$ ({\bf BETTER NOTATION???}) be any basis of $\ak{l}$, and let $\bar \chi^{(l)}_i$ denote the image of $\phi_l(\chi^{(l)}_i)$ in $B_n^k$. Clearly $\At{l}$ is spanned by $\{\bar \chi^{(l)}_i\}$, as it maps onto a basis of $\ak{l}$.

\noi Now let us define
\[
\overline V_{l,m} := \left\{ \a{i_1,j_1,l-2}{s_1} \a{i_2,j_2,l-4}{s_2} \o \a{i_m,j_m,l-2m}{s_m} \,\big|\, i_1 > i_2 > \ldots > i_m\;\text{ and } i_s\leq j_s\leq l-2s+1\right\}.
\]
Note that the $\a{ijl}{p}$ appearing in $\overline V_{l,m}$ need not satisfy $j\leq l$; in other words, the $\alpha$ chains need not contain any $e$'s. Also $V_{l,m}=\overline V_{l,m}E_{l,m}$, where
\[
E_{l,m}:=e_{l-1}e_{l-3}\ldots e_{l-2m+1}.
\]
In order to prove Theorem \ref{thm:span}, we need to show $\bnk$ is spanned by $V_{n,m}\widetilde{\X}_{n-2m,k} \overline V_{n,m}^*$. This will immediately follow as a corollary to the following result.

\begin{lemma}\label{Iprop}
Let $I_{l,m}=\overline V_{l,m}E_{l,m}\bl{l-2m}\overline V_{l,m}^*$.
\begin{enumerate}
\item[(a)] $I_{l,m}$ is a two-sided ideal of $\bl{l}$.
\item[(b)] For $l\geq2m+2$, we have
\[
    \overline V_{l,m}E_{l,m}\bl{l-2m}e_{l-2m-1}\bl{l-2m}\overline V_{l,m}^*\subseteq I_{l,m+1}.
\]
\item[(c)] For any fixed $M$, $I_{l,M} = \sum_{m \geq M} V_{l,m} \At{l-2m}\overline V_{l,m}^*$ and is spanned by elements of the form
\[
    \alpha_{i_1,j_1,l-1}^{s_1}\ldots\alpha_{i_m,j_m,l-2m+1}^{s_m} \chi(\alpha_{g_m,h_m,l-2m}^{t_m})^*\ldots(\alpha_{g_1,h_1,l-2}^{t_1})^*,
\]
where $m\geq M$, $i_1>i_2>\ldots>i_m$, $g_m<g_{m-1}<\ldots<g_1$, $i_p\leq j_p\leq l-2p+1$, $g_p\leq h_p\leq l-2p+1$, $s_p,t_p\in P$ and 
$\chi\in\widetilde{\X}_{l-2m,k}$.
%That is, $I_{l,M}$ is spanned by
%\[
%   \{a\bar{x}^{l-2m}_ib^*\mid m\geq M\text{ and }a,b\in V_{l,m}\}.
%\]
\end{enumerate}
\end{lemma}

\begin{proof}
\begin{enumerate}
\item[(a)]
    By Lemma \ref{lemma:VgV}, we have $I_{l,m}=V_{l,m}\bl{l-2m}\overline V_{l,m}^*=V_{l,m}^g\bl{l-2m}\overline V_{l,m}^*$. 
    Therefore $I_{l,m}$ is a 
    left ideal of $\bl{l}$ by Lemma \ref{lemma:Vgleftideal}. The subalgebra $\bl{l-2m}$ is preserved by $*$ and commutes with $E_{l,m}$, so 
    $I_{l,m}=I_{l,m}^*$ is also a right ideal.
\item[(b)]
Suppose $l\geq2m+2$. Since
\[
1=\alpha_{l-1,l-1,l-2}^0\alpha_{l-3,l-3,l-4}^0\ldots\alpha_{l-2m-1,l-2m-1,l-2m-2}^0\in\overline V_{l,m+1},
\]
we have
\[
E_{l,m}e_{l-2m-1}=E_{l,m+1}\in\overline V_{l,m+1}E_{l,m+1}\bl{l-2m-2}\overline V_{l,m+1}=I_{l,m+1}.
\]
But $I_{l,m+1}$ is a two-sided ideal in $\bl{l}$, so
\[
\overline V_{l,m}E_{l,m}\bl{l-2m}e_{l-2m-1}\bl{l-2m}\overline V_{l,m}^*
=\overline V_{l,m}\bl{l-2m}E_{l,m}e_{l-2m-1}\bl{l-2m}\overline V_{l,m}^*
\subseteq I_{l,m+1}. 
\]
\item[(c)]
    If $m \geq M$, then $l-2m \leq l-2M$ so the given elements are clearly contained in $I_{l,M}$. For a fixed $m$, they span the set
    \[
        \overline V_{l,m}E_{l,m}\At{l-2m}\overline V_{l,m}^*.
    \]
    It therefore suffices to prove that
    \begin{equation}\label{spanind}
        I_{l,M}\subseteq\sum_{m\geq M}\overline V_{l,m}E_{l,m}\At{l-2m}\overline V_{l,m}^*.
    \end{equation}
    We prove this statement by induction on $l-2M$. If $l-2M<2$ then  
\[
I_{l,M}=\overline V_{l,M}E_{l,M}\bl{l-2M}\overline V_{l,M}^*=V_{l,M}E_{l,M}\At{l-2M}\overline V_{l,M}^* \text{,\qquad by }(\ref{Bsplit01}).
\]
Now suppose $l-2M\geq2$ and assume $I_{l,M+1}\subseteq\sum_{m\geq M+1}\overline V_{l,m}E_{l,m}\At{l-2m}\overline V_{l,m}^*$. Then using (\ref{Bsplit}) and part (b) of this Lemma we have that
\begin{eqnarray*}
I_{l,M}
\= \overline V_{l,M}E_{l,M}\bl{l-2M}\overline V_{l,M}^* \\
&\rel{Bsplit}& \overline V_{l,M}E_{l,M}\At{l-2M}\overline V_{l,M}^*+\overline V_{l,M}E_{l,M}\bl{l-2M}e_{l-2M-1}\bl{l-2M}\overline V_{l,M}^*\\
&\subseteq& \overline V_{l,M}E_{l,M}\At{l-2M}\overline V_{l,M}^*+I_{l,M+1}\\
&\sstack{\text{ind. hypo.}}&\overline V_{l,M}E_{l,M}\At{l-2M}\overline V_{l,M}^*+\sum_{m\geq M+1}\overline V_{l,m}E_{l,m}\At{l-2m}\overline V_{l,m}^*\\
&=&\sum_{m\geq M}\overline V_{l,m}E_{l,m}\At{l-2m}\overline V_{l,m}^*,
\end{eqnarray*}
proving part (c).
\end{enumerate}
\end{proof}
In particular, $I_{n,0}=B_n^k$ by definition, so when $l = n$ and $m=0$, statement (c) of the previous Lemma implies that $\bnk$ is spanned by the set of elements
\[
\alpha_{i_1,j_1,n-1}^{s_1}\ldots\alpha_{i_m,j_m,n-2m+1}^{s_m} \chi^{(n-2m)}(\alpha_{g_m,h_m,n-2m}^{t_m})^*\ldots (\alpha_{g_1,h_1,n-2}^{t_1})^*,
\]
with conditions specified as above, giving Theorem \ref{thm:span}.

\section{The Admissibility Conditions} \label{sect:admissible}

In the previous section, we obtained a spanning set of $\bnk$ over an arbitrary ring $R$ and hence concluded the rank of $\bnk$ is \emph{at most} $k^n (2n-1)!!$. Before we can prove the linear independence of our spanning set, we must first focus our attention on the representation theory of the algebra $\bt(R)$. It is here that the notion of admissibility, as first introduced by H\"aring-Oldenburg~\cite{HO01}, arises. Essentially, it is a set of conditions on the parameters $A_0,\o, A_{k-1}, q_0, \o, q_{k-1}, q, \l$ in our ground ring $R$ which ensure the algebra $\bt(R)$ is $R$-free of the expected rank, namely $3k^2$.
 It turns out that, if $R$ is a ring with admissible parameters $A_0, \o, A_{k-1},q_0, \o, q_{k-1}, q, \l$ (see Definition \ref{defn:adm}) then the spanning set of Theorem \ref{thm:span} is actually a basis for \emph{general} $n$. 

We shall establish these admissibility conditions explicitly via a certain $\bt$-module $V$ of rank $k$. These results are contained in \cite{WY06}, in which the authors are able to use $V$ to then construct the regular representation of $\bt$ and provide an explicit basis of the algebra under the conditions of admissibility and the added assumption that $\delta$ is not a zero divisor.

It is non-trivial to show that there are any nonzero rings with admissible parameters; in other words, that the conditions we impose are consistent with each other. In Lemma \ref{lemma:Rsigma}, we demonstrate rings with admissible parameters and, in particular, construct a ``generic" ground ring $R_0$ with admissible parameters, in the sense that for every ring $R$ with admissible parameters there exists a unique map from $R_0$ to $R$ which respects the parameters (see Proposition \ref{prop:universalmap}). 

It is important to clarify the different notions of admissibility used in the literature. A comparison between the various definitions is offered at the end of the section. The proofs in this section are mostly the same as, if not a slight modification of, those in \cite{WY06} and \cite{Y07}, so we shall refer the reader to \cite{WY06} or \cite{Y07} for further details of proofs.

For this section, we simplify our notation by omitting the index $1$ of $X_1$ and $e_1$. Specifically, $\bt(R)$ is the unital associative $R$-algebra generated by $Y^{\pm 1}, X^{\pm 1}$ and~$e$ subject to the following relations.
\vspace{-1.5cm}
\begin{eqnarray}
Y^k \= \sum_{i=0}^{k-1} q_i Y^i \\
X - X^{-1} &=& \d(1-e) \label{eqn:xinv} \\
XYXY &=& YXYX \label{eqn:xyxy} \\
Xe &=& \l e \,\,\,=\,\,\, eX \label{eqn:xele}\\
YXYe &=& \lambda^{-1}e \,\,\,=\,\,\, eYXY \label{eqn:eyxy2} \\
eY^me &=& A_m e, \,\,\,\,\,\, \qquad \text{for }\, 0 \leq m \leq 
k-1. \label{eqn:eye2}
\end{eqnarray}

Recall, in Lemma \ref{lemma:eype}, we showed $X Y^p e = \l Y^{-p} e - \d \pum{s} Y^{p-2s} e + \d \pum{s} A_{p-s} Y^{-s} e$, for all $p \geq 0$. Using this and the $k^\mathrm{th}$ order relation on $Y$, it is straightforward to show the left ideal of $\bt$ generated by $e$ is the span of $\{Y^i e \mid 0 \leq i \leq k-1\}$. 
As a consequence of the results in Goodman and Hauschild \cite{GH06}, the set $\{Y^i e \mid i \in \mathbb{Z}\}$ is linearly independent in the affine BMW algebra and so it 
seems natural to expect that the set $\{Y^i e \mid 0 \leq i \leq k-1\}$ be linearly independent in the cyclotomic BMW algebra. If this 
were the case, the span of this set would be a $\bt$-module $V$ with the following properties:
\begin{equation}\begin{array}{l}\label{wantmodule}
\raisebox{0.2mm}{\text{\textbullet}} \,\,\, V\text{ has a basis }\{v_i\mid 0\leq i\leq k-1\};\\
\raisebox{0.2mm}{\text{\textbullet}} \,\,\, Yv_i=v_{i+1}\text{ for }0\leq i<k-1;\\
\raisebox{0.2mm}{\text{\textbullet}} \,\,\, YXYXv=v\text{ for }v\in V;\\
\raisebox{0.2mm}{\text{\textbullet}} \,\,\, Xv_0=\l v_0;\\
\raisebox{0.2mm}{\text{\textbullet}} \,\,\, ev_i=A_iv_0.
\end{array}\end{equation}
It is easy to see that these properties determine the action of X. The work of \cite{WY06} shows that the existence of such a module imposes additional restrictions 
on $A_0, \o, A_{k-1},q_0, \o, q_{k-1}, q, \l$. More precisely, write $k=2z-\e$ where $z:=\lceil k/2\rceil$ and $\e\in\{0,1\}$. Then (\ref{wantmodule}) implies that we must have
\[
	\beta=h_0=h_1=\ldots=h_{k-1}=0,
\]
where
\begin{eqnarray}
  \beta&:=&q_0\l-q_0^{-1}\l^{-1}+(1-\e)\d\label{eqn:alpha},\\
  h_0&:=&\l-\l^{-1}+\d(A_0-1)\label{eqn:h0}
\end{eqnarray}
and, for $l = 1, \o, k-1$,
\begin{equation} \label{eqn:hl}
  h_l:=\l^{-1}(q_l+q_0^{-1}q_{k-l})
    +\d\left[\sum_{r=1}^{k-l}q_{r+l}A_r
    -\sum_{i=\max(l+1,z)}^{\lfloor\frac{l+k}{2}\rfloor}\hspace{-5mm}q_{2i-l}
    +\sum_{i=\lceil\frac{l}{2}\rceil}^{\min(l,z-1)}\hspace{-5mm}q_{2i-l}\right].
\end{equation}
However, certain linear combinations of these elements are divisible by $\d$; a tedious calculation found in \cite{Y07} shows that, for $1\leq l \leq z - \e$,
\begin{equation}\label{hiprime}
    q_0^{-1}h_{k-l}-h_l+\beta q_0^{-1}q_l-h_0q_l=\d h_l',
\end{equation}
where
\begin{eqnarray}
    h_l'
        &:=&\sum_{r=1}^lq_0^{-1}q_{r+k-l}A_r-\sum_{r=0}^{k-l}q_{r+l}A_r\nonumber\\
        &&{}-\sum_{i=\lceil\frac{l}{2}\rceil}^{l-1}(q_0^{-1}q_{k-2i+l}+q_{2i-l})
            +\sum_{i=z}^{\lfloor\frac{l+k}{2}\rfloor}(q_0^{-1}q_{k-2i+l}+q_{2i-l}).\label{hiprime2}
\end{eqnarray}
It therefore seems sensible to work with rings $R$ in which we also require that $h_l'=0$. We aim to study the ``generic" ring $R_0$ (defined in Lemma \ref{lemma:Rsigma} below) in which all above relations hold. This 
will allow us to deduce results over an arbitrary such ring by proving them for $R_0$ first and then specialising. 
Before proceeding we first prove a simple lemma which will be used in a later proof to show $\d$ is not a zero divisor in certain rings.
\begin{lemma}\label{zerodiv}
Suppose a commutative ring $S$ contains elements $a$ and $b$, such that $a$ is not a zero divisor in $S$ and 
$b+aS$ is not a zero divisor in $S/aS$. Then $a+bS$ is not a zero divisor in $S/bS$.
\end{lemma}
\begin{proof}
Suppose $(a+bS)(x+bS)=0$ for some $x+bS\in S/bS$. Then $ax\in bS$, so $ax=by$ for some $y\in S$. Thus, as an element of $S/aS$, $(b+aS)(y+aS)=0$. This implies $y+aS=0$ since $b+aS$ is not a zero divisor in $S/aS$, by assumption. Hence, $y=az$ for some $z\in S$. Furthermore, $ax = by = azb$, so $x=zb$ since $a$ is not a zero divisor in $S$. Therefore $x+bS=0$ and $a+bS$ is not a zero divisor in $S/bS$.
\end{proof}

\noi It is easy to see that $\beta$ always factorises as $\beta_+\beta_-$, where if $k$ is odd,
\[
\beta_+=q_0\l-1 \qquad \text{and}\qquad \beta_-=q_0^{-1}\l^{-1}+1
\]
and when $k$ is even,
\[
    \beta_+=q_0\l-q^{-1} \qquad \text{and}\qquad \beta_-=qq_0^{-1}\l^{-1}+1.
\]
For convenience, we denote $\beta_0 := \beta$. At this point, we wish to remind the reader that for a subset $J \subseteq R$, we write $\la{J}_R$ to mean the \emph{ideal} generated by $J$ in $R$. The subscript may sometimes be omitted only if it is clear in the current context.

The following results exhibits rings with admissible parameters explicitly, in the sense of the definition following immediately after the Lemma. In particular, we introduce $R_0$, the ``generic" ring with admissible parameters. The results of the Lemma will play a key role in the arguments of Section \ref{sect:basis} for proving non-degeneracy of a trace map and thus the linearly independency of our spanning set.

\begin{lemma} \label{lemma:Rsigma}
Let
\[
\Omega := \Z[q^{\pm 1},\l^{\pm 1},q_0^{\pm 1},q_1,\ldots,q_{k-1},A_0,A_1,\ldots,A_{k-1}].
\]
For $\sigma\in\{0,+,-\}$, define $R_\sigma := \Omega/I_\sigma$ where
\[
    I_\sigma := \langle\beta_\sigma,h_0,h_1,\ldots,h_{z-\e},h_1',h_2',\ldots,h_{z-\e}'\rangle_{_{\Omega}}\subseteq\Omega.
\]
Then 
\begin{enumerate}
\item[(a)] the image of $\d$ is not a zero divisor in $R_\sigma$, for $\sigma\in\{0,+,-\}$;
\item[(b)] for $\sigma=\pm$,
\[
    R_\sigma[\d^{-1}]\cong\Z[q^{\pm 1},\l^{\pm 1},q_1,\ldots,q_{k-1}][\d^{-1}];
\]
\item[(c)] for $\sigma=\pm$, the ring $R_\sigma$ is an integral domain;
\item[(d)] $I_0=I_+\cap I_-$.
\end{enumerate}
\end{lemma}

\begin{proof}
\textbf{(a)} Since $\d= q - q^{-1} = q^{-1}(q-1)(q+1)$, to prove (a) it suffices to show that $q+\tau$ is not a zero divisor, for $\tau=\pm 1$. For $1\leq l\leq k-1$, let 
\begin{equation}\label{Bl}
    B_l :=\sum_{r=1}^{k-l}q_{r+l}A_r
        \,-\hspace{-5mm} \sum_{i=\max(l+1,z)}^{\lfloor\frac{l+k}{2}\rfloor}\hspace{-5mm}q_{2i-l}
        \,+\hspace{-3mm}\sum_{i=\lceil\frac{l}{2}\rceil}^{\min(l,z-1)}\hspace{-5mm}q_{2i-l}\quad \in\, \Omega.
\end{equation}
Then (\ref{eqn:hl}) says that
\begin{equation} \label{eqn:hBl}
    h_l=\l^{-1}(q_l+q_0^{-1}q_{k-l})+\d B_l, \quad \text{for } 1\leq l\leq k-1.
\end{equation}
Over $\Z[q^{\pm 1},\l^{\pm 1},q_0^{\pm 1},q_1,\ldots,q_{k-1}]$, the $B_l$ are related to 
the $A_l$ by an affine linear transformation; specifically, the column vector $(B_l)$, where $l = 1, \o,k-1$, is equal to a matrix $(s_{lr})_{l,r = 1}^{k-1}$ multiplied by the column vector $(A_l)$ plus a column vector of $q_i$'s. Moreover, $s_{lr} = 0$, unless $l + r \leq k$ and $s_{lr} = q_k = -1$, when $l+r = k$. Thus $(s_{lr})$ is triangular, with diagonal entries $q_k=-1$, 
so it is invertible. Therefore we may identify $\Omega$ with the polynomial ring
\begin{equation}\label{eqn:Lambdaalt}
\Omega=\Z[q^{\pm 1},\l^{\pm 1},q_0^{\pm 1},q_1,\ldots,q_{k-1},A_0,B_1,B_2,\ldots,B_{k-1}].
\end{equation}
Now, when $1 \leq l \leq z-1$, (\ref{hiprime2}) and (\ref{Bl}) implies that
\begin{align*}
  h_l' \- q_0^{-1} B_{k-l} \,\,+\,\, q_0^{-1} \sum_{i=\max(k-l+1,z)}^{\lfloor\frac{2k-l}{2}\rfloor}\hspace{-5mm}q_{2i-k+l}
         \,\,-\,\, q_0^{-1} \sum_{i=\lceil\frac{k-l}{2}\rceil}^{\min(k-l,z-1)}\hspace{-5mm}q_{2i-k+l} \\
\f - \,\sum_{r=0}^{k-l}q_{r+l}A_r \,\,-\sum_{i=\lceil\frac{l}{2}\rceil}^{l-1}(q_0^{-1}q_{k-2i+l}+q_{2i-l})
            \,\,+\sum_{i=z}^{\lfloor\frac{l+k}{2}\rfloor}(q_0^{-1}q_{k-2i+l}+q_{2i-l}) \\
&\rel{Bl}  q_0^{-1} B_{k-l} \,+\, q_0^{-1} \sum_{i=\max(k-l+1,z)}^{\lfloor\frac{2k-l}{2}\rfloor}\hspace{-5mm}q_{2i-k+l}
         \,-\, q_0^{-1} \sum_{i=\lceil\frac{k-l}{2}\rceil}^{\min(k-l,z-1)}\hspace{-5mm}q_{2i-k+l} - B_l - q_l A_0 \,\,-\!\!\!\!\!\!\! \sum_{i=\max(l+1,z)}^{\lfloor\frac{l+k}{2}\rfloor}\hspace{-5mm}q_{2i-l}\\
\f  
        +\!\!\!\sum_{i=\lceil\frac{l}{2}\rceil}^{\min(l,z-1)}\hspace{-3mm}q_{2i-l}\,\,\, -\sum_{i=\lceil\frac{l}{2}\rceil}^{l-1}(q_0^{-1}q_{k-2i+l}+q_{2i-l})\,\,
            +\sum_{i=z}^{\lfloor\frac{l+k}{2}\rfloor}(q_0^{-1}q_{k-2i+l}+q_{2i-l})
\end{align*}
Hence
\[
    h_l' \,\in\, q_0^{-1}B_{k-l}+\Z[q^{\pm 1},\l^{\pm 1},q_0^{\pm 1},q_1,\ldots,q_{k-1}, A_0,B_1,B_2,\ldots,B_{z-1}],
\]
for $1\leq l\leq z-1$. Thus
\[
    \Omega_1 := \Omega/\langle h_1',h_2',\ldots,h_{z-1}'\rangle_{_{\Omega}}
\cong\Z[q^{\pm 1},\l^{\pm 1},q_0^{\pm 1},q_1,\ldots,q_{k-1}, A_0,B_1,B_2,\ldots,B_{z-\e}].
\]
Indeed, if $1 \leq l \leq z-1$, then if $k$ is even (hence $\e = 0$ and $z+1 \leq k-l \leq k-1$), quotienting by the $h_l'$ expresses the elements $B_{k-1}, \o, B_{z+1}$, respectively, as elements of the ring $\Z[q^{\pm 1},\l^{\pm 1},q_0^{\pm 1},q_1,\ldots,q_{k-1}, A_0,B_1,B_2,\ldots, B_{z-1}]$; similarly, if $k$ is odd (hence $\e = 1$ and $z \leq k-l \leq k-1$), then $B_z, \o, B_{k-1}$ may be expressed in the quotient as elements of $\Z[q^{\pm 1},\l^{\pm 1},q_0^{\pm 1},q_1,\ldots,q_{k-1}, A_0,B_1,B_2,\ldots, B_{z-1}]$ .

In particular, $q+\tau$ is not a zero divisor in $\Omega_1$. By using Lemma \ref{zerodiv} recursively, we aim to show that $q+\tau$ does not become a zero divisor, as we quotient by further generators of $I_\sigma$. \\
If we set $q + \tau = 0$ then $\d = 0$, hence by (\ref{eqn:hl}), $h_i = \l^{-1}(q_i + q_0^{-1} q_{k-i})$, for any $i$. So for any $1\leq l\leq z-1$, we have that
\begin{eqnarray*}
    \Omega/\langle h_1',\ldots,h_{z-1}',h_1,\ldots,h_{l-1},q+\tau\rangle_{_{\Omega}}\hspace{-50mm}\\
        &\cong&\Omega_1/\langle(q_1+q_0^{-1}q_{k-1}),(q_2+q_0^{-1}q_{k-2})\ldots, (q_{l-1}+q_0^{-1}q_{k-l+1}),q+\tau\rangle_{_{\Omega_1}}\\
        &\cong&\Z[\l^{\pm 1},q_0^{\pm 1},q_1,\ldots,q_{k-l}, A_0,B_1,B_2,\ldots,B_{z-\e}].
\end{eqnarray*}
Certainly $h_l\equiv\l^{-1}(q_l+q_0^{-1}q_{k-l})$ is not a zero divisor in this ring, so repeated application of Lemma \ref{zerodiv} proves that $q+\tau$ is not a zero divisor in
\[
    \Omega_2 := \Omega_1/\la{h_1, \o, h_{z-1}}_{_{\Omega_1}} = \Omega/\langle h_1',h_2',\ldots,h_{z-1}',h_1,h_2,\ldots,h_{z-1}\rangle_{_{\Omega}}.
\]
Moreover, the above argument (with $l = z$) says that
\[
    \Omega_2/\langle q+\tau\rangle_{_{\Omega_2}}\cong\Z[\l^{\pm 1},q_0^{\pm 1},q_1,\ldots,q_{z-\e}, A_0,B_1,B_2,\ldots,B_{z-\e}].
\]
Observe that $\langle h_0,\beta_\sigma\rangle_{_{\Omega}}=\langle h_0',\beta_\sigma\rangle_{_{\Omega}}$, where
\[
    h_0':=h_0-q_0^{-1}\beta=\l^{-1}(q_0^{-2}-1)+\d(A_0-1-(1-\e)q_0^{-1}).
\]
Suppose first that $k$ is odd, then $R_\sigma = \Omega_2/\la{h_0,\beta_\sigma} = \Omega_2/\la{h_0',\beta_\sigma}$. Certainly, we know that $h_0'\equiv\l^{-1}(q_0^{-2}-1)$ is not 
a zero divisor in the polynomial ring $\Omega_2/\langle q+\tau\rangle_{_{\Omega_2}}$. Moreover,
\begin{eqnarray*}
    \Omega_2/\langle h_0',q+\tau\rangle_{_{\Omega_2}}
        &\cong&\Z[\l^{\pm 1},q_0^{\pm 1},q_1,\ldots,q_{z-\e}, A_0,B_1,B_2,\ldots,B_{z-\e}]/\langle q_0^{-2}-1\rangle\\
        &=&(\Z[q_0^{\pm 1},q_1,\ldots,q_{z-\e}, A_0,B_1,B_2,\ldots,B_{z-\e}]/\langle q_0^{-2}-1\rangle)[\l^{\pm 1}]
\end{eqnarray*}
is a Laurent polynomial ring in $\l$. Now, in this ring, $\beta_\sigma$ is one of $q_0\l-q_0^{-1}\l^{-1}$, $q_0\l-1$ or $q_0^{-1}\l^{-1}-1$. In every case, 
it has an invertible leading coefficient as a polynomial in $\l$, so it is not a zero divisor in  $\Omega_2/\langle h_0',q+\tau\rangle_{_{\Omega_2}}$. Now, because $q+\tau$ is not a zero divisor in $\Omega_2$ and $h_0'$ is not 
a zero divisor in the polynomial ring $\Omega_2/\langle q+\tau\rangle_{_{\Omega_2}}$, Lemma \ref{zerodiv} implies that $q+\tau$ is not a zero divisor in $\Omega_2/\la{h_0'}_{_{\Omega_2}}$. Then applying Lemma \ref{zerodiv} again shows that $q+\tau$ is not a zero divisor in $\Omega_2/\la{h_0',\beta_\sigma}_{_{\Omega_2}} = R_\sigma$, so we have proven (a) when $k$ is odd.

Now suppose $k$ is even, then $R_\sigma = \Omega_2/\la{\beta_\sigma,h_0,h_z,h_z'}_{_{\Omega_2}}$. Equation (\ref{hiprime}) then gives
\begin{align*}
    \d h_z'
        &\h (q_0^{-1}-1)h_z+\beta q_0^{-1}q_z-h_0q_z\\
        &\rel{eqn:hBl} \d\left[(q_0^{-1}-1)B_z-(A_0-1-q_0^{-1})q_z\right] \text{ in } \Omega.
\end{align*}
Since $\d$ is not a zero divisor in $\Omega$, this implies
\begin{equation}\label{hmprimep4}
    h_z'=(q_0^{-1}-1)B_z-(A_0-1-q_0^{-1})q_z.
\end{equation}
We first aim to show that $q+\tau$ is not a zero divisor in $\Omega_2/\langle h_z,h_z',h_0'\rangle_{_{\Omega_2}}$. Suppose
\[
    (q+\tau)x=ah_z+bh_z'+ch_0',
\]
for some $a,b,c\in\Omega_2$. For $d\in\Omega_2$, let $\bar d$ denote the image of $d$ in $\Omega_2/\langle q+\tau\rangle_{_{\Omega_2}}$. Then we have $\bar{h}_z = \l^{-1} (q_z + q_0^{-1} q_z)$, $\bar{h}_z' = (q_0^{-1}-1)B_z-(A_0-1-q_0^{-1})q_z$, and $h_0' = \l^{-1}(q_0^{-2} - 1)$. \\ Thus
\[
    \bar a\l^{-1}(q_0^{-1}+1)q_z+\bar b\left[(q_0^{-1}-1)B_z-(A_0-1-q_0^{-1})q_z\right]+\bar c\l^{-1}(q_0^{-2}-1)=0.
\]
In particular, $q_0^{-1}+1$ divides $\bar b\left[(q_0^{-1}-1)B_z-(A_0-1-q_0^{-1})q_z\right]$. Because $\Omega_2/\langle q+\tau\rangle_{_{\Omega_2}}$ is just 
the polynomial ring $\Z[\l^{\pm 1},q_0^{\pm 1},q_1,\ldots,q_z,  A_0,B_1,B_2,\ldots,B_z]$ (shown above), $q_0^{-1}+1$ divides $\bar b$. Thus $\bar b=\bar b_1\l^{-1}(q_0^{-1}+1)$, for some $b_1\in\Omega_2$, and so 
\[
    \bar a q_z+\bar b_1\left[(q_0^{-1}-1)B_z-(A_0-1-q_0^{-1})q_z\right]+\bar c(q_0^{-1}-1)=0.
\]
Rearranging then gives
\[
    (q_0^{-1}-1)\left[\bar b_1 B_z+\bar c\right]=q_z\left[\bar b_1(A_0-1-q_0^{-1})-\bar a\right].
\]
Now $q_0^{-1}-1$ and $q_z$ are coprime as elements of $\Omega_2/\la{q+\tau}_{\Omega_2}$, so there exists a $c_1\in\Omega_2$ such that \begin{align*}
    \bar b_1 B_z+\bar c \- \bar c_1q_z\quad \text{ and}\\
    \bar b_1(A_0-1-q_0^{-1})-\bar a \- \bar c_1(q_0^{-1}-1).
\end{align*}
We may now write
\begin{align*}
    a \- b_1(A_0-1-q_0^{-1})-c_1(q_0^{-1}-1)+(q+\tau)a_2,\\
    b \- \l^{-1}(q_0^{-1}+1)b_1+(q+\tau)b_2 \quad \text{ and}\\
    c \- c_1q_z-b_1B_z+(q+\tau)c_2,
\end{align*}
for some $a_2,\,b_2,\,c_2\in\Omega_2$. Thus
\begin{align*}
    (q+\tau)x
        \- ah_z+bh_z'+ch_0'\\
        \- b_1\left[(A_0-1-q_0^{-1})h_z+\l^{-1}(q_0^{-1}+1)h_z'-B_zh_0'\right]\\
        \f +c_1\left[q_zh_0'-(q_0^{-1}-1)h_z\right]
            +(q+\tau)(a_2h_z+b_2h_z'+c_2h_0').
\end{align*}
Using the definition of $h_0'$ and equations (\ref{eqn:hBl}) and (\ref{hmprimep4}), it is straightforward to verify that $\left[(A_0-1-q_0^{-1})h_z+\l^{-1}(q_0^{-1}+1)h_z'-B_zh_0'\right] = 0$. Also, by definition of $h_z'$ and $h_0'$, we know that $\d h_z' = (q_0^{-1}-1)h_z - q_zh_0'$. Hence the above reduces to
\begin{align*}
(q+\tau)x \- - \d c_1h_z'+(q+\tau)(a_2h_z+b_2h_z'+c_2h_0') \\
\- (q+\tau) \left[ -q^{-1} (q-\tau) c_1 h_z' + a_2h_z+b_2h_z'+c_2h_0' \right].
\end{align*}
as $\d = q^{-1}(q+\tau)(q-\tau)$.
Earlier we showed that $q+\tau$ is not a zero divisor in $\Omega_2$, so
\[
    x=a_2h_z+b_2h_z'+c_2h_0'- q^{-1}(q-\tau)c_1h_z'\in\langle h_z,h_z',h_0'\rangle_{_{\Omega_2}}.
\]
That is, $q+\tau$ is not a zero divisor in $\Omega_2/\langle h_z,h_z',h_0'\rangle_{_{\Omega_2}}$. Finally, by a similar reasoning as in the odd case, $\beta_\sigma$ is not a zero divisor in $\Omega_2/\langle h_z,h_z',h_0',q+\tau\rangle_{_{\Omega_2}}$. A final application of Lemma \ref{zerodiv} 
therefore shows that $q+\tau$ is not a zero divisor in $\Omega_2/\la{\beta_\sigma,h_0,h_z,h_z'}_{_{\Omega_2}}$, thereby completing the proof of (a).

%%%%%%%%%%%%%%%%%%%%%%%%%%%%%%%%%%%%%%%%%%

\textbf{(b)} For the moment, $\sigma \in \{0,+,-\}$. We now give a concrete realisation of the ring $R_\sigma[\d^{-1}]$. Let $I_\sigma[\d^{-1}]$ denote the ideal of $\Omega[\d^{-1}]$ generated by $I_\sigma$.
%; that is the ring of polynomials in $\d^{-1}$ with coefficients in $I_\sigma$. 
A standard argument shows that
\[
    R_\sigma[\d^{-1}]=(\Omega/I_\sigma)[\d^{-1}]\cong \Omega[\d^{-1}]/I_\sigma[\d^{-1}].
\]

By equation (\ref{hiprime}), $h_{k-l} = \d q_0 h_l' + q_0 h_l - \beta q_l - q_0 h_0q_l \in I_\sigma[\d^{-1}]$,  for $1\leq l\leq z-\e$. Thus the ideal in $\Omega[\d^{-1}]$ generated by $\beta_\sigma$ and all $h_0,h_1,\o,h_{k-1}$ must also be contained in $I_\sigma[\d^{-1}]$.
%\[
%    \langle\beta_\sigma,h_0,h_1,\ldots,h_{k-1}\rangle\subseteq I_\sigma[\d^{-1}].
%\]
Conversely, since $\d$ is invertible in $\Omega[\d^{-1}]$, equation (\ref{hiprime}) also shows that
\[
    h_l'\in\langle\beta_\sigma,h_0,h_1,\ldots,h_{k-1}\rangle_{_{\Omega[\d^{-1}]}}\text{,\quad for } 1\leq l\leq z-\e. 
\]
Thus
\[
    I_\sigma[\d^{-1}]=\langle\beta_\sigma,h_0,h_1,\ldots,h_{k-1}\rangle_{_{\Omega[\d^{-1}]}}.
\]
In $\Omega[\d^{-1}]$, the equations $h_l=0$ are equivalent to
\begin{eqnarray*}
    A_0&=& \d^{-1}\l^{-1}-\d^{-1}\l + 1,\\
    \text{and\quad}B_l&=&-\d^{-1}\l^{-1}(q_l+q_0^{-1}q_{k-l})\text{,\quad  for }1\leq l\leq k-1.
\end{eqnarray*}
Let $\Omega_3 := \Omega[\d^{-1}]/\langle h_0,h_1,\ldots,h_{k-1}\rangle_{_{\Omega[\d^{-1}]}}$. Then we have expressed $A_0$ and $B_1, \ldots, B_{k-1}$ in $\Omega_3$ as polynomials in $q^{\pm 1},\l^{\pm 1},q_0^{\pm 1},q_1,\ldots,q_{k-1}$ and $\d^{-1}$.
Therefore, by (\ref{eqn:Lambdaalt}),
\[
 \Omega_3 \cong \Z[q^{\pm 1},\l^{\pm 1},q_0^{\pm 1},q_1,\ldots,q_{k-1}][\d^{-1}].
\]
Moreover,
\[
R_\sigma[\d^{-1}] \cong \Omega[\d^{-1}]/I_\sigma[\d^{-1}] = \Omega[\d^{-1}]/\la{\beta_\sigma,h_0,h_1,\ldots,h_{k-1}}_{_{\Omega[\d^{-1}]}} \cong
\Omega_3/\langle\beta_\sigma\rangle_{_{\Omega_3}}.
\]
Now suppose $\sigma=\pm$. Then $\beta_\sigma$ can be ``solved'' for $q_0^{\pm 1}$, so
 \[
    R_\sigma[\d^{-1}]\cong\Z[q^{\pm 1},\l^{\pm 1},q_1,\ldots,q_{k-1}][\d^{-1}],
\]
completing the proof of (b).

\textbf{(c)} Observe that the ring  $\Z[q^{\pm 1},\l^{\pm 1},q_1,\ldots,q_{k-1}][\d^{-1}]$ above is obtained from an integral domain via localisation, hence $R_\sigma[\d^{-1}]$ is also an integral domain. 
Now, we have already proven in part (a) that $\d$ is not a zero divisor in $R_\sigma$. Therefore the map $R_\sigma\rightarrow R_\sigma[\d^{-1}]$ is injective, and statement (c) now follows immediately.

\textbf{(d)} Because $\beta_0 = \beta_+ \beta_-$, it is clear that $\beta_0 \in I_+ \cap I_-$, thus all generators of the ideal $I_0$ are in $I_+ \cap I_-$. Hence $I_0 \subseteq I_+ \cap I_-$. 

Finally, note that $\Omega_3 \cong \Z[q^{\pm 1},\l^{\pm 1},q_0^{\pm 1},q_1,\ldots,q_{k-1}][\d^{-1}]$ is obtained 
from a UFD (unique factorisation domain) by localisation, and is therefore also a UFD. Suppose that $x\in I_+\cap I_-\subseteq\Omega$. 
Then $x$ vanishes in $R_\pm$ and hence in $R_\pm[\d^{-1}] \cong \Omega_3/\la{\beta_\pm}_{_{\Omega_3}}$. So the image $\bar x$ of $x$ 
in $\Omega_3$ satisfies
\[
    \bar x\in\langle\beta_+\rangle_{_{\Omega_3}} \cap \langle\beta_-\rangle_{_{\Omega_3}} = \langle\beta_+\beta_-\rangle_{_{\Omega_3}} = \langle\beta\rangle_{_{\Omega_3}}, 
\]
since $\beta_+$ and $\beta_-$ are coprime in $\Omega_3$. Thus $x$ maps to $0$ in $R_0[\d^{-1}] \cong \Omega_3/\la{\beta}_{_{\Omega_3}}$. Since $R_0$ embeds into $R_0[\d^{-1}]$, the image 
of $x$ in $R_0$ must also be $0$. That is, $x\in I_0$, hence $I_+\cap I_-\subseteq I_0$, completing the proof of (d).
\end{proof}

\noi Note that, in $R_0$, we have
\[
    h_l=q_0h_{k-l}-\beta q_{k-l}+h_0q_0q_{k-l}+\d q_0h_{k-l}'=0,
\]
for all $z-\e<l\leq k-1$, by (\ref{hiprime}). We have just shown that $\d$ is not a zero divisor in $R_0$, so we may conclude from Lemma 3.4 of \cite{WY06} the 
existence of a $\bt(R_0)$-module $V(R_0)$ satisfying (\ref{wantmodule}). Since $R_0$ is the ``generic" ring in which the above relations hold, we 
may now specialize this result to the following class of ground rings:

\begin{defn} \label{defn:adm}
Let $R$ be as in the definition of $\bnk$ (see Definition \ref{defn:bnk}). The family of parameters $\left( A_0, \ldots, A_{k-1}, q_0, \ldots, q_{k-1}, q, \l\right)$ is called \textbf{\emph{admissible}} if 
\[
\beta=h_0=h_1=\ldots=h_{z-\e}=h_1'=h_2'=\ldots=h_{z-\e}'=0.
\]
\end{defn}

\textbf{Remark:} The admissibility conditions given in Definition \ref{defn:adm} are the most general conditions offered in the literature and are \emph{necessary and sufficient} for the freeness and isomorphism results in Corollary \ref{cor:basis}. Also, in \cite{WY06}, we give a definition of admissibility in the special case when $\delta$ is not a zero divisor. It is 
straightforward to see that the two definitions are equivalent in this special case. Moreover, under the assumption that $\delta$ is not a zero divisor, we proved that admissibility is equivalent to the existence of a module satisfying (\ref{wantmodule}), and to $\bt$ being a free $R$-module of rank $3k^2$ 
(see Corollary 4.5 of \cite{WY06}). With the above definition, this result still holds in this more general context. Indeed freeness implies the existence of $V$, as shown in Corollary 4.5 of~\cite{WY06},
and admissibility will imply freeness by Corollary \ref{cor:basis}. The fact that the
existence of $V$ implies admissibility follows from a laborious calculation which we omit, as we do not require the result here.
\vskip 0.2cm
For all $\sigma \in \{0,+,-\}$, $R_\sigma$ is a ring with admissible parameters, by definition of $I_\sigma$ in Lemma~\ref{lemma:Rsigma}. 
Moreover $R_0$ is the generic admissible ring in the following sense.
\begin{prop} \label{prop:universalmap}
  Let $R$ be as in Definition \ref{defn:bnk} with admissible parameters $A_0$, \ldots, $A_{k-1}$, $q_0, \ldots$, $q_{k-1}$, $q$ and $\l$. Then there exists a unique map $R_0\rightarrow R$ which respects the parameters.
\end{prop}
\begin{proof}
  There is a unique ring map $\rho: \Omega\rightarrow R$ which respects the parameters. Furthermore, the admissibility of the parameters in $R$ is equivalent to \[
\rho(\langle\beta,h_0,h_1,\ldots,h_{z-\e},h_1',\ldots,h_{z-\e}'\rangle _\Omega)=0.
\]
That is, the map $\rho$ kills $I_0 \subseteq \Omega$ and hence factors through $R_0$.
\end{proof}

In addition to the Remarks made after Definition \ref{defn:bnk}, on the differences in the assumptions on the parameters of the ground ring, it is also important to relate the various notions of admissibility used in the literature. 
Before proceeding with the comparison, we first draw particular attention to ``weak admissibility", used by Goodman and Hauschild Mosley in \cite{GH108,GH208}. (For the specific history of these conditions, please see \cite{Y07}). Parameters that are admissible, in the sense of this paper, are also weakly admissible. Weak admissibility may be viewed as a minimal condition for which the algebra does not collapse. Certainly, in the absence of weak admissibility, the $e_i$ become torsion elements and, over a field, $\bnk$ would reduce to the Ariki-Koike algebras. Also, it turns out that weak admissibility of the parameters in $R$ is sufficient enough to ensure the algebra at $n=0$ is non-trivial and to produce a (nondegenerate) trace function $\bnk(R) \rightarrow R$ in Section \ref{sect:traceconstruction}.

Suppose the parameters in $R$ are admissible. We may tensor our $\bt(R_0)$-module 
$V(R_0)$, satisfying (\ref{wantmodule}), with $R$ and use the natural homomorphism $\bt(R)\rightarrow R\otimes_{R_0}\bt(R_0)$ to obtain a $\bt(R)$-module
\[
	V=R\otimes_{R_0}V(R_0)
\]
which satisfies (\ref{wantmodule}).

%%%%%%%%%%%%%%%%%%%%%%%%%%%%% adm => weak adm

Let us define elements $v_i\in V$ for $i<0$ and $i\geq k$ by $v_i=Y^iv_0$. Also there are unique elements $A_i\in R$ for 
$i<0$ and $i\geq k$ such that
\begin{equation}\label{extendA}
	\sum_{i=0}^kq_iA_{i+j}=0
\end{equation}
for all $j\in\Z$. It follows that $ev_i=A_iv_0$ for $i\in\Z$. Also note that $Y_2'=XYX$ acts as $Y^{-1}$ on $V$. 
\pagebreak
Therefore applying (\ref{eqn:magic}) (with $i=1$) to $v_0$, we have
\[
	Xv_p=\l A_{-p}v_0-\d\pum{s}(v_{p-2s}-Y^{-s}A_{p-s}v_0).
\]
Now applying $e$, and taking the coefficient of $v_0$,
\begin{equation}\label{weakA}
	\l A_p=\l A_{-p}-\d\pum{s}(A_{p-2s}-A_{-s}A_{p-s}),
\end{equation}
for all $p\geq1$. These relations amongst the parameters give rise to another class of ground rings, defined by Goodman and Hauschild Mosley \cite{GH108,GH208}:
\begin{defn} \label{defn:GHadm}
Let $R$ be as in the definition of $\bnk$ (see Definition \ref{defn:bnk}). The family of parameters 
$\left( A_0, \ldots, A_{k-1}, q_0, \ldots, q_{k-1}, q, \l\right)$ is called \textbf{\emph{weakly admissible}} if $h_0=0$ and (\ref{weakA}) holds 
for all $p\geq1$, where the $A_i$ are defined for $i<0$ and $i\geq k$ by (\ref{extendA}).
\end{defn}

\textbf{Remarks:}\quad (1) \,\, In fact, Goodman and Hauschild Mosley begin with a ring containing elements $\{A_i\mid i\geq0\}$ then impose infinitely many relations which are polynomials in the $A_i$'s by defining $A_i$ for $i<0$ using (\ref{weakA}), and say the parameters are weakly admissible if (\ref{extendA}) holds. It is easy to see that such choices of parameters correspond exactly to the above definition.

(2) \,\, The above discussion shows that an admissible family of parameters is also weakly admissible.

(3) \,\, In \cite{GH108,GH208}, Goodman and Hauschild Mosley defines admissibility to be the collective conditions used in \cite{WY06}; they assume the parameters of $R$ satisfy $\beta = h_0 = \ldots = h_{k-1} = 0$ and that $\delta$ is a not a zero divisor in $R$. As noted earlier, this coincides with Definition~\ref{defn:adm} of admissibility together with the extra assumption that $\delta$ is not a zero divisor. 

(4) \,\, It is worth noting here that, strongly influenced by the methods and proofs for the cyclotomic Nazarov-Wenzl algebra in Ariki et al.~\!\cite{AMR06}, Goodman and Hauschild Mosley \cite{GH208} and Rui-Si-Xu \cite{RX08,RS08} find that the existence of an (irreducible) $k$-dimensional $\bt$-module leads them to also formulate a stronger third type of admissibility condition called $\mathbf{u}$-admissibility (differing from each other). In \cite{GH208}, it is used as a step towards constructing rings with admissible parameters; a ring $R$ satisfies $\mathbf{u}$-admissibility if it is an integral domain with admissible parameters in the sense of \cite{WY06}, and the roots of the $k^\mathrm{th}$ order relation and their inverses are all distinct.  Hence, if the parameters of $R$ are $\mathbf{u}$-admissible in the sense of \cite{GH208}, then they are also admissible in the sense of Definition \ref{defn:adm}.

On the other hand, $\mathbf{u}$-admissibility in \cite{RX08,RS08} is closer to (and implies) the weak admissibility conditions defined in this paper. Freeness and cellularity of $\bnk$ are established over these conditions in \cite{RX08,RS08}. In particular, if the parameters of $R$ are $\mathbf{u}$-admissible in the sense of \cite{RX08,RS08}, then $\bt(R)$ is free of rank $3k^2$ hence they are also admissible in the sense of Definition~\ref{defn:adm}.
For more details, we refer the reader to \cite{AMR06,GH208,RX08,RS08}.

\section{The Cyclotomic Kauffman Tangle Algebras} \label{sect:ckt}

Recall in the introduction we mentioned that the BMW algebra $\bmw_n$ is isomorphic to the Kauffman tangle algebra $\mathbb{KT}_n$ and is of rank $(2n-1)!! = (2n-1) \cdot (2n-3) \cdot \o \cdot 1$, the same as that of the Brauer algebras.
The Brauer algebras were introduced by Brauer~\cite{B37} as a device for studying the representation theory of the symplectic and orthogonal groups, and are typically defined to have a basis consisting of Brauer diagrams, which are basically permutation diagrams except horizontal arcs between vertices in the same row are now also permitted. The BMW algebras are a deformation of the Brauer algebras comparable to the way the Iwahori-Hecke algebras of type $A_{n-1}$ are a deformation of the group algebras of the symmetric group $\mathfrak{S}_n$. Alternatively, the Brauer algebra is the ``classical limit" of the BMW algebra or Kauffman tangle algebra in the sense one just ``forgets" the notion of over and under crossings in tangles diagrams and so tangle diagrams consisting only of vertical strands degenerate into permutations. In fact, by starting with the set of Brauer $n$-diagrams (see beginning of Section \ref{section:cyclobrauer}), together with a fixed ordering of the vertices and a rule for which strands cross over which, one may easily write down a diagrammatic basis of the BMW algebra $\mathscr{C}_n$; for more details on this construction, we refer the reader to Morton and Wasserman \cite{MW89} and Halverson and Ram \cite{HR95}.

The affine and cyclotomic Brauer algebras were introduced by H\"aring-Oldenburg in \cite{HO01} as classical limits of their BMW analogues in the above sense. The cyclotomic case is studied by Rui and Yu in \cite{RY04} and Rui and Xu in \cite{RX07}. (There is also the notion of a $G$-Brauer algebra for an arbitrary abelian group $G$, introduced by Parvathi and Savithri in \cite{PS02}). The cyclotomic Brauer algebra is free of rank $k^n (2n-1)!!$, by definition. Therefore one would expect the cyclotomic BMW algebras to be of this rank too.
The main aim of this section is to establish the linear independence of our spanning set, obtained in Section \ref{sect:spanning}, initially over $R_0$, the generic ground ring with admissible parameters constructed in Lemma \ref{lemma:Rsigma}. To achieve this goal, we require the existence of a trace on $\bnk(R_0)$ which, using its specialisation into the cyclotomic Brauer algebra, is shown to be nondegenerate.

We begin this section with a brief introduction on tangles and affine tangles and define the cyclotomic Kauffman tangle algebras through the affine Kauffman tangle algebras. From here, we introduce the cyclotomic Brauer algebras and associate with it a trace map. Using the nondegeneracy of this trace, given by a result of Parvathi and Savithri \cite{PS02}, we are then able to establish the nondegeneracy of a trace on the cyclotomic BMW algebras over a specific quotient ring of $R_\sigma$, for $\sigma \in \{0,+,-\}$. From this, we then deduce that the same result holds for these three rings. In particular, this implies that we have a basis of $\bnk(R_0)$, thereby proving $\bnk(R)$ is $R$-free for \emph{all} rings $R$ with admissible parameters. Moreover, as a consequence of these results, we also prove the cyclotomic BMW algebras are isomorphic to the cyclotomic Kauffman tangle algebras. %These results verify the claims made in Goodman and Hauschild Mosley \cite{GH107} for rings with \emph{weakly admissible} parameters. (<-GUESS SHOULD NOT REFER TO THE OLD PREPRINTS ANYMORE???)

%\begin{comment}
%\begin{defn}
%  Fix points $a_1,\o,a_n$ in the interval $I=[0,1]$. An \emph{$n$-tangle} is a piece of a KNOT/LINK? diagram consisting of a union of curves $F$ in $I \times I$ such that \\
%(1) \vspace{-7mm} \begin{align*} 
%  F \cap \partial(I \times I) \- F \cap (I \times \{0,1\}) \\
%\- \{(a_1,1), \o, (a_n,1)\} \cup \{(a_1,0), \o, (a_n,0)\} 
%\end{align*}
%and \\
%(2) \quad any curve in $F$ intersects the boundary of $I\times I$ transversally.
%\end{defn} 
%\end{comment}

\begin{defn}%[Alternative defn of tangle taken from Xi \cite{X00} - NOTE THIS INCLUDES CYCLES]
  An \emph{\textbf{$n$-tangle}} is a piece of a link diagram,  consisting of a union of arcs and a finite number of closed cycles, in a rectangle in the plane such that the end points of the arcs consist of $n$ points located at the top and $n$ points at the bottom in some fixed position.
\end{defn}

An $n$-tangle may be diagrammatically presented as two rows of $n$ vertices and $n$ strands connecting the vertices so that every vertex is incident to precisely one strand, and over and under-crossings and self-intersections are indicated. In addition, this diagram may contain finitely many closed cycles.
 
\begin{figure}[h!]
  \begin{center}
  \begin{tabular}{rc}
\raisebox{.5cm}{I}\qquad\qquad\qquad &   \includegraphics{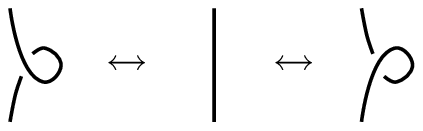} \\[0.4cm]
\raisebox{.55cm}{II}\qquad\qquad\qquad &   \, \includegraphics{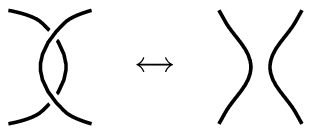} \\[0.4cm]
\raisebox{.65cm}{III}\qquad\qquad\qquad &    \, \includegraphics{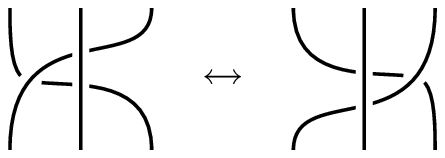}
\end{tabular}
    \caption{Reidemeister moves of types I, II and III.} \label{fig:reidemeister}
  \end{center}
\vspace{-8mm}
\end{figure}

\begin{defn}
  Two tangles are said to be \emph{\textbf{ambient isotopic}} if they are related by a sequence of \emph{\textbf{Reidemeister moves}} of types I, II and III (see Figure \ref{fig:reidemeister}), together with an isotopy of the rectangle which fixes the boundary. They are \emph{\textbf{regularly isotopic}} if the Reidemeister move of type I is omitted from the previous definition.
\end{defn}
One obtains a monoid structure on the regular isotopy equivalence classes of $n$-tangles where composition is defined by concatenation of diagrams.
As mentioned in the introduction, the algebra of these tangles together with Kauffman skein relations is a diagrammatic formulation of the original BMW algebra. The tangles which appear in the topological interpretation of the affine and cyclotomic BMW algebras feature in type $B$ braid/knot theory. Affine braids or braids of type $B$ are essentially just ordinary braids (of type $A$) on $n+1$ strands in which the first strand is pointwise fixed; this single fixed line is usually presented as a ``flagpole'', a thickened vertical segment on the left, and the other strands may loop around this flagpole. They are also commonly depicted as braids in a (slightly thickened) cylinder or a (slightly thickened) annulus; see Lambropoulou~\cite{L93}, tom~Dieck \cite{T94} and Allcock \cite{A02}.
 
\begin{defn}
  An \emph{\textbf{affine $n$-tangle}} is an $n+1$-tangle with a monotonic path joining the first top and bottom vertex, which is diagrammatically presented by the flagpole mentioned above.
\end{defn}
\begin{figure}[h!]
  \begin{center}
    \includegraphics[width=35mm,height=30mm] {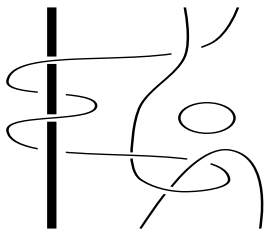}
\caption{Affine $2$-tangle diagram.}
\label{fig:afftang}
  \end{center}
\end{figure}
\vskip -5mm 
Two affine $n$-tangles are ambient (regularly, respectively) isotopic if they are ambient (regularly, respectively) isotopic as $n+1$-tangles. Note that, as the flagpole is required to be a fixed monotonic path, the Reidemeister move of type I is never applied to the flagpole in an ambient isotopy. As with ordinary tangles, the equivalence classes of affine $n$-tangles under regular isotopy carry a monoid structure under concatenation of tangle diagrams. Let $\wh{\mathbb{T}}_n$ denote the monoid of the regular isotopy equivalence classes of affine $n$-tangles.

For $j \geq 0$, let us denote by $\Theta_j$ (the regular isotopy equivalence class of) the non-self-intersecting closed curve which winds around the flagpole in the `positive sense' $j$ times. These are special affine $0$-tangles which will feature in definitions later. Observe that $\Theta_0$ is represented by a closed curve that does not interact with the flagpole. The following figure illustrates $\Theta_3$:
  \begin{center}
    \includegraphics[scale=0.85]{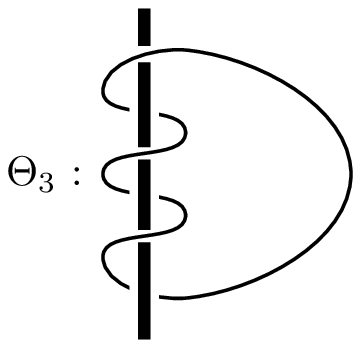}
  \end{center}

Using the monoid algebra of affine $n$-tangles, we may now define the affine and cyclotomic Kauffman tangle algebras. We remark here that the definitions differ slightly to those given in Goodman and Hauschild Mosley \cite{GH108,GH208}; the difference is that their initial ground ring involves an infinite family of $A_j$'s, for every $j \geq 0$. Instead, here we define the algebras over a ring $R$ under the same assumptions as in the definition of $\bnk$ and take $A_j$, where $j \geq k$, to be elements of $R$ defined by equation (\ref{extendA}).

The figures in relations given in the following two definitions indicate affine tangle diagrams which differ locally only in the region shown and are identical otherwise. 

\begin{defn} \label{defn:affineKT}
Let $R$ be as in Definition \ref{defn:bnk}; that is, a commutative unital ring containing units $q_0, \ldots, q_{k-1},q,\l$ and further elements $A_0, A_1, \ldots, A_{k-1}$ such that $\l - \l^{-1} = \d(1-A_0)$ holds, where $\d = q - q^{-1}$. Moreover, let $A_j \in R$, for all $j \geq k$, be defined by equation (\ref{extendA}). \\
The \emph{\textbf{affine Kauffman tangle algebra}}  $\akt{n}(R)$ is the monoid $R$-algebra $R\wh{\mathbb{T}}_n$ modulo the following relations:
\begin{enumerate}
\item[(1)] (Kauffman skein relation)
  \begin{center}
    \includegraphics{pspictures/kauff_scaled.eps}
  \end{center}
\item[(2)] (Untwisting relation) 
  \begin{center}
    \includegraphics{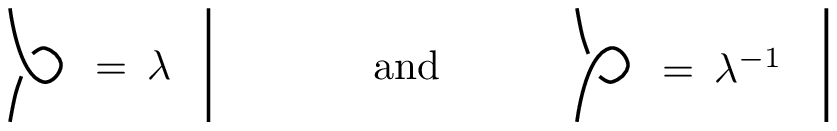}   
  \end{center}
\item[(3)] (Free loop relations) \\
For $j \geq 0$, 
\[ T \amalg \Theta_j = A_j T, \]
where $T \amalg \Theta_j$ is the diagram consisting of the affine $n$-tangle $T$ and a copy of the loop $\Theta_j$ defined above, such that there are no crossings between $T$ and $\Theta_j$.
\end{enumerate}
\end{defn}

\begin{figure}[h!]
  \begin{center}
    \includegraphics[scale=0.8]{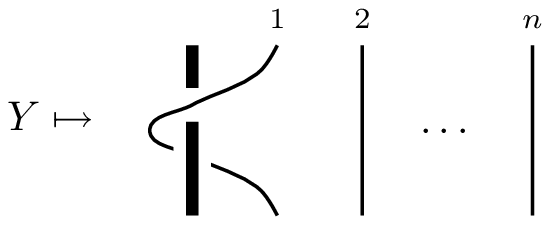} \vskip 0.65cm
    \includegraphics[scale=0.8]{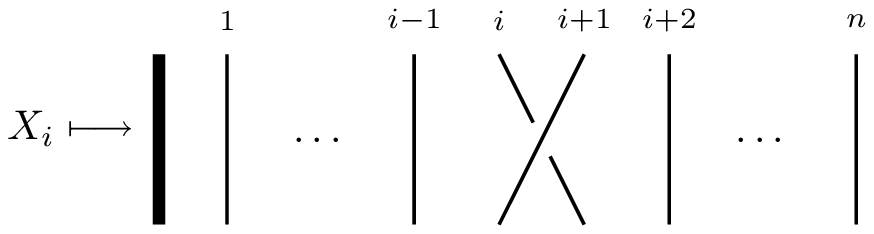}
\vskip 0.65cm
    \includegraphics[scale=0.8]{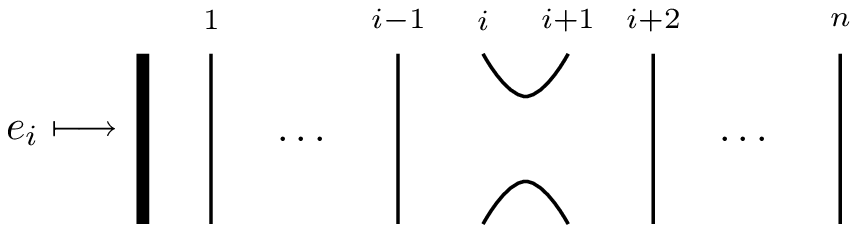}
    \caption{The isomomorphism between the affine BMW and affine Kauffman tangle algebras.}
\label{fig:dghom}
  \end{center}
\end{figure}

\begin{figure}[h!]
\vskip -2mm
  \begin{center}
    \includegraphics[scale=0.77]{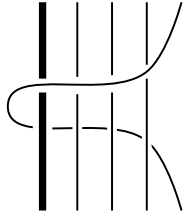}
    \caption{The affine $4$-tangle associated with the element $Y_4'$.}
    \label{fig:y4dash}
  \end{center}
\end{figure}

\textbf{Remark:} %(\textbf{actually how can this still be true for any ring $R$ given that the $\Theta_j$'s are no longer ``algebraically independent''; see p.7 of affine BMW algebra paper}) 
An important case to consider is the affine $0$-tangle algebra $\akt{0}(R)$. If relation~(3) is removed from the above definition for $j\geq1$, a result of Turaev \cite{T88} shows that the affine $0$-tangle algebra is freely generated by the $\Theta_j$, where $j \geq 1$, and embeds in the center of the affine $n$-tangle algebra. This motivates relation~(3) in the above definition. This then shows $\akt{0}(R) \cong R$.  We shall see later that weakly admissible conditions are sufficient to prove the analogous result in the cyclotomic case.
\vskip 0.2cm
Recall the affine BMW algebra $\wh{\mathscr{B}}_n(R)$ over $R$ is simply the cyclotomic BMW algebra $\bnk(R)$ with the $k^\mathrm{th}$ order polynomial relation on the generator $Y$ omitted. Goodman and Hauschild~\cite{GH06} prove the maps given in Figure \ref{fig:dghom} determine an $R$-algebra isomorphism $\wh{\dg}$ between the affine BMW and affine Kauffman tangle algebras.
We write $\mathcal{Y}$, $\mathcal{X}_i$ and $\mathcal{E}_i$ for the images of the generators $Y$, $X_i$ and $e_i$ under $\wh{\dg}$, respectively.
In particular, the image of $Y_i'$ is exemplified in Figure \ref{fig:y4dash}.
\vskip -4mm
\begin{defn} \label{defn:cycloKT}
  Let $R$ be as in Definition \ref{defn:affineKT}.
The \emph{\textbf{cyclotomic Kauffman tangle algebra}} $\ckt{n}(R)$ is the affine Kauffman tangle algebra $\akt{n}(R)$ modulo the cyclotomic skein relation:
  \begin{center}
\includegraphics[scale=1.2]{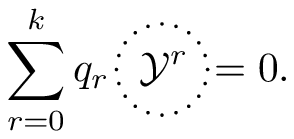}
  \end{center}
The interior of the disc shown in the cyclotomic skein relation represents part of an affine tangle diagram isotopic to the affine $1$-tangle $\mathcal{Y}^r$, where $\mathcal{Y}$ here is the tangle diagram illustrated in Figure \ref{fig:dghom} when $n=1$.
The sum in this relation is over affine tangle diagrams which differ only in the interior of the disc shown and are otherwise identical. 
\end{defn}

By definition, there is a natural projection $\pi_t: \akt{n}(R) \twoheadrightarrow \ckt{n}(R)$ and a natural projection $\pi_b: \wh{\mathscr{B}}_n(R) \twoheadrightarrow \bnk(R)$. Moreover, $\hat{\dg}: \wh{\mathscr{B}}_n(R) \rightarrow \akt{n}(R)$ induces an $R$-algebra homomoprhism $\dg: \bnk(R) \rightarrow \ckt{n}(R)$ such that the following diagram of $R$-algebra homomorphisms commutes. 
\[
\begin{CD}
\wh{\mathscr{B}}_n(R) @>\hat{\dg}>> \akt{n}(R) \\
@V{\pi_b}VV    @VV{\pi_t}V \\
\bnk(R) @>\dg>> \ckt{n}(R)
\end{CD}
\]
\\
Moreover, because $\wh{\dg}$ is an isomorphism, this implies $\dg: \bnk(R) \rightarrow \ckt{n}(R)$ is surjective. \label{defn:dg}
Furthermore, the homomorphism $\dg$ commutes with specialisation of rings. More precisely, given a parameter preserving ring homomorphism $R_1 \rightarrow R_2$, we can consider $\bnk(R_2)$ as an $R_1$-algebra and construct the $R_1$-algebra homomorphism $\eta_b: \bnk(R_1) \rightarrow \bnk(R_2)$, which sends generator to generator. The ring homomorphism also extends to a $R_1$-algebra homomorphism $\eta_t: \ckt{n}(R_1) \rightarrow \ckt{n}(R_2)$. Then it is easy to verify that $\dg \circ \eta_b = \eta_t \circ \dg$ holds on the generators of $\bnk(R_1)$, hence we have the following commutative diagram of $R_1$-algebra homomorphisms:
\begin{equation} \label{cd:dgspec}
\begin{CD}
\bnk(R_1) @>\eta_b>> \bnk(R_2) \\
@V{\dg}VV    @VV{\dg}V \\
\ckt{n}(R_1) @>\eta_t>> \ckt{n}(R_2)
\end{CD}
\end{equation}

Also, the tangle analogue of the (\ref{eqn:star}) anti-involution described at the beginning of Section \ref{sec:prelims} is then just  the anti-automorphism of $\ckt{n}(R)$ which flips diagrams top to bottom. In particular, it fixes $\mathcal{Y}$, $\mathcal{X}_i$ and $\mathcal{E}_i$. 
Furthermore, the map of affine tangles that reverses all crossings, including crossings of strands with the flagpole, determines an isomorphism from $\ckt{n}(q,\l,A_i,q_i)$ to $\ckt{n}(q^{-1},\l^{-1},A_{-i},-q_{k-i}q_0^{-1})$.

\textbf{Remark:} For integers $n,m\geq0$, one may similarly construct the affine Kauffman tangles $\akt{n,m}(R)$ and cyclotomic Kauffman 
tangles $\ckt{n,m}(R)$ with $n$ strands at the top of the diagram and $m$ at the bottom (not including the flagpole). (Note that $\akt{n}(R) = \akt{n,n}(R)$). These $R$-modules do not 
have an algebra structure, but we do have bilinear maps
\begin{eqnarray*}
	\akt{n,m}(R)\times\akt{m,l}(R)&\rightarrow&\akt{n,l}(R),\\
	\ckt{n,m}(R)\times\ckt{m,l}(R)&\rightarrow&\ckt{n,l}(R).
\end{eqnarray*}

\section{Construction of a Trace on $\bnk$} \label{sect:traceconstruction}

We now work our way towards a Markov-type trace on the cyclotomic BMW algebras via maps on the affine BMW algebras described in Goodman and Hauschild \cite{GH06}.
Observe that, there is a natural inclusion map $\iota$ from the set of affine $(n-1)$-tangles to the set of affine $n\text{-tangles}$ defined by simply adding an additional strand on the right without imposing any further crossings, as illustrated below.  
  \begin{center}
    \includegraphics{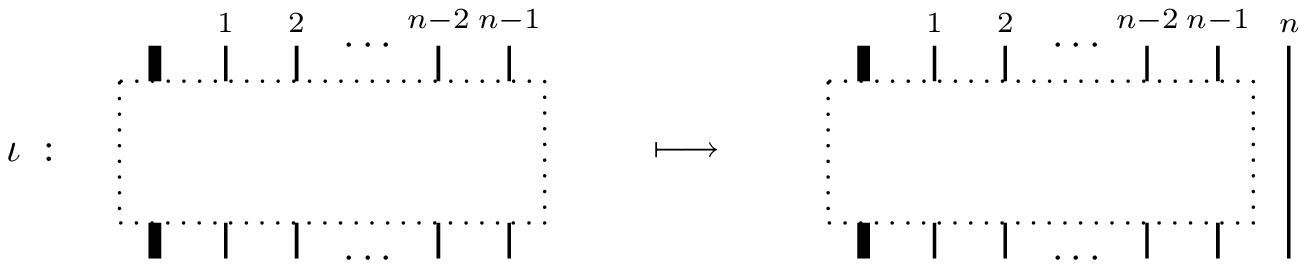}   
  \end{center}
  
\noi Furthermore, the map $\iota$ respects regular isotopy, composition of affine tangle diagrams and the relations of $\akt{n}$, so induces an $R$-algebra homomorphism $\iota: \akt{n-1}(R) \rightarrow \akt{n}(R)$. Moreover, it respects the cyclotomic skein relation, hence induces an $R$-algebra homomorphism $\iota: \ckt{n-1}(R) \rightarrow \ckt{n}(R)$. 

There is also a ``closure'' map $\cl{n}$ for the affine Kauffman tangle algebras, from the set of affine $n$-tangles to the set of affine $(n-1)$-tangles, given by closure of the rightmost strand, as illustrated below. We extend $\cl{n}$ to an $R$-linear map $\varepsilon_n: \akt{n}(R) \rightarrow \akt{n-1}(R)$. Taking these closure 
maps recursively produces a \emph{trace} map on $\akt{n}$. 
We define $\varepsilon_{_{\mathbb{T}}}: \akt{n}(R) \rightarrow \akt{0}(R) \cong R$ by
\[
\varepsilon_{_{\mathbb{T}}} := \varepsilon_1 \circ \cdots \circ \varepsilon_n.
\]
Moreover, these closure and trace maps respect the cyclotomic skein relation, therefore induce analogous maps $\varepsilon_n: \ckt{n}(R) \rightarrow \ckt{n-1}(R)$ and $\varepsilon_{_{\mathbb{T}}}: \ckt{n}(R) \rightarrow \ckt{0}(R)$.
\vskip 0.5cm
  \begin{center}
    \includegraphics{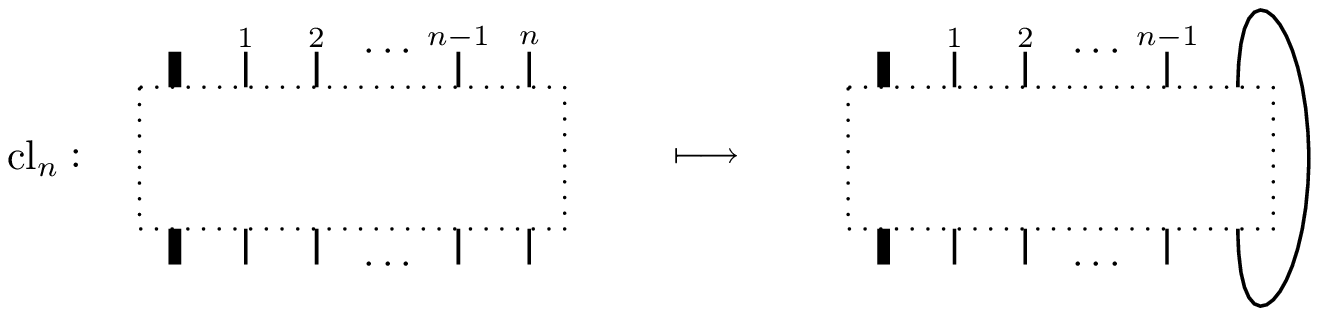}   
  \end{center}
\begin{center}
  \raisebox{2cm}{{\Large $\varepsilon_{\mathbb{T}}$\,:}} \quad \includegraphics[scale=0.8]{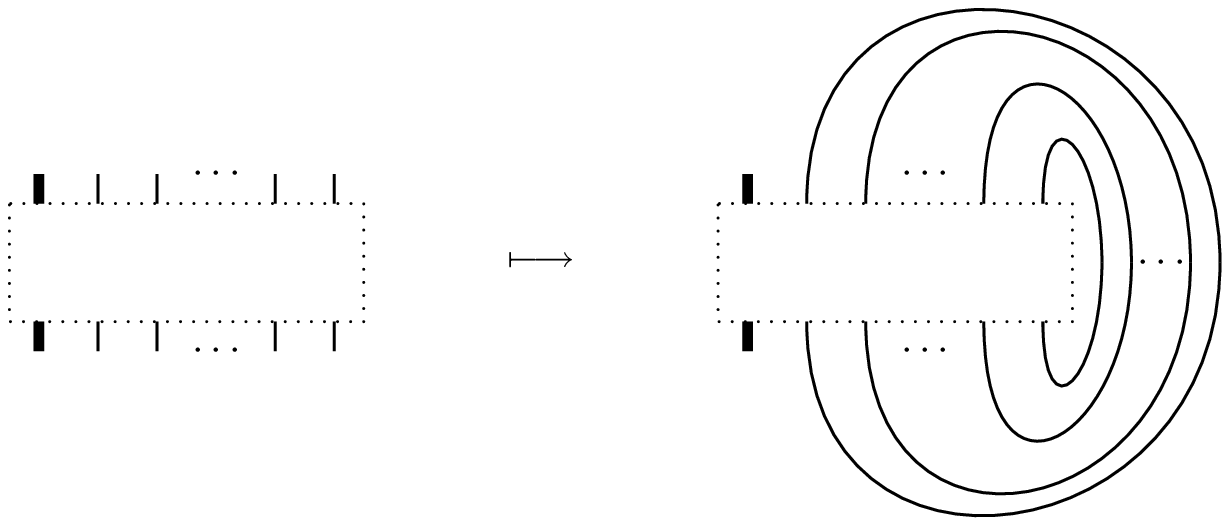} 
\end{center}

Our next step is to prove $\ckt{0}(R) \cong R$, where $R$ is a ring with weakly admissible parameters, so that $\varepsilon_{_{\mathbb{T}}}$ really yields a (unnormalised) trace map of $\ckt{n}$ over $R$, which will lead to one for $\bnk$. This requires Turaev's result on the affine level and further topological arguments. Weak admissibility may be thought of as a minimal condition for which this is true and in fact, over a field, the algebra $\ckt{0}$ would collapse in the absence of weak admissibility. 

Let $\cap\in\akt{0,2}(R)$ and $\cup\in\akt{2,0}(R)$ denote the tangles with a single strand without 
self-intersections, and let $\mathcal{Y}$ denote the tangle diagram in $\akt{2}(R)$ given by $\wh{\dg}(Y)$, as seen above. For $j\geq0$, we have
\[
	\cap\mathcal{Y}^j\cup=\Theta_j=A_j,
\]
by relation (3) in the definition of $\akt{0}(R)$. It can be shown (see Lemma 2.8 of \cite{GH06}) that the above equation also 
holds for $j<0$, provided we define $A_j$ for $j<0$ using (\ref{weakA}).
\begin{lemma} \label{lemma:akt20}
The $R$-module $\akt{2,0}(R)$ is spanned by $\{\mathcal{Y}^i\cup\mid i\in\Z\}$.
\end{lemma}
\begin{proof}
Clearly every tangle diagram in $\akt{2,0}(R)$ can be expressed as a affine $2$-tangle times $\cup$, so right multiplication by 
$\cup$ gives a surjection $\akt{2}(R)\twoheadrightarrow\akt{2,0}(R)$. Now Corollary 5.9 of \cite{GH06} shows that $\akt{2}(R)$ is spanned by
\[
	\{\mathcal{Y}^i\mathcal{E}_1\mathcal{Y}^j,\,\mathcal{Y}^i\mathcal{X}_1\mathcal{Y}^j,
		\mathcal{Y}^i\mathcal{X}_1\mathcal{Y}^j\mathcal{X}_1^{-1}\mid i,j\in\Z\},
\]
so $\akt{2,0}(R)$ is spanned by $\{\mathcal{Y}^i\mathcal{E}_1\mathcal{Y}^j\cup,\,\mathcal{Y}^i\mathcal{X}_1\mathcal{Y}^j\cup,
\mathcal{Y}^i\mathcal{X}_1\mathcal{Y}^j\mathcal{X}_1^{-1}\cup\mid i,j\in\Z\}$.
However,
\[
	\mathcal{Y}^i\mathcal{E}_1\mathcal{Y}^j\cup=A_j\mathcal{Y}^i\cup
	\hspace{10mm}\text{and}\hspace{10mm}
	\mathcal{Y}^i\mathcal{X}_1\mathcal{Y}^j\mathcal{X}_1^{-1}\cup
		=\l^{-1}\mathcal{Y}^i\mathcal{X}_1\mathcal{Y}^j\cup.
\]
Moreover a calculation exactly analogous to the proofs of (\ref{eqn:magic}) and (\ref{eqn:m2}) shows that
\[
	\mathcal{X}_1\mathcal{Y}^j\cup\in\la{\mathcal{Y}^k\cup\mid k\in\Z}.
\]
The result now follows.
\end{proof}
\begin{prop}
Suppose the parameters in $R$ are weakly admissible (see Definition \ref{defn:GHadm}). Then the $R$-algebra $\ckt{0}(R)$ is isomorphic to $R$.
\end{prop}
\begin{proof}
Recall that $\akt{0}(R)\cong R$, by Turaev's result. It therefore suffices to show the canonical surjection $\akt{0}(R)\twoheadrightarrow\ckt{0}(R)$ is injective. By definition, the kernel is spanned by
\vspace{-3mm}
\[
	\sum_{i=0}^kq_i\mathcal{T}_i,
\]
\vskip -0.2cm
\noi where the $\mathcal{T}_i$ are tangle diagrams that differ only in a disc, in which $\mathcal{T}_i$ is isotopic to $\mathcal{Y}^i$. 

We perform 
the following regular isotopy to the $\mathcal{T}_i$: first we shrink the disc until it is very close to the flagpole. Specifically the 
diameter of the disk should be smaller than the distance between the flagpole and any crossing or point with vertical tangent in 
$\mathcal{T}_i$ (not in the disc). Now we ``drag" the disc up the flagpole, passing inside any strands of $\mathcal{T}_i$ that wrap around 
the flagpole, until the disc is above the rest of the diagram. The following diagrams illustrate what happens when we encounter a strand 
crossing the flagpole:
  \begin{center}
    \includegraphics[scale=0.8]{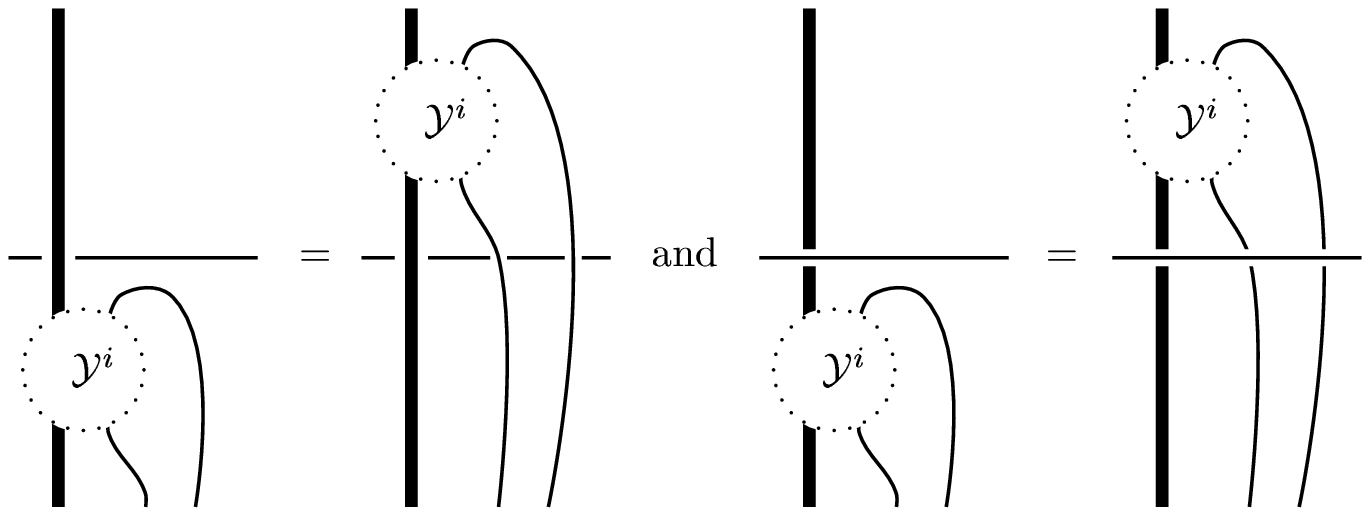}
  \end{center} 
This isotopy turns $\mathcal{T}_i$ into $\cap\mathcal{Y}^i\mathcal{T}'$, where $\mathcal{T}'\in\akt{2,0}(R)$. Note that the 
isotopy did not depend on the contents of the disc, so $\mathcal{T}'$ is independent of $i$. The kernel is therefore spanned by
\[
	\sum_{i=0}^kq_i\cap\mathcal{Y}^i\mathcal{T}',
\]
where $\mathcal{T}'\in\akt{2,0}(R)$. Now by Lemma \ref{lemma:akt20}, the kernel is spanned by
\[
	\sum_{i=0}^kq_i\cap\mathcal{Y}^i\mathcal{Y}^j\cup=\sum_{i=0}^kq_iA_{i+j}.
\]
But we have assumed that the parameters are weakly admissible, so the $A_i$ defined by (\ref{weakA}) are the same as those defined by 
(\ref{extendA}). Therefore the RHS is zero, hence the kernel is zero, completing the proof.
\end{proof}

Under the assumption of weak admissibility, we can therefore compose our earlier map $\varepsilon_{_{\mathbb{T}}}: \ckt{n}(R) \rightarrow \ckt{0}(R)$ with 
$\ckt{0}(R)\stackrel{{}_\sim}{\rightarrow}R$ to obtain a trace map
\[
\varepsilon_{_{\mathbb{T}}}: \ckt{n}(R) \rightarrow R.
\]
%This remark only works if A_0 is invertible... Sorry! But we get this inclusion later anyway :)
%\textbf{Remark:} Observe that $\varepsilon_n \circ \iota$ is the identity on $\ckt{n-1}$, hence the map $\iota: \ckt{n-1} \rightarrow \ckt{n}$ is injective and we may naturally regard $\ckt{n-1}$ as a subalgebra of $\ckt{n}$.
%However, it is \emph{not} clear, \emph{a priori}, that the analogous map $\iota: \bmk \rightarrow \bnk$ defined by $Y \mapsto Y$, $X_i \mapsto X_i$ and $e_i \mapsto e_i$ is injective.

\textbf{Remark:} The map $\varepsilon_n$ is sometimes called a ``\emph{conditional expectation}''. Also, for a ring $S$ with admissible or weakly admissible parameters, the trace map $\varepsilon_{_{\mathbb{T}}}: \ckt{n}(S) \rightarrow S$, up to multiplication by a power of $A_0$, satisfies the properties of a \emph{Markov trace}. This terminology originated in Jones \cite{J83}.

Taking the closure of the usual braids of type $A$ or tangles produces links in $S^3$. In the type $B$ case, closure of affine braids and affine tangles yields links in a solid torus. This leads to the study of Markov traces on the Artin braid group of type $B$ and, moreover, invariants of links in the solid torus. Various invariants of links in the solid torus, analogous to the Jones and HOMFLY-PT invariants (for links in $S^3$), have been discovered using Markov traces on the cyclotomic Hecke algebras; for example, see Turaev \cite{T88}, tom Dieck \cite{T93}, Lambropoulou \cite{L99} and references therein.
Kauffman-type invariants for links in the solid torus can be recovered from Markov traces on the affine BMW algebras, and similarly, the cyclotomic BMW algebras. This is discussed by Goodman and Hauschild in~\cite{GH06}. 
\vskip 0.2cm
The trace map on $\ckt{n}$ commutes with specialisation of ground rings, in the following sense.
Suppose we have two rings $S_1$ and $S_2$ with weakly admissible parameters and there is a parameter preserving ring homomorphism $\nu: S_1 \rightarrow S_2$. Then the following diagram of $S_1$-linear maps commutes.
\begin{equation} \label{cd:ckttracespec}
\begin{CD}
\ckt{n}(S_1) @>\varepsilon_{_{\mathbb{T}}}>> S_1 \\
@V{\eta_t}VV  @VV{\nu}V \\
\ckt{n}(S_2) @>\varepsilon_{_{\mathbb{T}}}>> S_2
\end{CD}
\end{equation}
Indeed, because $\ckt{n}(S_1)$ is spanned over $S_1$ by affine $n$-tangle diagrams, it suffices to check $\nu \circ \varepsilon_{_{\mathbb{T}}} = \varepsilon_{_{\mathbb{T}}} \circ \eta_t$ on affine $n$-tangle diagrams, by the $S_1$-linearity of $\varepsilon_{_{\mathbb{T}}}$ and definition of $\eta_t$. This then follows because $\cl{n}$ and the isomorphism $\ckt{0}(S) \cong S$ commutes with specialisation. The former is easy to verify, as it suffices to check on diagrams, and the latter is clear since the isomorphism $\ckt{0}(S) \rightarrow S$ is inverse to the natural inclusion map $S \rightarrow \ckt{0}$. 

Using this trace on the cyclotomic Kauffman tangle algebras, we are now able to define a trace on the cyclotomic BMW algebras, over any ring $S$ with weakly admissible parameters, by taking its composition with the diagram homomorphism $\dg: \bnk(S) \rightarrow \ckt{n}(S)$ described on page \pageref{defn:dg}.

\begin{defn}
  For any ring $S$ with weakly admissible parameters, define the $S$-linear map 
\[\varepsilon_{\!_{\mathscr{B}}} := \varepsilon_{_{\mathbb{T}}} \circ \dg: \bnk(S) \rightarrow S.
\]
\end{defn}
%Note that, in particular, we now have a trace map on $\bnk$ over any ring with \emph{admissible} parameters, by Proposition \ref{prop:WYimpliesGH}.

Now, by the commutative diagrams (\ref{cd:dgspec}) and (\ref{cd:ckttracespec}), it is easy to see that $\varepsilon_{\!_{\mathscr{B}}}$ commutes with specialisation of rings as well. In other words, if $S_1$ and $S_2$ are two rings with weakly admissible parameters and $\nu: S_1 \rightarrow S_2$ is a parameter preserving ring homomorphism, then the following diagram of $S_1$-linear maps commutes.
\begin{equation} \label{cd:bnktracespec}
\begin{CD}
\bnk(S_1) @>\varepsilon_{\!_{\mathscr{B}}}>> S_1 \\
@V{\eta_b}VV  @VV{\nu}V \\
\bnk(S_2) @>\varepsilon_{\!_{\mathscr{B}}}>> S_2
\end{CD}
\end{equation}
In the next section, we will show that for a particular ring with admissible (and hence weakly admissible) parameters, this trace map $\varepsilon_{\!_{\mathscr{B}}}$ is in fact nondegenerate. For this, we consider its relationship with a known nondegenerate trace on the cyclotomic Brauer algebras.

%%%%%%%%%% CYCLOTOMIC BRAUER ALGEBRA %%%%%%%%%%

\section{Cyclotomic Brauer Algebras} \label{section:cyclobrauer}

For a fixed $n$, consider the set of all partitions of the set $\{1,2,\o,n,1',2',\o,n'\}$ into subsets of size two. Any such partition can be represented by a \emph{Brauer $n$-diagram}; that is, a graph on $2n$ vertices with the $n$ top vertices marked by $\{1,2,\o,n\}$ and the bottom vertices $\{1',2',\o,n'\}$ and a strand connecting vertices $i$ and $j$ if they are in the same subset.
Given an arbitrary unital commutative ring $U$ and an element $A_0 \in U$, the \emph{Brauer algebra} $\EuScript{B}_n$ is defined to be the $U$-algebra with $U$-basis the set of Brauer $n$-diagrams. The multiplication rule in the algebra is defined as follows. Given two Brauer $n$-diagrams $D_1$ and $D_2$, let us define $D_3$ to be the Brauer $n$-diagram obtained by removing all closed loops formed in the concatenation of $D_1$ and $D_2$. Then the product of $D_1$ and $D_2$ is defined to be $A_0^r D_3$, where $r$ denotes the number of loops removed.

The Brauer algebras and their representation theory have been studied extensively in the literature. For example, the generic structure of the algebra and a criterion for semisimplicity of $\EuScript{B}_n$ have been determined; see Wenzl \cite{W88}, Rui \cite{R05}, Enyang \cite{E07} and references therein.

The Brauer algebra is the ``classical limit'' of the BMW algebra in the sense that $\EuScript{B}_n$ is a specialisation of the BMW algebra $\bmw_n$ obtained by sending the parameter $q$ to~$1$. Under this specialisation, $\delta$ is sent to $0$ and so the element $X_i$ is identified with its inverse; this is equivalent to removing the notion of over and under-crossings in $n$-tangles.
%The classical limit of a tangle algebra is essentially obtained by taking the specialisation in which one forgets the notion of over and under-crossings.
In a similar fashion, the cyclotomic Brauer algebras may be thought of as the ``classical limit'' of the cyclotomic Kauffman tangle algebras $\ckt{n}$. 

\begin{defn} (cf. \cite{PS02,GH108})
  A \emph{\textbf{$k$-cyclotomic Brauer $n$-diagram}} (or $\Z_k$-Brauer $n$-diagram) is a Brauer $n$-diagram, in which each strand is endowed with an orientation and labelled by an element of the cyclic group $\Z_k = \Z/k \Z$. Two diagrams are considered the same if the orientation of a strand is reversed and the $\Z_k$-label on the strand is replaced by its inverse (in the group~$\Z_k$).
\end{defn}

\noi An example of a $\Z_6$-Brauer $5$-diagram is given below. 
  \begin{center}
    \includegraphics[scale=0.9]{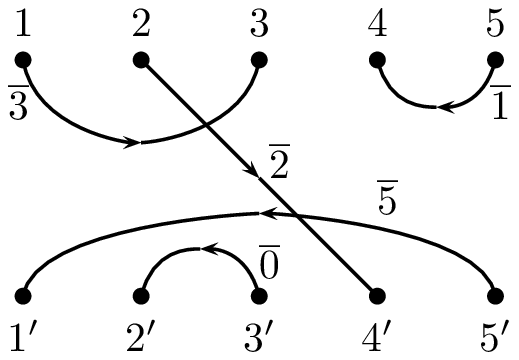}
  \end{center}
Now let $R_c$ denote the polynomial ring $\mathbb{Z}[A_0^{\pm1}, A_1,\ldots, A_{\floor{k/2}}]$.
The following rules define a multiplication for these diagrams.
Firstly, given two $\Z_k$-Brauer $n$-diagrams $D_1$ and $D_2$, concatenate them as one would for ordinary Brauer $n$-diagrams. In the resulting diagram, horizontal strands, vertical strands and closed loops are formed. For each composite strand $s$, we arbitrarily assign an orientation to $s$ and make the orientations of the components of $s$ from the two diagrams agree with the orientation of $s$ by changing the $\Z_k$-labels for each component of $s$ accordingly. The label of $s$ is then the sum of its consisting component labels. Finally, for $d = 0, 1, \ldots, \floor{\frac{k}{2}}$, let $r_d$ be the number of closed loops with label $\pm d$, mod $k$. Let $D_1 \circ D_2$ denote the $\Z_k$-Brauer $n$-diagram obtained by removing all closed loops and define 
\[D_1 \cdot D_2 := \left( \prod_d A_d^{r_d} \right) D_1 \circ D_2.\]
\begin{defn}
  The \emph{\textbf{cyclotomic Brauer algebra}} (or $\mathbb{Z}_k$-Brauer algebra) $\EuScript{CB}_n^k := \EuScript{CB}_n^k(R_c)$ is the unital associative $R_c$-algebra with $R_c$-basis the set of $\Z_k$-Brauer diagrams, with multiplication defined by $\cdot$ above.
\end{defn}

We now proceed to show that $\varepsilon_{\!_{\mathscr{B}}}\!: \bnk(R_c) \rightarrow R_c$ is a nondegenerate trace, using an analogously defined trace on $\EuScript{CB}_n^k(R_c)$. 
%Rather then considering this trace directly in \cite{GH107}, Goodman and Hauschild Mosley instead consider the trace $\varepsilon_{_{\mathbb{T}}}$ using a so-called `connector' map from the cyclotomic Kauffman tangle algebras to the cyclotomic Brauer algebras.

As in the context of tangle algebras, there is similarly a closure map $\mathrm{cl}_n$ from $\Z_k$-Brauer $n$-diagrams to $\Z_k$-Brauer $(n-1)$-diagrams given by joining up the vertices $n$ and $n'$. In addition, any concatenated strands formed in the resulting diagram are labelled according to the same rule as for multiplication in the algebra. Also, if $n$ and $n'$ are joined in the original diagram, then closure will result in a closed loop with some label $\bar{d} \in \Z_k$, where $d = 0,1,\o,\floor{\frac{k}{2}}$, which is then removed and replaced by the coefficient $A_d$. We  extend $\cl{n}$ to a linear map $\varepsilon_n: \EuScript{CB}_n^k(R_c) \rightarrow \EuScript{CB}_{n-1}^k(R_c)$.  Then 
$\varepsilon_c := \varepsilon_1 \circ \cdots \circ \varepsilon_n$ is a trace map
\[
\varepsilon_c: \EuScript{CB}_n^k(R_c) \rightarrow \EuScript{CB}_0^k(R_c) = R_c.
\]
\noi The following lemma is due to the work of Parvathi and Savithri \cite{PS02} on $G$-Brauer algebra traces, for finite abelian groups $G$. 
%A sketch of its proof is given in Goodman and Hauschild Mosley \cite{GH107}.
\begin{lemma} \label{lemma:cycloBrauernondeg}
The trace $\varepsilon_c: \EuScript{CB}_n^k(R_c) \rightarrow R_c$ is nondegenerate. That is, for every $d_1 \in \EuScript{CB}_n^k(R_c)$, there exists a $d_2 \in \EuScript{CB}_n^k(R_c)$ such that $\varepsilon_c(d_1d_2) \neq 0$. Equivalently, it says that the determinant of the matrix $\left( \varepsilon_c(D\!\cdot \!D')\right)_{D,D'}$, where $D,D'$ vary over all $\Z_k$-Brauer $n$-diagrams, is nonzero in~$R_c$.
\end{lemma}

\section{The Freeness of $\bnk$} \label{sect:basis}
Let us fix $\sigma$ to be $+$ or $-$.
In the ring $R_c$, put $q := 1$, $\l := \pm 1$, depending on the sign of~$\sigma$, $q_0 := 1$ and $q_i : =0$, for all $i = 1, \o, k-1$. Also, let $A_j$, where $j \notin \{0, \o, \floor{\frac{k}{2}}\}$, be such that $A_m = A_{m+k}$ and $A_{-m} = A_m$ hold for all $m \in \mathbb{Z}$.
We have a homomorphism from the polynomial ring $\Omega$, defined in Lemma \ref{lemma:Rsigma}, to $R_c$, defined on the generators by
\begin{eqnarray*}
\varsigma: \Omega &\rightarrow& R_c \\ 
\l &\mapsto& \pm 1 \\
q &\mapsto& 1 \quad (\Rightarrow \d \mapsto 0) \\
q_0 &\mapsto& 1 \\
q_i &\mapsto& 0 \\
A_j &\mapsto& A_j.
\end{eqnarray*}
% Certainly, $\varsigma$ defines a map $\Omega \rightarrow R_c$. 

It is easy to verify that the image of $\beta_\sigma$ and $h_l$ under $\varsigma$ are zero, for all $l = 0,1, \o, z-~\epsilon$. Also, by (\ref{hiprime2}), the $h_l'$ are mapped to $A_{k-l} - A_l$ in $R_c$, but these are simply all zero, as $A_l = A_{-l} = A_{-l+k}$ in $R_c$. Thus the generators of $I_\sigma$ 
%indeed 
vanish under the map $\varsigma$, hence the above defines a ring homomorphism $\varsigma_\sigma: R_{\sigma} \rightarrow R_c$. As an 
immediate consequence of this, we also have a map $R_0 \rightarrow R_c$, which factors through $R_{\sigma}$. Observe that, by 
Definition \ref{defn:adm}, the existence of this map  shows that $R_c$ is admissible, and in particular weakly admissible.
% is equivalent to the admissibility of the parameters in $R_c$ chosen at the beginning of this section.
% Hence, by Proposition \ref{prop:WYimpliesGH}, this gives an alternative proof that $R_c$ is a ring with weakly admissible parameters.

\begin{figure}[h!]
  \begin{center}
    \includegraphics[scale=0.9]{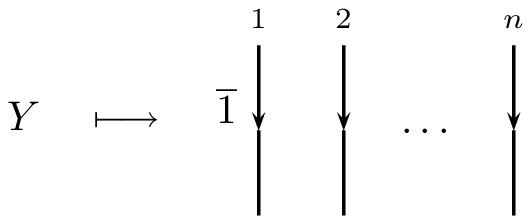} \vskip 0.7cm
    \includegraphics[scale=0.9]{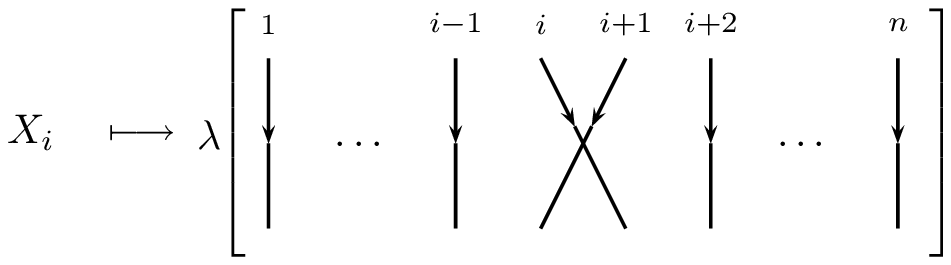} \vskip 0.7cm
    \includegraphics[scale=0.9]{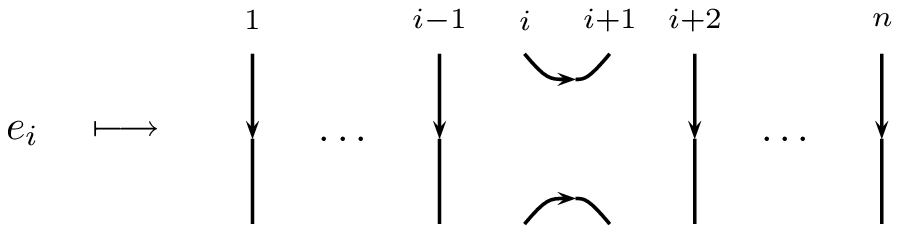}
    \caption{The isomorphism $\xi: \bnk(R_c) \rightarrow \EuScript{CB}_n^k(R_c)$.}
\label{fig:bnk2brauer}
  \end{center}
\vspace{-0.5cm}
\end{figure}

We have an $R_c$-algebra homomorphism $\xi: \bnk(R_c) \rightarrow \EuScript{CB}_n^k(R_c)$, given by Figure \ref{fig:bnk2brauer}, in which only non-zero labels on strands have been indicated.
It is known that the cyclotomic Brauer algebra $\EuScript{CB}_n^k$ is generated by the diagrams given in Figure \ref{fig:bnk2brauer}. (We refer the reader to Rui and Xu \cite{RX07} for a full presentation). Hence $\xi$ is surjective.
Moreover, we already have a spanning set of size $k^n (2n-1)!!$ of $\bnk(R_c)$, given by Theorem \ref{thm:span}. Since $\xi$ is surjective, this maps onto a spanning set of $\EuScript{CB}_n^k$. But $\EuScript{CB}_n^k$ is of rank $k^n (2n-1)!!$, by definition, hence the image of our spanning set of $\bnk(R_c)$ is in fact a $R_c$-\emph{basis} of $\EuScript{CB}_n^k$. Thus, as $R_c$-algebras, $\bnk(R_c) \cong \EuScript{CB}_n^k$, under $\xi$.
We now compile the above information into the following diagram.
\begin{equation} \label{cd:bmwbrauer}
\begin{CD}
\bnk(R_c) @>{\xi}>{\cong}> \EuScript{CB}_n^k(R_c) \\
@V{\dg}VV  @VV{\varepsilon_c}V \\
\ckt{n}(R_c) @>\varepsilon_{_{\mathbb{T}}}>> R_c
\end{CD}
\end{equation}
Let us now prove that the diagram above commutes. Take a diagram $D \in \EuScript{CB}_n^k$. Then the map $\varepsilon_{_{\mathbb{T}}} \circ \dg \circ \xi^{-1} (D)$ is essentially the trace of the diagram obtained by `separating' all the non-zero labels from the rest of the diagram and replacing them with appropriate analogous $\mathcal{Y}$-type diagrams. In $\mathbb{KT}_0^k(R_c)$, over and under-crossings do not matter, so $\varepsilon_{_{\mathbb{T}}} \circ \dg \circ \xi^{-1} (D)$ is reduced to disjoint closed loops (around the flagpole), which are then all identified (as $\mathbb{KT}_0^k(R_c) \cong R_c$) with a product of $A_i$'s, where $i$ depends on the original labels in $D$. This produces precisely the same result as taking the trace $\varepsilon_c$ of $D$, as required.

%We now construct a map $c_k: \ckt{n}(R_c) \rightarrow \EuScript{CB}_n^k(R_c)$. Let T be an affine tangle diagram with OH NO.

We may now use the nondegeneracy of $\varepsilon_c$ and the specialisation maps $R_0\rightarrow R_\sigma\rightarrow R_c$ to deduce that our 
trace map $\varepsilon_{\!_{\mathscr{B}}}\!: \bnk(R_0) \rightarrow R_0$ is nondegenerate over the generic ring, and that the spanning set 
constructed previously is a basis.

\begin{thm} \label{thm:thebigbasisone} 
Let $\mathbb{B}_{R_0}$ denote the spanning set of $\bnk(R_0)$ given in Theorem \ref{thm:span}.
\begin{enumerate}
\item[(1)] $\det\left(\varepsilon_{\!_{\mathscr{B}}}(x_1x_2)\right)_{x_1,x_2\in\mathbb{B}_{R_0}}$ is not a zero divisor in $R_0$.
\item[(2)] $\mathbb{B}_{R_0}$ is a basis for $\bnk(R_0)$. Thus $\bnk(R_0)$ is $R_0$-free of rank $k^n (2n-1)!!$.
\item[(3)] $\dg: \bnk(R_0) \rightarrow \ckt{n}(R_0)$ is an $R_0$-algebra isomorphism. 
\end{enumerate}
\end{thm}
\begin{proof}
Let $\omega := \det\left(\varepsilon_{\!_{\mathscr{B}}}(x_1x_2)\right)_{x_1,x_2\in\mathbb{B}_{R_0}}\in R_0$. Recall the specialisation 
maps $\varsigma_\sigma:R_\sigma\rightarrow R_c$. Let $\pi_\sigma:R_0\rightarrow R_\sigma$ denote the canonical surjections. Then 
$\varsigma_\sigma \circ \pi_\sigma:R_0\rightarrow R_c$ gives an $R_0$-algebra homomorphism
\[
	\eta_\sigma:\bnk(R_0)\rightarrow\bnk(R_c)
\]
such that $\varsigma_\sigma \circ \pi_\sigma \circ \varepsilon_{\!_{\mathscr{B}}}=\varepsilon_{\!_{\mathscr{B}}} \circ \eta_\sigma$. Moreover $\eta_\sigma$ 
induces an isomorphism of $R_c$-algebras
\[
	R_c\otimes_{R_0}\bnk(R_0)\stackrel{{}_\sim}{\rightarrow}\bnk(R_c),
\]
so in particular, $\{\eta_\sigma(x)\mid x\in\mathbb{B}_{R_0}\}$ spans $\bnk(R_c)$. Since $\xi:\bnk(R_c)\rightarrow\EuScript{CB}_n^k(R_c)$ is 
surjective,
\[
	\{\xi \circ \eta_\sigma(x)\mid x\in\mathbb{B}_{R_0}\}
\]
spans $\EuScript{CB}_n^k(R_c)$. But $|\mathbb{B}_{R_0}|=k^n (2n-1)!!$ is the rank of $\EuScript{CB}_n^k(R_c)$, so this set is a basis. 

\noi In particular, since $\varepsilon_c$ is nondegenerate,
\[
	\det\left(\varepsilon_c(\xi\eta_\sigma(x_1)\xi\eta_\sigma(x_2))\right)_{x_1,x_2\in\mathbb{B}_{R_0}}\neq0
\]
in $R_c$. Hence
\begin{eqnarray*}
	\varsigma_\sigma \circ \pi_\sigma(\omega)
		&=&\varsigma_\sigma \circ \pi_\sigma[\det\left(\varepsilon_{\!_{\mathscr{B}}}(x_1x_2)\right)_{x_1,x_2\in\mathbb{B}_{R_0}}]\\
		&=&\det\left(\varsigma_\sigma\pi_\sigma\varepsilon_{\!_{\mathscr{B}}}(x_1x_2)\right)_{x_1,x_2\in\mathbb{B}_{R_0}}\\
		&=&\det\left(\varepsilon_{\!_{\mathscr{B}}}(\eta_\sigma(x_1x_2))\right)_{x_1,x_2\in\mathbb{B}_{R_0}}\\
		&=&\det\left(\varepsilon_c(\xi\eta_\sigma(x_1x_2))\right)_{x_1,x_2\in\mathbb{B}_{R_0}}\\
		&\neq&0.
\end{eqnarray*}
Thus $\pi_\sigma(\omega)\neq0$. Now suppose $\omega x = 0$, for some $x\in R_0$. Then $\pi_\sigma(\omega)\pi_\sigma(x)=0$, but 
$R_\sigma$ is an integral domain by part (c) of Lemma \ref{lemma:Rsigma}. Hence $\pi_\sigma(x)=0$ for $\sigma=\pm$. Now part (d) of Lemma~\ref{lemma:Rsigma} shows that $x=0$. This proves (1). Since $\mathbb{B}_{R_0}$ spans $\bnk(R_0)$, (2) now follows.

Finally, by definition of $\varepsilon_{\!_{\mathscr{B}}}$,
\[
	\omega=\det\left(\varepsilon_{\!_{\mathscr{B}}}(x_1x_2)\right)_{x_1,x_2\in\mathbb{B}_{R_0}}
		=\det\left(\varepsilon_{_{\mathbb{T}}}(\dg(x_1)\dg(x_2))\right)_{x_1,x_2\in\mathbb{B}_{R_0}}
\]
is not a zero divisor in $R_0$. This implies that the set $\{\dg(x)\mid x\in\mathbb{B}_{R_0}\}$ is linearly independent over~$R_0$, so 
$\psi$ is injective. We already know it is surjective, hence $\psi$ is an isomorphism, proving~(3).
\end{proof}

\noi We may now specialise this result to an arbitrary admissible ring.
\begin{cor}\label{cor:basis}
Let $R$ be a ring with admissible parameters $A_0, \ldots, A_{k-1}, q_0, \ldots, q_{k-1}$, $q$ and $\l$. Let $\mathbb{B}_R$ denote the set of all
\[
\alpha_{i_1j_1,n-1}^{s_1}\ldots\alpha_{i_mj_m,n-2m+1}^{s_m} \chi^{(n-2m)}(\alpha_{g_mh_m,n-2m}^{t_m})^*\ldots(\alpha_{g_1h_1,n-2}^{t_1})^*,
\]
where $m = 1,2,\ldots, \lfloor \frac{n}{2} \rfloor$, $i_1>i_2>\ldots>i_m$, $g_m<g_{m-1}<\ldots<g_1$ and, for each $f = 1,2,\ldots m$, $1 \leq i_f\leq j_f \leq n-2f+1$, $1 \leq g_f\leq h_f \leq n-2f+1$, $s_f, t_f \in \left\{ \floor{\frac{k}{2}} - (k-1), \o, \floor{\frac{k}{2}} \right\}$ and $\chi^{(n-2m)}$ is an element of $\widetilde{\X}_{n-2m,k}$.

Then $\mathbb{B}_R$ is an $R$-basis of $\bnk(R)$, and its image under $\dg: \bnk(R) \rightarrow \ckt{n}(R)$ is an $R$-basis of $\ckt{n}(R)$. 
In particular, $\dg: \bnk(R) \rightarrow \ckt{n}(R)$ is an isomorphism of $R$-algebras.
\end{cor}
\begin{proof}
By Theorem \ref{thm:thebigbasisone}, $\bnk(R)\cong R\otimes_{R_0}\bnk(R_0)$ and $\ckt{n}(R)\cong R\otimes_{R_0}\ckt{n}(R_0)$ are both free $R$-modules 
of rank $k^n(2n-1)!!$. However, $\mathbb{B}_R$ and $\{\dg(x)\mid x\in\mathbb{B}_R\}$ are spanning sets of these algebras of cardinality $k^n(2n-1)!!$. 
They must therefore be bases and, since $\dg$ maps the former basis to the latter basis, $\dg$ is an isomorphism.
\end{proof}

Recall we noted earlier that it is not clear \emph{a priori} that $\bmk(R)$ is a subalgebra of $\bnk(R)$. This now follows as a direct consequence of the isomorphism $\bnk(R) \cong \ckt{n}(R)$ and we may identify $\bl{l}$ with $\blk$. 

\begin{figure}[h!]
 \begin{center}
   \includegraphics[scale=0.8]{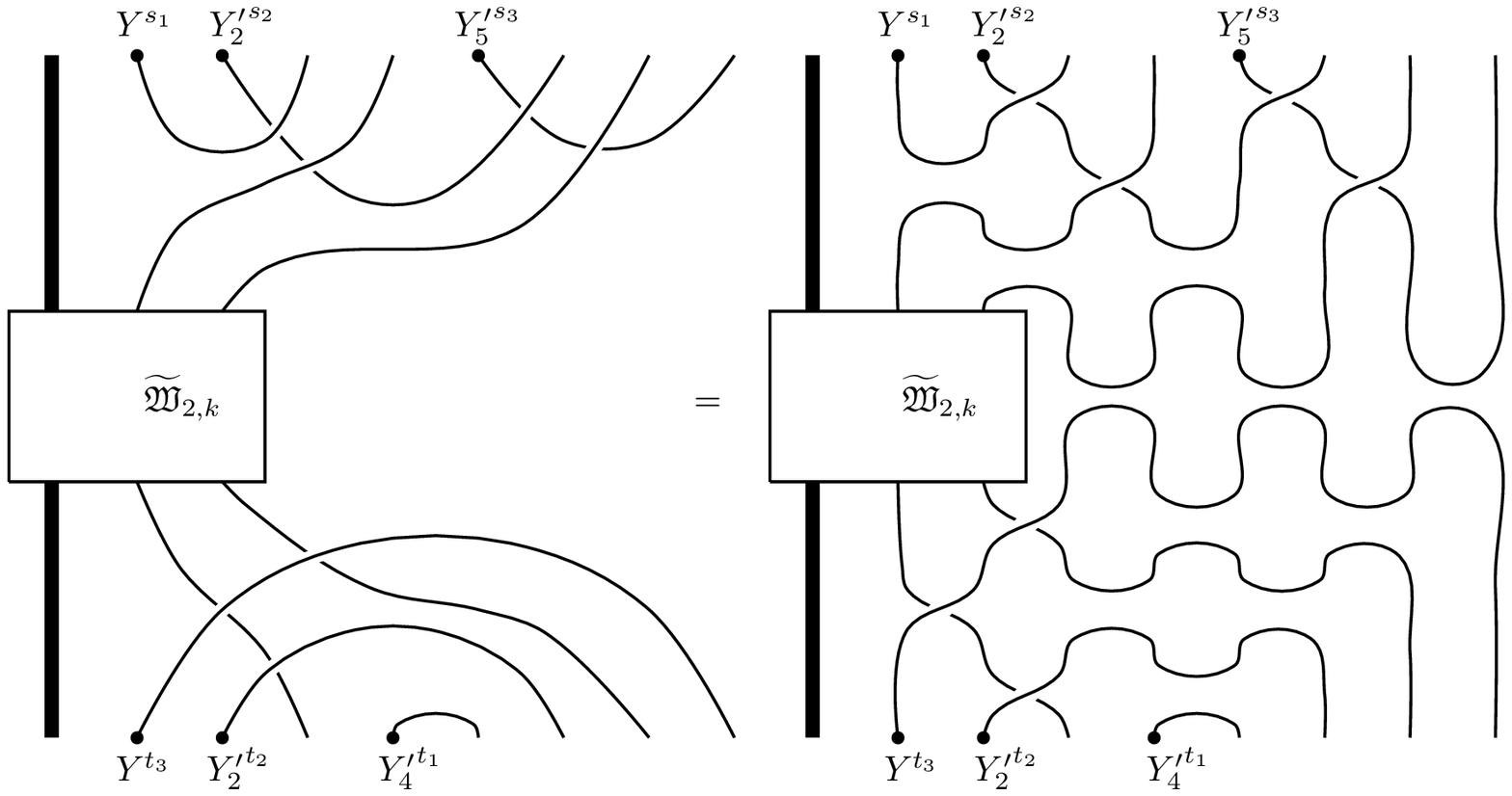}
   \caption{The diagrammatic interpretation of basis element $\a{5,7,7}{s_1}\a{2,4,5}{s_2}\a{1,1,3}{s_3} \,\chi^{(2)}\, (\a{1,3,3}{t_3})^*(\a{2,3,5}{t_2})^*(\a{4,4,7}{t_1})^*$ of $\mathscr{B}_8^k$.} \label{fig:inflatedbasis}
\end{center}
\end{figure}

Finally, to help the reader visualise our basis of the cyclotomic BMW algebras, recall from Figure \ref{fig:onealpha} that an $\alpha$ chain in a basis element corresponds to a horizontal arc in the corresponding tangle diagram. These bases may then be viewed as a kind of ``inflation" of a basis of a smaller Ariki-Koike algebra (which is first lifted to $\widetilde{\X}_{n-2m,k}$ in $\bnk$) by `dangles', as seen in Xi \cite{X00}, with powers of these Jucy-Murphy type elements $Y_i'$ attached. This is illustrated in Figure~\ref{fig:inflatedbasis}, in which the $\y{i}{s}$ elements appearing should of course be replaced with there tangle counterparts (in other words, their images under $\dg$; see Figure \ref{fig:y4dash}).

When $k=1$, $Y$ disappears and no admissibility relations are required on the parameters of $R$ except for $h_0 = 0$. In this case, taking the standard basis of the Iwahori-Hecke algebras of type~$A$ in the basis of Corollary \ref{cor:basis} gives an explicit algebraic description of a known diagrammatic basis of the BMW algebra due to \cite{MW89,HR95,X00}. Also, by the Remarks made after Definition \ref{defn:GHadm} regarding the relationships between the different kinds of admissibility used in the literature, the freeness results of \cite{G08,RX08,RS08} are implied by Corollary \ref{cor:basis}.
Furthermore, with Figure \ref{fig:inflatedbasis} in mind and a carefully constructed lifting map from $\mathfrak{h}_{l,k}$ to $\blk$ (that is, $\overline{\X}_{l,k}$ may not be chosen arbitrarily), it is proven in \cite{WY09b,Y07} that the cyclotomic BMW algebras are cellular, in the sense of Graham and Lehrer~\cite{GL96}. 
%%%%%%%%%%%%%%%%%%%%%%%%%%%%%%%%%%%%%%%%%%%%%%%%%%%
\vskip 0.5cm
\noi \textbf{\large Acknowledgements}

The authors thank Robert Howlett for providing pstricks diagrams appearing in this paper.

\providecommand{\noopsort}[1]{}

\bibliographystyle{elsarticle-harv}
\end{document}